\documentclass[10pt]{article}
\usepackage{amsmath,epsfig,bm,float}
\usepackage{amssymb}
\usepackage{graphicx}
\usepackage{color}
\DeclareMathAlphabet{\itbf}{OML}{cmm}{b}{it}

\newcommand{\sq}{\sqrt}

\newcommand{\RR}{\mathbb{R}}

\newcommand{\ri}{\rightarrow}
\newcommand{\p}{\partial}
\newcommand{\q}{\quad}

\newcommand{\ds}{\displaystyle}

\newcommand{\erfc}{\mbox{erfc}}
\newcommand{\f}{\frac}

\newcommand{\be}{\begin{eqnarray}}
\newcommand{\ee}{\end{eqnarray}}

\newcommand{\bea}{\begin{eqnarray*}}
\newcommand{\eea}{\end{eqnarray*}}

\def\ds{\displaystyle}

\begin{document}

\title{Preferred Frequencies for Coupling of Seismic Waves and Vibrating Tall Buildings}
\author
{Darko Volkov and Sergey Zheltukhin \thanks{Department of Mathematical Sciences, Worcester Polytechnic
Institute,  Worcester MA 01609
({\tt darko@wpi.edu}, {\tt sergey@wpi.edu})} }

\maketitle

\begin{abstract}
We study a model for the so called "city effect"
in which an earthquake can be locally enhanced by the collective
 response of tall buildings in a large city.
We use a set of equations  coupling vibrations in buildings to motion under the ground.
These equations were previously studied exclusively in the case of a finite set 
of identical, equally spaced, buildings.
%We study in this paper two new cases. First, the case of an infinite collection of identical buildings, arranged
%periodically. 
These two restrictions are lifted in this paper. We may now simulate geometries involving 
infinitely many buildings as long as an initial pattern of buildings is repeated.
Our new method using periodic domains and periodic Green's functions 
yields much faster computations. This is the main reason why we are now
able to study systems of buildings
of variable height, mass, and rigidity. 
We show how solving for the wavenumber in a non-linear equation involving the integral
of a function solution to an adequate integral equation, we are able to find 
resonant frequencies  coupling seismic waves and vibrating tall buildings.
Interestingly, in the case of non identical buildings, our simulations indicate that the response 
to this coupling phenomenon
may differ drastically from one building to another.
%We can assume that the number of buildings is infinite as long as a 
%given pattern is repeated periodically. This is the subject of section 5.  

\end{abstract}

\normalsize

%\begin{keywords} 
%Electromagnetic scattering,
%\end{keywords} 

\section{Introduction}
The traditional approach to evaluating seismic risk in urban areas is to consider  seismic waves in the underground 
 as the only cause for  motion above ground. In earlier  studies, seismic wave propagation was evaluated 
in an initial step and in a second step  impacts on  man made structures  were inferred. However, observational 
evidence has since then suggested that when an earthquake strikes a large city,  seismic activity 
may in turn be
 altered by the  response of  the buildings. 
This phenomenon is referred to as the ``city-effect''
and has been studied by many authors, see \cite{millikan, shuttles, ullevi}.\\
 Many occurrences of this city-effect have been documented. 
In 1970, vibrations of the Millikan library on the Caltech campus triggered
 by roof actuators were registered by 
seismographs located a few kilometers away
 (see \cite{millikan}
for a related phenomenon recorded in 2002). 
On multiple instances
 a shock wave formed as the NASA space shuttle  was returning in the atmosphere 
and  as a result many  seismic 
stations near Edwards Air Force Base, California ( see \cite{shuttles}),
recorded tremors. 
In particular, 
a shock wave created by the re-entry of the Columbia space
shuttle in the atmosphere hit high buildings in Los Angeles and induced seismic waves which were recorded by 
seismographs in Pasadena, situated about 15 kilometers away. 
During the  2001 attacks on the Twin Towers in New York City, shortly after
the towers were hit,   tremors were recorded  tens of kilometers away. 
%This phenomenon is called 
%structure-soil interaction. 
It is hypothesized that in these three examples
the transmission of waves to the underground
 occurred because the natural  vibration frequencies of the buildings  occurred to be 
very close to the frequencies of the ground layers in these areas.
%A more complicated phenomenon occurred in Sweden, when the audience of a rock concert at the Ullevi stadium in Gothenburg 
%started to jump in accord with the beat. Resulting waves transmitted to the soil, were trapped in it, and propagated back 
%as the surface waves. These waves in turn, resonating with the crowd jumping, made the stadium shake (see \cite{ullevi}).
%
%The 1985 Michoacan earthquake in Mexico City led Wirgin and Bard, \cite{michoacan},
% to hypothesize that city buildings may collectively affect the ground motion during an earthquake. That idea was supported by several other technical and computational studies, see \cite{site_city_interaction, urbanization_effect, oscillators_on_free}.
Possibly, the most significant observation of the city-effect
happened during the Michoacan earthquake that struck Mexico city in 1985 (see, for example, 
\cite{michoacan_web}). Classical models and 
computational methods failed to explain all the seismic features recorded during that event. This 
 led Wirgin and Bard in \cite{michoacan} to suggest that  ground motion may be significantly altered 
by the presence of
 buildings and this effect may have been enhanced in Mexico City 
since it is such a densely urbanized area. Bard et al. showed experimentally in 
\cite{site_city_interaction} that
this hypothesis certainly has  some relevance and is worth exploring further.\\
 %This is called ``city-effect''.
%One of the crucial common factors for all these observations was that the frequency of the structures coincided with the frequency of the soil layers. Our main objective in this thesis will be studying city frequencies. \\
More recently, \cite{ghergu},
 Ghergu and Ionescu
% proposed a model derived from the equations of physics and a solution algorithm relying on solid mathematics. Our main contribution 
have derived a model for the city effect based on the equations of solid mechanics and appropriate
coupling of the different parts involved in the physical set up of the problem.
They then proposed a clever way to compute a numerical solution to their system of equations.
%
%
%is to extend their work and to provide a mathematical analysis for proving the existence of preferred frequencies coupling vibrations of buildings to underground seismic waves.
%
This way, in \cite{ghergu},
Ghergu and Ionescu were able to compute a city frequency constant:
 given the geometry and the specific physical constants of an idealized two dimensional city, 
they computed a frequency that  leads to the coupling between vibrating buildings and underground seismic waves. 
This is quite an interesting first result, 
but it is limited in scope insofar as
that city frequency constant was obtained by simply increasing the number of buildings at the expense of
 solving larger and larger systems. 
The main point of this study is to show that if 
 instead we use a periodic Green's function and perform computations on a single period,
we can compute coupling frequencies more efficiently. 
This allows for much faster computations, and in turn  
makes it possible to consider more complex geometries within a single period.\\
Here is an outline of this paper. In section 2 we introduce the physics of the problem
under consideration. We restrict this study to anti plane shearing which allows us to use 
scalar displacement fields depending on two spatial variables.
We take into account the mass densities and the shear rigidities of the buildings and of the underground.
Tall buildings are then essentially modeled to be one dimensional with  different values for 
displacements at their foundations and at their tops. These two different displacements are then related in
a balance equation using the elastic moduli of the buildings. We then narrow our focus on time harmonic
solutions, which we non dimensionalize. We obtain in this way a PDE which is just
a Helmholtz equation in a half plane. At the boundary of the half plane a no force
condition is imposed outside the building foundations, while at the foundations we obtain an integral
 differential equation  expressing how the force acting on the buildings relates 
to the effective motion of these buildings. \\
In section 3 we cover the case where the number of buildings is finite and all the buildings are identical. 
This is a shorter section where we essentially recover  results obtained previously by Ghergu and Ionescu, 
\cite{ghergu}. Their method relies on the symmetry of a boundary operator $T$: we prove that
$T$ is symmetric in Appendix A.
In section 4 we cover the case of an infinite collection of identical buildings, arranged
periodically. We examine how our new results relate to Ghergu and Ionecu's results as the number of buildings
 in their case 
grows large. Some technical details regarding our numerical method for computing resonant frequencies
in the infinite, periodic, case are presented in appendix B.\\
Our new method using periodic domains and periodic Green's function 
results in much faster computations. This is the main reason why we are now
able to study the case of buildings
of variable height, mass, and, rigidity. We can assume that the number of buildings is infinite as long as a 
given pattern is repeated periodically. This is the subject of section 5.  

\section{The physics of the city-effect 
problem and non dimentionalization of the associated eigenvalue formulation}
%\section{Cities with a finite number of equal-sized, uniformly spaced, buildings}

\begin{figure}[H]
  \centering
  \includegraphics[width=0.8\textwidth]{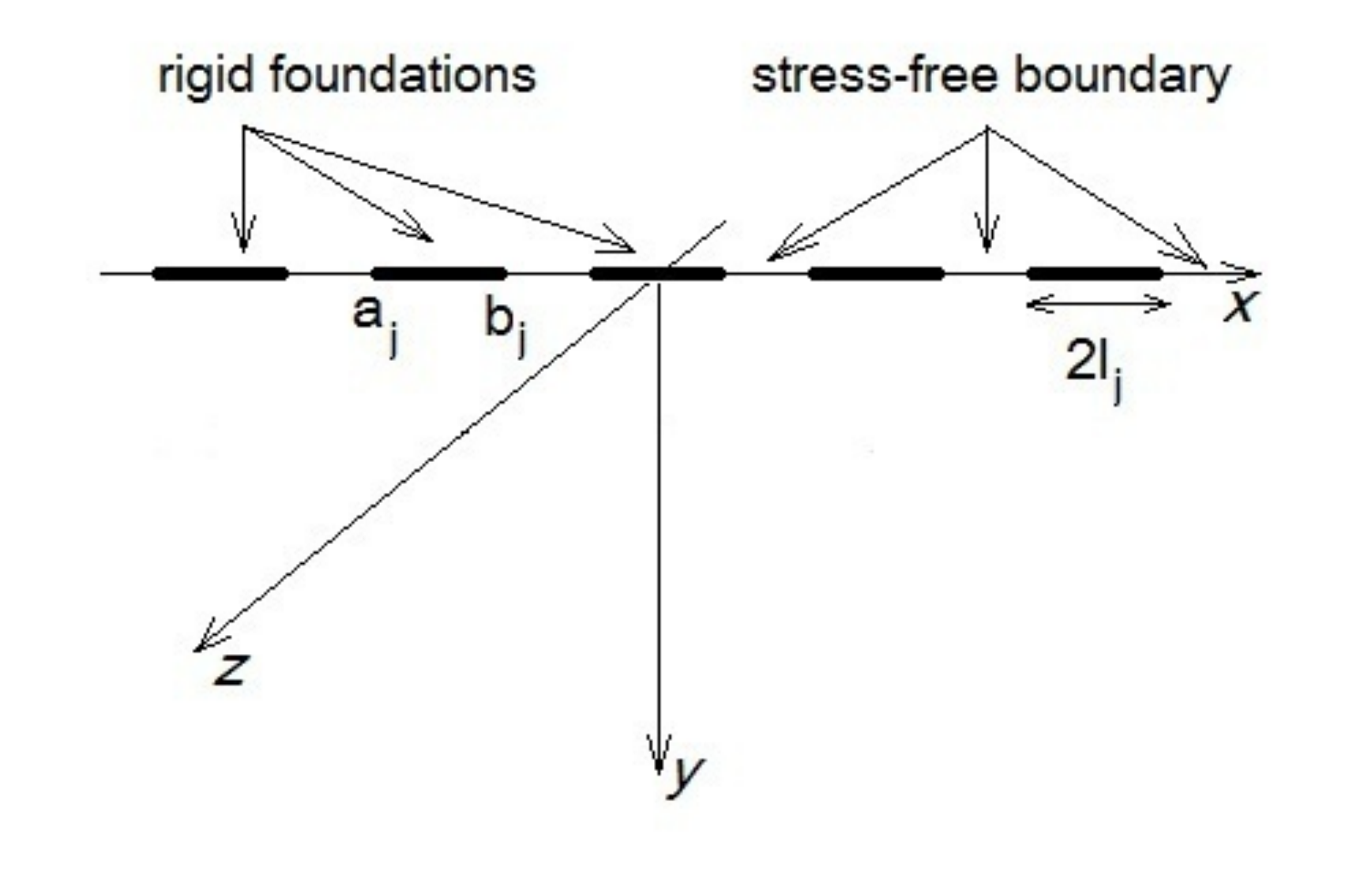}
  \caption{Model geometry: the building foundations lie on the intervals
	$[a_j,b_j]$ which appear in bold.}
  \label{fig:city_layout}
\end{figure}
%----------------------------------------------------------------------------------------------------
\subsection{Physical model}
\label{sec:physical_model}
%We consider soil as 
We model the underground to be 
an elastic half-space 
$y>0$ in the $xyz$ space.
%$\Omega \times \mathbb{R}$, where $\Omega = \mathbb{R} \times (0,\infty)$. 
We only consider the \textbf{anti-plane shearing} case:
all displacements  occur in the $z$ direction 
and are independent of $z$.
We denote by $\Omega$ the half plane $y>0, z=0$.
Consider  $N$ buildings with width $2l_{j}$ and height $h_j$  standing on the
$x$ axis.
%, and the distance between two consecutive buildings is $space$. 
In this study $N$ is either a finite number or it is infinite.
In the first case the index $j$ will range from 1 to $N$, and in the second case 
$j$ will take all values in the set of integers $\mathbb{Z}$.
The rigid building foundations  are all 
located along the $x$-axis and denoted by $\Gamma_j = [a_j,b_j] \times \{0\}$ in $\Omega$. Let us 
introduce the following notations
\begin{itemize}
  \item $\Gamma = \ds \bigcup_{j=1}^{N} \Gamma_j$, if $N$ is finite, or
  $\Gamma = \ds \bigcup_{j=- \infty}^{\infty} \Gamma_j$, if $N$ is infinite:
  the set of building foundations
  \item $\Gamma_{free} = \mathbb{R} \times \{0\} \backslash \Gamma$: the stress free soil boundary
  \item $w(t,x,y)$: the scalar displacement field
\end{itemize}
The other physical parameters relevant to our problem are 
\begin{itemize}
	\item $\rho$, $\rho_j$: the  mass densities of the underground and of building $j$
	\item $S$, $S_j$: the   shear rigidities of the underground and of building $j$
	\item $\ds \beta = \sqrt{S/\rho}$, $\ds \beta_j = \sqrt{S_j/\rho_j}$: the  shear velocities
	 of the underground and of building $j$
	\item $k_j$: the elastic modulus of  building $j$
	\item $u_j(t)$: the displacement of the rigid building foundation $\Gamma_j$
	\item $v_j(t)$: the displacement of the top of the  building $j$
	\item $m_{1,j}, \: m_{0,j}$: the mass at the top and at the foundation of building  $j$
	\item $R_j(w)$: the underground force acting on the building foundation $\Gamma_j$
\end{itemize}
\begin{figure}[H]
  \centering
  \includegraphics[width=0.8\textwidth]{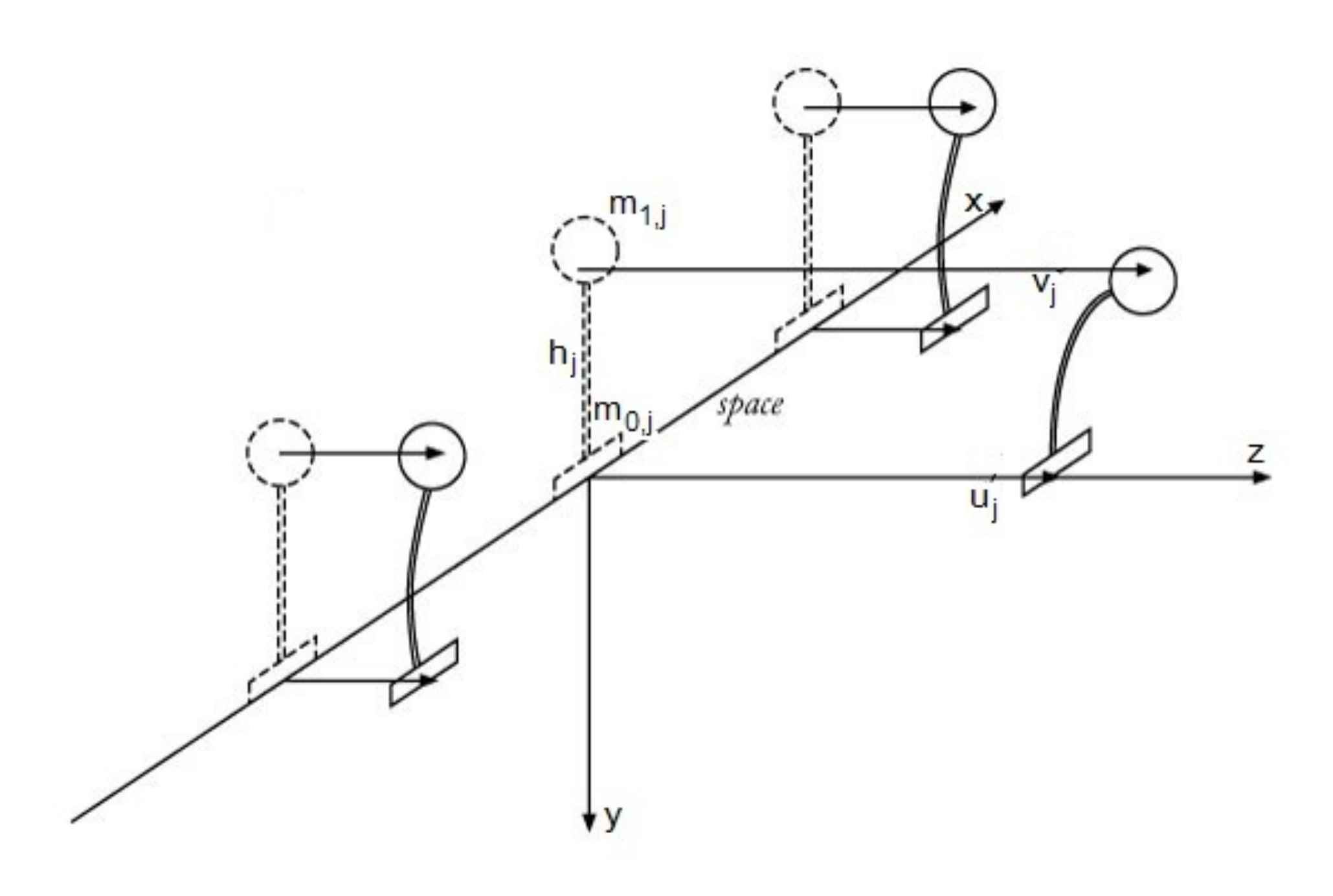}
  \caption{A sketch of  city buildings subject to anti plane shearing. Here $N=3$.
  Building foundations are represented by rectangles and tops are represented by circles.
  Initial positions are sketched in dashed lines, displaced positions are sketched in solid lines.}
  \label{fig:sample_city}
\end{figure}
Using fundamental laws of solid physics, we obtain the following time dependent equations, where $t$ is the
 time variable,
\begin{align}
  \label{wave_equation}
  \rho \ddot{w}(t) &= S \Delta w(t) \mbox{ in } \Omega \times \RR, \\ 
  \label{boundary_displacement}
  \displaystyle w(t,x,0) &= u_j(t)  \mbox{ for all } (x,0)\in \Gamma_j \mbox{, } \frac{\partial w}{\partial y} (t,x,0) = 0 \; \mbox{ for all } (x,0) \mbox{ in } \Gamma_{free}, \\ 
  \label{newton_law_top}
  m_{1,j} \ddot{v}_j(t) &= -k_j(v_j(t)-u_j(t)), \\
  \label{force} 
  \displaystyle R_j(w) &= \int_{\Gamma_j} {S \frac{\partial w}{\partial y}(t,s,0)} ds, \\
  \label{newton_law_found}
  \displaystyle m_{0,j} \ddot{u}_j(t) &= R_j(w) + k_j(v_j(t)-u_j(t)).
\end{align}
Note that  (\ref{wave_equation}) is the wave equation for the displacement $w$; (\ref{boundary_displacement}) shows that the displacement is constant for each rigid foundation and the space between the buildings is stress free; (\ref{newton_law_top}) 
comes from Newton's law of motion for the top of building $j$;
equation  (\ref{force}) expresses continuity of forces between
the underground and our one dimensional building foundations $\Gamma_j$;
and (\ref{newton_law_found}) comes from Newton's law of motion for the bottom of building $j$.
%and integral (\ref{force}) expresses continuity of forces between
%the underground and our one dimensional building foundations $\Gamma_j$.
%----------------------------------------------------------------------------------------------------
\subsection{The associated spectral problem}
We focus on   time harmonic solutions to the system (\ref{wave_equation}-\ref{newton_law_found}). 
Accordingly, we set
\be
\label{time_harmonic_solution}
 w(t,x,y) = \mbox{Re }(\Phi (x,y) e^{-i\omega t}).
\ee 
In other words  $\omega > 0$ is the associated frequency
of a time harmonic vibration. Let us denote by $\alpha_j, \: \eta_j $  the displacements of the foundation and the top of each building. 
%After performing the necessary differentiation in (\ref{wave_equation})-(\ref{newton_law_found}), we obtain the corresponding eigenvalue problem:
After substituting (\ref{time_harmonic_solution}) in (\ref{wave_equation}-\ref{newton_law_found}) we obtain
\begin{align}
    -S \Delta \Phi &= \rho \omega^2 \Phi \mbox{ in } \Omega, \label{dim1}
  \\ k_j(\eta_j - \alpha_j) = \omega^2 m_{1,j} \eta_j, & \; \; -R_j(\Phi) - k_j(\eta_j - \alpha_j) 
  = \omega^2 m_{0,j} \alpha_j, \label{before_non_dim_2}
  \\ \displaystyle \Phi = \alpha_j \mbox{ on } \Gamma_j, & \; \; \frac{\partial \Phi} {\partial y} = 0 \mbox{ on } \Gamma_{free}. \label{dim3}
\end{align}

We now move on  to non dimentionalize equations (\ref{dim1}-\ref{dim3}).  
We proceed exactly as in \cite{ghergu}. We carefully cross checked all the corresponding 
calculations and our results agree with those from \cite{ghergu}, save for the second identity
in equation (\ref{before_non_dim_2}) of our present paper. 
This discrepancy is most likely due to a typographical error in \cite{ghergu}.
%The last step is the non-dimensionalization of the problem above. 
We introduce a characteristic length $l$. The non-dimensional spatial coordinates are: 
\begin{equation}
  \label{non_dimensional_coordinates}
  \displaystyle x^{'} = \frac{x}{l}, \; y^{'} = \frac{y}{l}
\end{equation}
Accordingly, the non dimensional frequency comes out as
\be
\xi = \omega \frac{l}{\beta}
\ee
From now on  we will omit primes and write $x$ and $y$ 
in place of $x'$ and $y'$
%but these are the \textbf{non-dimensional} coordinates. Sets $\Omega, \; \Gamma, \; \Gamma_{free}$ change accordingly, but we will keep the notation. 
for ease of notation.
For each building $j$ we
introduce the non dimensional  parameters 
\begin{equation} 
  \label{non_dimensional_parameters}
  \displaystyle \gamma_j = \frac{m_{1,j}}{m_{0,j}}, \; f_j = \frac{l_j}{h_j}, \; c_j = \frac{l_j}{l}, \; r_j = \frac{\rho_j}{\rho}, \; {\cal{B}}_j = \frac{\beta_j}{\beta}.
\end{equation}
Note that $m_{1,j}, l_j, h_j, \rho_j$ are related
by the equation $$m_{1,j} = 2 l_j h_j \rho_j,$$
and $k_j $ is related to the shear rigidity $S_j$ through
$$
k_j = \f{2S_j l_j}{h_j}.
$$
After a long calculation we arrive at
\be
\eta_j = - \frac{{\cal{B}}_j^2 f_j^2 \alpha_j} {p_j(\xi^2)}, \mbox{ where }
p_j(\xi^2) = c_j^2 \xi^2 - {\cal{B}}_j^2 f_j^2
\ee
Further calculations lead to the system of equations
\begin{align}
  \label{non_dim_helmholtz}
  \Delta \Phi + \xi^2 \Phi &= 0 \mbox{ in } \Omega,
  \\ \label{non_dim_boundary}
  \displaystyle \frac{\partial \Phi}{\partial y} &= 0 \mbox{ on } \Gamma_{free},
  \\ \label{non_dim_eigenvalue}
  \displaystyle q_j(\xi^2) \Phi(x,0) = p_j(\xi^2) \int_{\Gamma_j} {\frac{\partial \Phi} {\partial y} (s,0)} ds & \mbox{ for } (x,0) \in \Gamma_j, 
\end{align}
where
\begin{equation}
\label{p_and_q}
  %\ds p(\xi^2) = c_b^2 \xi^2 - b^2 f_b^2, \; \; 
  q_j(\xi^2) = \frac{2r_j c_j^2 \xi^2} {f_j} \big( c_j^2 \xi^2 - \frac{\gamma_j+1}{\gamma_j}p_j(\xi^2) \big).
\end{equation}
Note that system (\ref{non_dim_helmholtz}) implies that $\Phi$ is constant
on each $\Gamma_j$. We have denoted these constants by $\alpha_j$: we will keep using this notation throughout this paper.
Note also that equation (\ref{non_dim_eigenvalue}) must hold for all integers $j$ between 1 and $N$,
if there are $N$ buildings, or for all $j$ in $\mathbb{Z}$ if there are infinitely  many buildings.

%we can calculate $\alpha, \; \eta$ as
%\begin{center}
%  $\ds \alpha_j = \Phi(x,0) \mbox{ for (x,0)} \in \Gamma_j, \; \eta_j = - \frac{b^2 f_b^2 \alpha_j} {p(\xi^2)}.$
%\end{center}

%Our main objective is to find values of the non-dimensionalized frequency $\xi$ for which system (\ref{non_dim_helmholtz})-(\ref{non_dim_eigenvalue}) will be solvable. In Chapter III we will describe the mathematical solution of this problem and its numerical implemetation, both given by Ghergu and Ionescu. Then we will modify the model to extend it to a wider variety of city types and to improve the speed and the accuracy of the numerical algorithm.
%----------------------------------------------------------------------------------------------------
%\begin{figure}
%\begin{center}
%\includegraphics[width=0.45\textwidth]{carte_Mexique.pdf}
%\includegraphics[scale=.5]{carte_Mexique.pdf}
%\caption{The Guerrero area of Mexico. The subduction zone studied in this paper
%meets the sea floor of the Pacific ocean along the dashed line. The large triangles mark
%the locations of the GPS stations whose data will be used in subsequent sections.
%They are named ACAP, ACYA, CAYA, COYU, CPDP, 
%DEMA, DOAR, IGUA, MEZC, OAXA, 
%PINO, UNIP, YAIG, ZIHP
%} \label{geography}
%\end{center}
%\end{figure}

\section{The case of finitely many buildings which are all identical and equally spaced}
In this case we assume that 
\begin{itemize}
\item $N$ is finite
\item all the physical parameters of the buildings $m_{1,j}$, $m_{0,j}$, 
$l_j$, $h_j$, $\rho_j$, $\beta_j$ are independent of $j$
\item $b_j - a_j$ (the length of building $j$)and $a_{j+1} - b_j$ (the distance from
building $j$ to building $j+1$) are independent of $j$
\end{itemize}
This was the only case considered in \cite{ghergu}.

\subsection{Using linearity to reduce computational time}
Finding a solution to problem (\ref{non_dim_helmholtz}-\ref{non_dim_eigenvalue})
becomes increasingly difficult as $N$ increases, however, Ghergu et al. were able to
combine integral equation techniques to 
eigenvalues of relevant symmetric matrices  in order to reduce computational time.
Let us now overview their computational method.
Fix $\alpha=(\alpha_1, .., \alpha_N)$ in $\RR^N$ and introduce the partial differential equation.
\begin{align}
  \label{helmholtz_1} 
    \Delta \Psi + \xi^2 \Psi &= 0 \mbox{ in } \Omega, \\ 
  \label{helmholtz_2} 
    \ds \Psi = \alpha_j \mbox{ on } \Gamma_j, & \; \frac{\partial \Psi} {\partial y} = 0 \mbox{ on } \Gamma_{free}, \\ 
  \label{sommerfeld_cond} 
    \ds \frac{\partial \Psi}{\partial r} - i\xi \Psi &= o(r^{-1/2}) \mbox{ as } r=|x| \rightarrow +\infty.
\end{align}
%Clearly, a physical solution to our problem must decay at infinity. We therfore
%add the condition
%\be
%ds \frac{\partial \Phi}{\partial r} - i\xi \Phi &= o(|x|^{-1/2}) \mbox{ as } |x| \rightarrow \infty.
%\label{at infinity}
%\ee
For any fixed $\xi>0$, problem (\ref{helmholtz_1}-\ref{sommerfeld_cond}) is uniquely solvable:
this can be shown using standard PDE techniques, see \cite{DVandSergey_theory}.
We now introduce the $N \times N$ matrix $T(\xi^2)$ whose $k \, l$ entry is 
\be
\label{operator_matrix}
  T(\xi^2)_{k \,l} = Re \int_{\Gamma_k} {\frac{\partial \Psi_{\xi,e_l}} {\partial y} (s,0)} ds,
\ee
where $e_l$ is the $l^{th}$ basis vector in $\RR^N$, 
and $\Psi_{\xi,e_l}$ solves (\ref{helmholtz_1}-\ref{sommerfeld_cond})
for $\alpha = e_l$. \\
Gherghu et al. observed  in \cite{ghergu}  that the matrix $T(\xi^2)$ is symmetric.
We provide a proof in Appendix A.
Since $T(\xi^2)$ is symmetric, it can be diagonalized, and all its eigenvalues are real. For every $\xi^2>0$ we denote the eigenvalues of $T(\xi^2)$ by
\begin{equation}
\label{xi_equation}
  \tau_1(\xi^2) \leq \tau_2(\xi^2) \leq \ldots \leq \tau_N(\xi^2),
\end{equation} 
and corresponding normalized eigenvectors $\theta_1(\xi^2), \ldots, \theta_N(\xi^2)$. \\
Since in the present case 
all the non-dimensional parameters (\ref{non_dimensional_parameters}) are the same for all the buildings, the
functions $p_i$ and $q_i$ do not depend on $i$:
 they will therefore be simply denoted $p$ and $q$ in this section.
Assume that $\xi$ is such that for some $i$
\begin{equation}
  \label{eigenvector_equation}
  p(\xi^2) \tau_i (\xi^2) = q(\xi^2).
\end{equation} 
Consider the 
corresponding
 solution $\ds \Psi_{\xi,\theta_i(\xi^2)}$ to (\ref{helmholtz_1}-\ref{sommerfeld_cond})
for $\alpha=\theta_i(\xi^2)$.
Given that, $T(\xi^2) \theta_i(\xi^2) = \tau_i(\xi^2) \theta_i(\xi^2)$,
the coordinates of $\theta_i(\xi^2)$ are real and $\Psi$ solution to (\ref{helmholtz_1}-\ref{sommerfeld_cond})
is linear in $\alpha$,
it follows that
$$
\mbox{Re }\Psi_{\xi,\theta_i(\xi^2)} = 
\sum_{l=1}^N (\theta_i(\xi^2) \cdot e_l  ) \mbox{Re }\Psi_{\xi,e_l},
$$
and
$$
\mbox{Re } \int_{\Gamma_k} \f{\p \Psi_{\xi,\theta_i(\xi^2)}}{\p y } (s,0) ds
= \tau_i(\xi^2) \theta_i(\xi^2) \cdot e_k,
$$
for all $k$ from 1 to $N$. Now using  relation (\ref{eigenvector_equation})
it follows that 
$$
p(\xi^2)\mbox{Re } \int_{\Gamma_k} \f{\p \Psi_{\xi,\theta_i(\xi^2)}}{\p y } (s,0) ds
= q(\xi^2) \theta_i(\xi^2) \cdot e_k,
$$
for all $k$ from 1 to $N$.
% By linearity, $\ds \Psi_{\xi,\theta_i(\xi^2)}$  also satisfies (\ref{non_dim_eigenvalue}) and 
% the system (\ref{non_dim_helmholtz}-\ref{non_dim_eigenvalue}). 
We have thus found a wavenumber $\xi$
and boundary conditions $\alpha = \theta_i(\xi^2)$ such that 
(\ref{non_dim_helmholtz}-\ref{p_and_q}) hold for $\Phi = \mbox{Re } \Psi_{\xi,\theta_i(\xi^2)}$. \\
To find in practice a $\xi $ satisfying (\ref{eigenvector_equation})
we fix $i$ and we proceed to solve numerically (\ref{eigenvector_equation}) as a non linear equation in $\xi$.
The actual existence of a solution is a rather involved theoretical question. We address this question in a separate
study, \cite{DVandSergey_theory}.

%and solve nonlinear equation (\ref{eigenvector_equation}) for the unknown $\xi$ using $Matlab$ solver. We find $\ds \Psi_{\xi,\theta_i(\xi^2)}$  which satisfies the system (\ref{helmholtz_1})-(\ref{sommerfeld_cond}). Then $\ds \Phi = Re \Psi_{\xi, \theta_i(\xi^2)}$ will solve the non-linear eigenvalue problem (\ref{non_dim_helmholtz})-(\ref{non_dim_eigenvalue}).

\subsection{An integral equation for solving the system (\ref{helmholtz_1}-\ref{sommerfeld_cond})}
Clearly, the crucial step in this symmetric matrix  method is the ability 
to solve  the  Helmholtz problem (\ref{helmholtz_1}-\ref{sommerfeld_cond}).
This is most easily done through an integral equation formulation.
Accordingly, we set the solution to 
(\ref{helmholtz_1}-\ref{sommerfeld_cond}) to be 
\begin{equation}
  \label{single_layer_potential}
  \displaystyle \Psi(x,y) = \int_{\Gamma}  G(x-s, y) \psi(s) ds,
\end{equation}
where $G$ is the usual fundamental solution of the Helmholtz equation: $G(x,y) = \frac{i}{4} H_0^{(1)}(\xi
\sq{x^2 + y^2})$.
%\begin{equation}
%  \label{free_space_fund_sol}
%  \displaystyle G(z) = \frac{i}{4} H_0^{(1)}(z).
%\end{equation}
We explain in  \cite{DVandSergey_theory} why formulation (\ref{single_layer_potential}) is valid.
Using potential theory we can argue that the function $\psi$ involved in (\ref{single_layer_potential}) 
is in fact equal on $\Gamma$ to $\ds \lim_{y \ri 0^+}\ds -2 \f{\p \Psi}{\p y} (x,y)$.
Due to equation (\ref{non_dim_helmholtz}-\ref{non_dim_eigenvalue}) $\psi$ must satisfy
\begin{equation}
  \label{single_layer_potential2}
  \int_{\Gamma_i}  G(x-s, y) \psi(s) ds = \alpha_i, \q i=1, \ldots, N .
\end{equation}
We present in appendix a numerical method for solving integral equation (\ref{single_layer_potential}).
\subsection{Numerical illustration for $N=1$}
\begin{figure}[H]
  \centering
  \includegraphics[scale=.6]{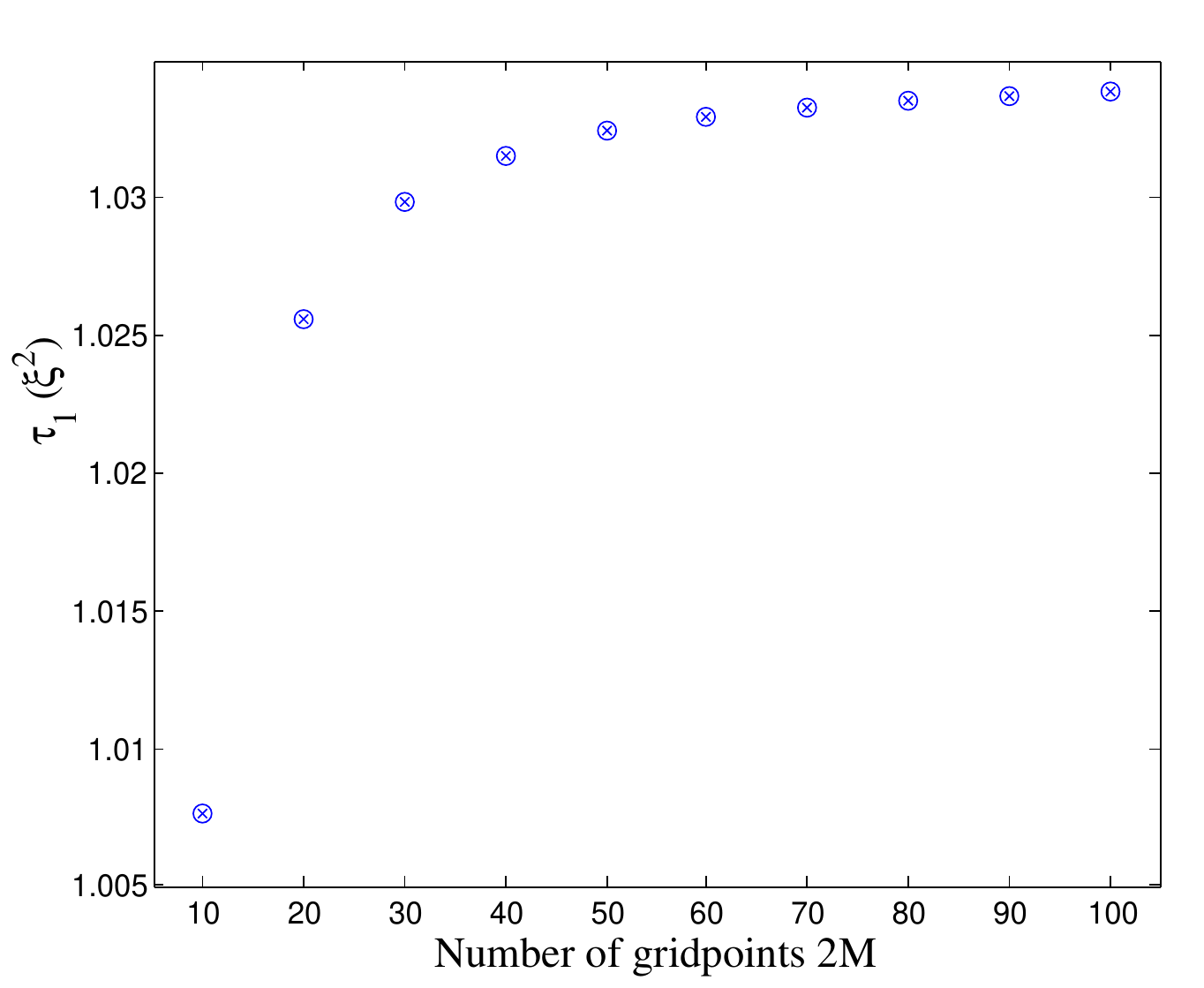}%[width=0.8\textwidth,height=0.4\textheight]
  \caption{One-building city: $l_b = 1$, $b_{1} - a_1 = 0.4$, $l = l_1$, $\gamma_1 = 1.5$, $f_1 = 0.5$, $c_1 = 1$, $r_1=0.1$, ${\cal{B}}_1=1.5$. 
  Observed numerical convergence of the eigenvalue $\tau_1(\xi^2)$  as the number of gridpoints $2M $ increases.
  Note that for $M=5$ the numerical error appears to be under $3 \%$.}
  \label{fig:Fig10_Ghergu}
\end{figure}

We  verify on an example that we can recover the numerical values obtained
by Gherghu et al. in the case of a single building. 
In order to do this, we repeat a calculation from \cite{ghergu}. 
More specifically, we treat the case relative to Figure  10 from \cite{ghergu} 
where the  city and underground parameters were set to $l_1 = 1$, 
$b_{1} - a_1 = 0.4$, $l = l_1$, $\gamma_1 = 1.5$, $f_1 = 0.5$, $c_1 = 1$, $r_1=0.1$, ${\cal{B}}_1=1.5$.
Note that in this case the matrix $T$ reduces to a scalar.
Denote by $2M$ the number of grid points for that building.
For $M=10, 20, ..., 100$, we solve the non linear equation in 
$\xi$ given by (\ref{xi_equation}) where $i=N=1$ and 
$\ds \tau_1(\xi^2) = \mbox{Re } \int_{\Gamma_1} {\frac{\partial \Psi_{\xi,e_1}} {\partial y} (s,0)}$.
The numerical convergence of the solution as $M$ grows large is clearly
observed in 
Figure \ref{fig:Fig10_Ghergu}: on this example we were able to  replicate 
the results given by Figure 10 in \cite{ghergu}. We also note that for $M=5$ the numerical error appears to be under $3 \%$.
% $M=5$, 
%we observed the same numerical convergence of the smallest eigenvalue $\xi_1$ to $0.7792$ as number of buildings $N$ grows, just the same as on Fig. 5 of \cite{ghergu}. 
%For one-building system, the computed eigenvalue $\tau(\xi^2)$ of (\ref{operator_matrix}) exhibits convergent behavior identical to Fig. 10 from \cite{ghergu} as number of gridpoints $2M$ increases.   

\section{The large number of identical buildings case: convergence to periodic structures as $N$ grows large}\label{matrix_method}

\subsection{Ghergu et al.'s large computations and our novel idea on how to drastically cut computational time}
In this case we assume that all the physical parameters of the buildings $m_{1,j}$, $m_{0,j}$, 
$l_j$, $h_j$, $\rho_j$, $\beta_j$ are independent of $j$
and $b_j - a_j$ and $a_{j+1} - b_j$ are independent of $j$, that is the length of foundations is constant, and the buildings
are equally spaced. If $N$ is finite, that case was covered in \cite{ghergu}. 
Ghergu et al. were actually interested in identifying a collective behavior of a large set
of buildings, which they called "city effect".
To do so, they simply applied the computation outlined in the previous section and let the parameter $N$ grow large.
Despite their brilliant idea consisting of using the matrix $T$ defined in the previous section, the computational time 
involved in their method can become prohibitive since it involves \textbf{ for each step of the search 
in $\xi$ } solving a $2 M N$ by $ 2 M N$ linear system for $N$ different right hand sides.  Our novel  idea is to set
directly $N = \infty$ and to introduce a periodic PDE. \\
Let us first introduce notations for our new periodic computational domain, where $2P$
is the period in the $x$ direction 
\be
     \begin{array}{l}
     \Gamma_{per} = [-l_b,l_b] \nonumber \\ 
		\Gamma_{per}^{free} =  (-P, P)\setminus \overline{\Gamma_{per}}  \nonumber
		 \\
     \Omega_{per} = (-P,P) \times (0, \infty)
     \end{array}
     %\Gamma_{per}^{free} &=  [-l_b-\frac{1}{2}space, -l_b] \cup [l_b,l_b+\frac{1}{2}space], \\
     %\Omega_{per} &= [-l_b-\frac{1}{2}space,l_b+\frac{1}{2}space] 
\ee
%In particular the distance between two consecutive buildings is $space$.
%In the case of multiple, equally spaced
%identical buildings 
Accordingly, in place  of (\ref{helmholtz_1}-\ref{sommerfeld_cond})
 we now solve a periodic Helmholtz equation in a  domain
which is bounded in the $x$ direction. We formulate this partial 
differential equation  as follows:
\begin{align}
     \label{helmholtz_per_1} \Delta \Psi + \xi^2 \Psi &= 0 \mbox{ in } \Omega_{per},
  \\ \label{helmholtz_per_2} \displaystyle \Psi = 1 \mbox{ on } \Gamma_{per}, & \; \; \frac{\partial \Psi} {\partial y} = 0 \mbox{ on } \Gamma_{per}^{free}.
\end{align}
augmented by the decay condition 
\be \label{per_decay}
 \frac{\partial \Psi}{\partial y} - i\xi \Psi &= o(y^{-1/2}) \mbox{ as } y \rightarrow +\infty,
\ee
and the periodic boundary condition
\be \label{per_cond}
\Psi(-P, y) = \Psi(P, y)  \mbox{ for all } y > 0
\ee
The non linear equation in $\xi$ to be solved is now reduced to
\be
  \label{eigenvalue_per}
  \ds q(\xi^2) = p(\xi^2) \mbox{Re } \int_{\Gamma_{per}} {\frac{\partial \Psi} {\partial y} (s,0)} ds,
\ee
Evidently this new computational technique relies on the ability
to solve efficiently PDE (\ref{helmholtz_per_1}-\ref{per_cond}).
This can be done by setting
\begin{equation}
  \label{layer_per}
  \displaystyle \Psi(x,y) = \int_{\Gamma_{per}}  G_{per}(x-s, y) \psi(s) ds,
\end{equation}
where $G_{per}$ is the adequate Green's function relative
 to the periodic problem (\ref{helmholtz_per_1}-\ref{per_cond}), and solving for $\psi$ the integral equation
\begin{equation}
  \label{single_layer_potential2_per}
  \int_{\Gamma_{per}}  G_{per}(x-s, y) \psi(s) ds = 1,
\end{equation}
where $x$ and $s$ are in $\Gamma_{per}$.
Defining and computing $G_{per}$ is a vast subject: we discuss it in Appendix B.

%and the matrix (\ref{operator_matrix}) reduces to a scalar. In this case we do not have several eigenvalues and do not have to calculate an $N \times N$ matrix. We still have a nonlinear equation (\ref{eigenvalue_per}) for the unknown $\xi$ which we can solve using Newton's method. On each iteration we still will have to solve a linear system of equations corresponding to condition $\Psi = 1$, but only $2M \times 2M$ size. \\
%The solution to (\ref{helmholtz_per_1})-(\ref{helmholtz_per_2}) is given by the single-layer potential (\ref{single_layer_potential}), except that the integration kernel will be different. The fundamental solution in this case is given by the periodic Green's function (\ref{periodic_Green_func_imag_sum}) with period $d = 2l_b+space$, though we will have to use Ewald representation (\ref{periodic_Green_func_Ewald}) for the numerical calculations. The analog of (\ref{bessel_decomp}) is
%$$ \frac{i}{4} \sum_{n=-\infty}^{\infty}{H_0^{(1)}(\xi r_n)} $$.

%----------------------------------------------------------------------------------------------------

\subsection{Numerical results and  comparison to  Ghergu et al.'s results}
%We first compare 
It is instructive to compare values of resonant frequencies 
$\xi_j$ solving for some $j$ in $\{1, .., N \}$, $p(\xi^2)\tau_j (\xi^2) = q(\xi^2) $
%$\xi_j, \, j=1 .. N$ 
coupled to  the $N$ identical, equally spaced, building
problem (\ref{helmholtz_1}-\ref{eigenvector_equation})
to values of resonant frequencies $\xi_{per}$ solving the periodic problem
(\ref{helmholtz_per_1}-\ref{eigenvalue_per}). We do that  in the case
where solutions  are sought in the vicinity of the initial guess $\xi_0=1$.
%In this section, for ease of notation,
%we denote by $\xi_{min}$ the solution to $p(\xi^2)\tau_1 (\xi^2) = q(\xi^2) $,
%and $\xi_{max}$ the solution to $p(\xi^2)\tau_N (\xi^2) = q(\xi^2) $.
Our numerical simulations point to the following observations:
\begin{itemize}
\item $\xi_j$ is increasing in $j$: $\xi_1 \leq .. \leq \xi_j \leq \xi_{j+1} \leq .. \leq \xi_N$
\item all the resonant frequencies  $\xi_j$ 
%solutions to  $p(\xi^2)\tau_j (\xi^2) = q(\xi^2) $
%
for  $j$ in $\{1, .., N \}$ vary within a narrow range.
 More precisely $\f{\xi_{N} - \xi_{1}}{\xi_{N}} $ and $\f{\xi_{N} - \xi_{1}}{\xi_{1}} $ are both small.
\item $\xi_{1} \leq \xi_{per} \leq \xi_{N}$ (up to numerical accuracy)
\item for smaller values of spacing between buildings $a_{j+1} - b_j$, $\xi_{per} \sim \xi_{N}$, and for
 larger values, $\xi_{per} \sim \xi_{1}$
\end{itemize}
Computed values of $\xi_{1}, \xi_{per}, \xi_{N}$ are given in Table \ref{compare_per_and_free}. 
$N$ is here set to be 51, and $M$ had to be set to as low as 5, for computational time to
be reasonable.
  The variable $space$ is for the spacing between buildings $a_{j+1} - b_j$.
The building half-widths are $l_j=1$. As previously, 
 $l = l_j$, $\gamma_j = 1.5$, $f_j = 0.5$, $c_j = 1$, $r_j=0.1$, ${\cal{B}}_j=1.5$.
%Below we provide the table comparing solutions of periodic problem with maximum and minimum eigenvalues of the corresponding free-space case. For all the city patterns $l_b=1$, $space$ is different. Free-space case: $N=51, \; M=5$; periodic case: $M=10$. We highlight the values if they look to coincide. \\
\begin{table}[H]
  \begin{tabular}{|c|c|c|c|c|c|c|c|c|c|}
    \hline
    $space$ & 0.5 & 1 & 1.3 & 1.4 & 1.5 & 2 & 3 \\
    \hline
    $\xi_{1}$ & 0.7821 & 0.7990 & 0.8156 & 0.8264 & 0.8418 & \textbf{0.8222} & \textbf{0.7933} \\
    \hline
		$\xi_{per}$ & \textbf{1.0864} & \textbf{0.9420} & 0.8873 & 0.8737 & 0.8619 & \textbf{0.8225} & \textbf{0.7934} \\
    \hline
    $\xi_{N}$ & \textbf{1.0844} & \textbf{0.9408} & 1.0391 & 1.0514 & 1.0602 & 1.0635 & 0.9772 \\
    \hline
  \end{tabular}
  \caption{Numerical values for $\xi_1, \xi_{per}, \xi_N$ defined above.
	Here $N=51$, $M=5$ ($2M$ is the number of grid points on each building).
	The variable $space$ indicates the spacing between buildings.
	%The variable $space$ is for the spacing between buildings $a_{j+1} - b_j$.
The building half-widths are $l_j=1$. As previously, 
 $l = l_j$, $\gamma_j = 1.5$, $f_j = 0.5$, $c_j = 1$, $r_j=0.1$, ${\cal{B}}_j=1.5$.}
	%Comparison of the free-space and periodic frequencies $\xi_{per}$  vs $\xi_{1}, \; \xi_{N}$; building halfwidth $l_b=1$; the number of buildings $N=51$, $M=5$ ($2M$ is the number of grid points on each building).}
  \label{compare_per_and_free}
\end{table} 
%It can be seen that when $space$ is small enough, $\xi_{per}$ coincides with $\xi_{max}$, and when it is big enough - with $\xi_{min}$. Tables \ref{tab:DifferentEigenvectors1} and \ref{tab:DifferentEigenvectors2} show that the corresponding eigenvectors keep the sign. Then, for some range of $space$ $\xi_{min} < \xi_{per} < \xi_{max}$. It was noticed before that for such cities no eigenvector keeps the sign.

%We would like to look deeper into this matter. Is there any connection between $\xi_{per}$ and eigenvalues for such cities? First, let us consider the city where $l_b=1, \; space = 1.7$ and find all of its eigenvalues. We will do it for different number of buildings $N=11; \; 31; \; 51$. Again, for every simulation $M=5$. It appears that $\xi_{per} \approx \xi_3$.  Also, we show all the eigenvalues for $space = 0.5, \; space =3$.  All the results are presented in Table \ref{tab:all_eigenvalues}. \\
\begin{table}[H]
  {\small
  \begin{tabular}{c }
 %   \includegraphics[scale=.35]{all_eigenvalues_space17.pdf} & \includegraphics[scale=.35]{all_eigenvalues_space17_zoom.pdf} 
%    \\
 %   \multicolumn{2}{c} {$space=1.7$} \\
    \includegraphics[scale=.5]{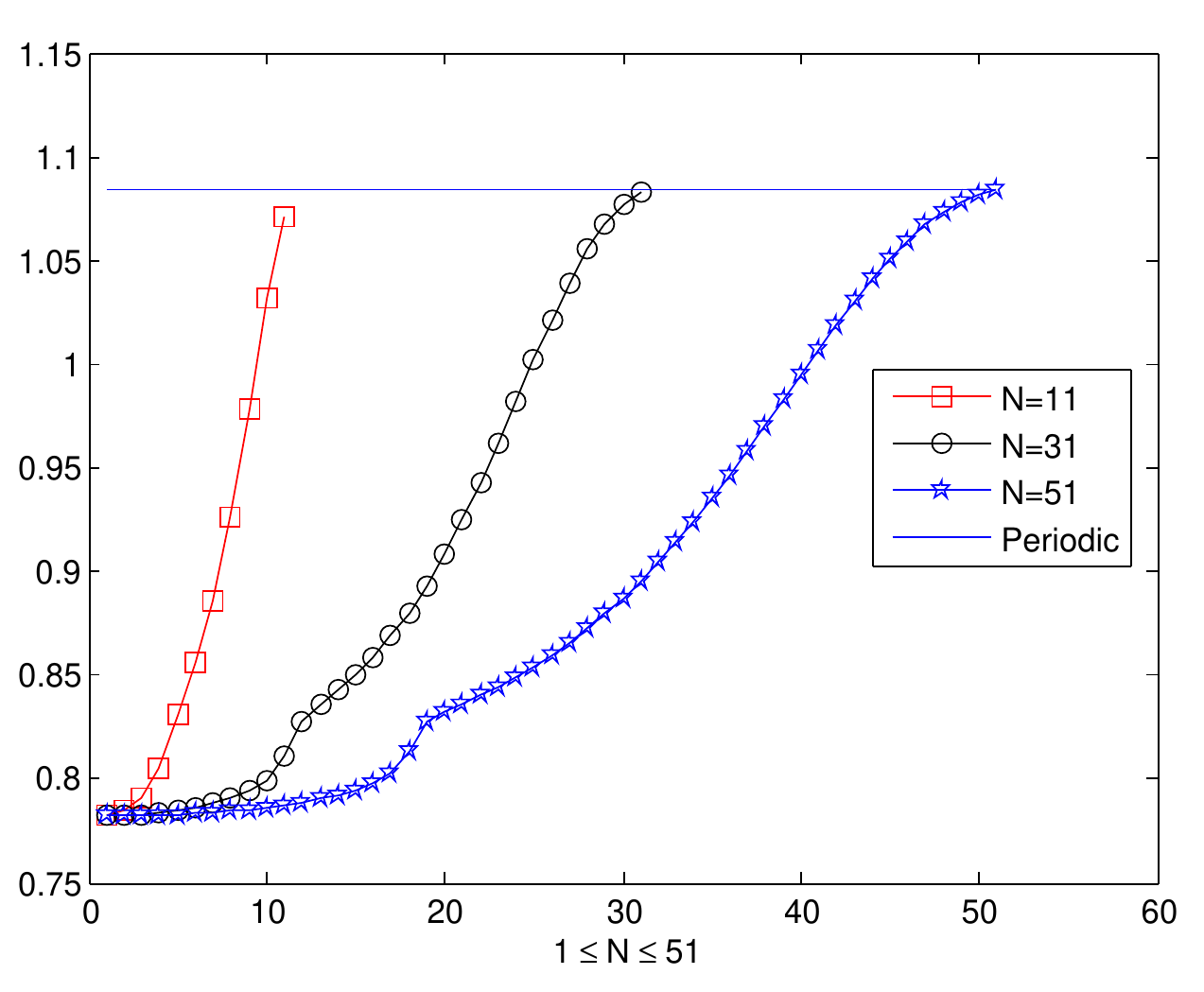} %& \includegraphics[scale=.35]{all_eigenvalues_space05_zoom.pdf} 
    %\multicolumn{2}{c} 
     
    \includegraphics[scale=.5]{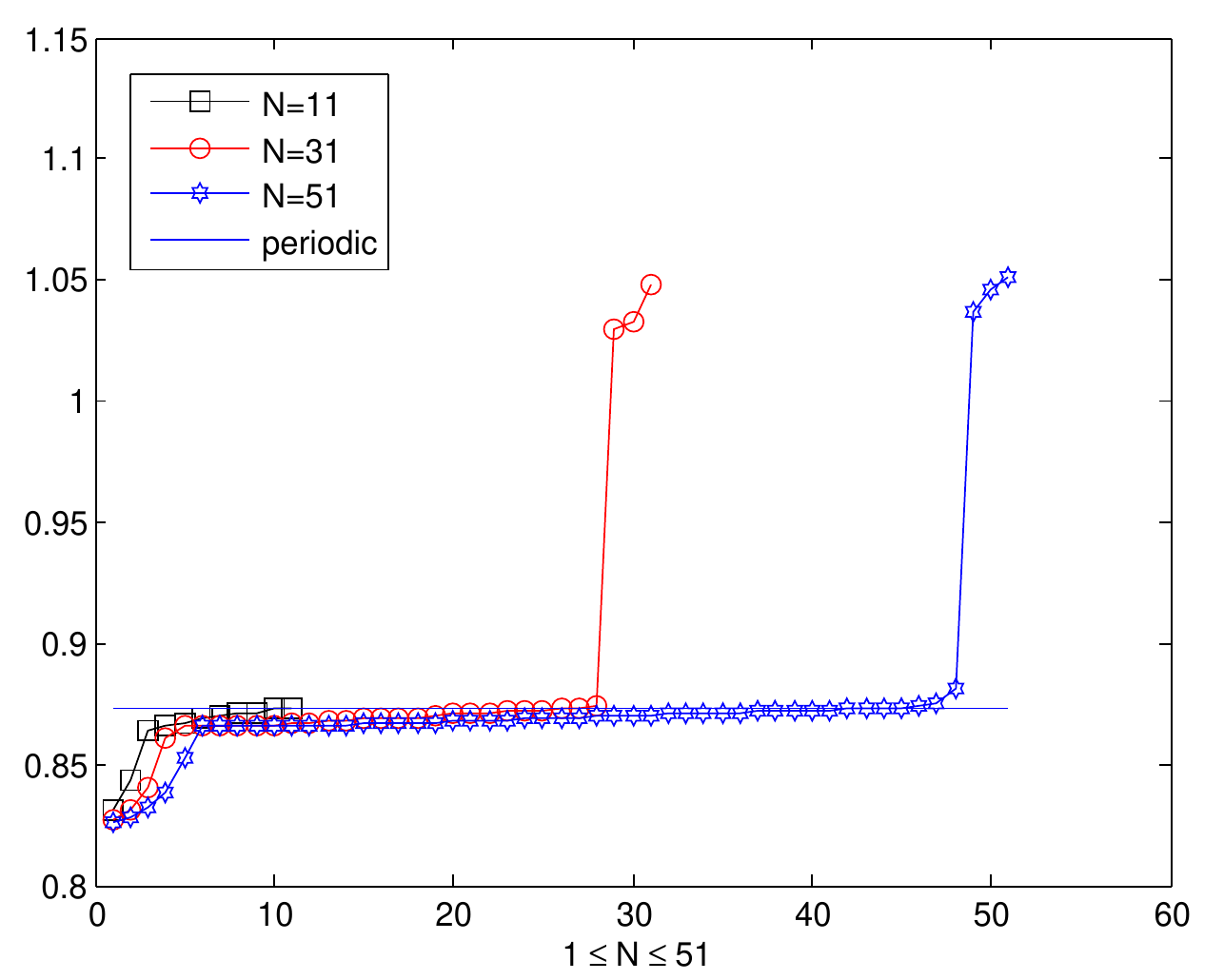}  \\ %\includegraphics[scale=.35]{all_eigenvalues_space05_zoom.pdf} 
 % \\
    %\multicolumn{2}{c} {
   \hskip .3in  $space=0.5$  \hskip 2in $space=1.7$  \\
    \includegraphics[scale=.5]{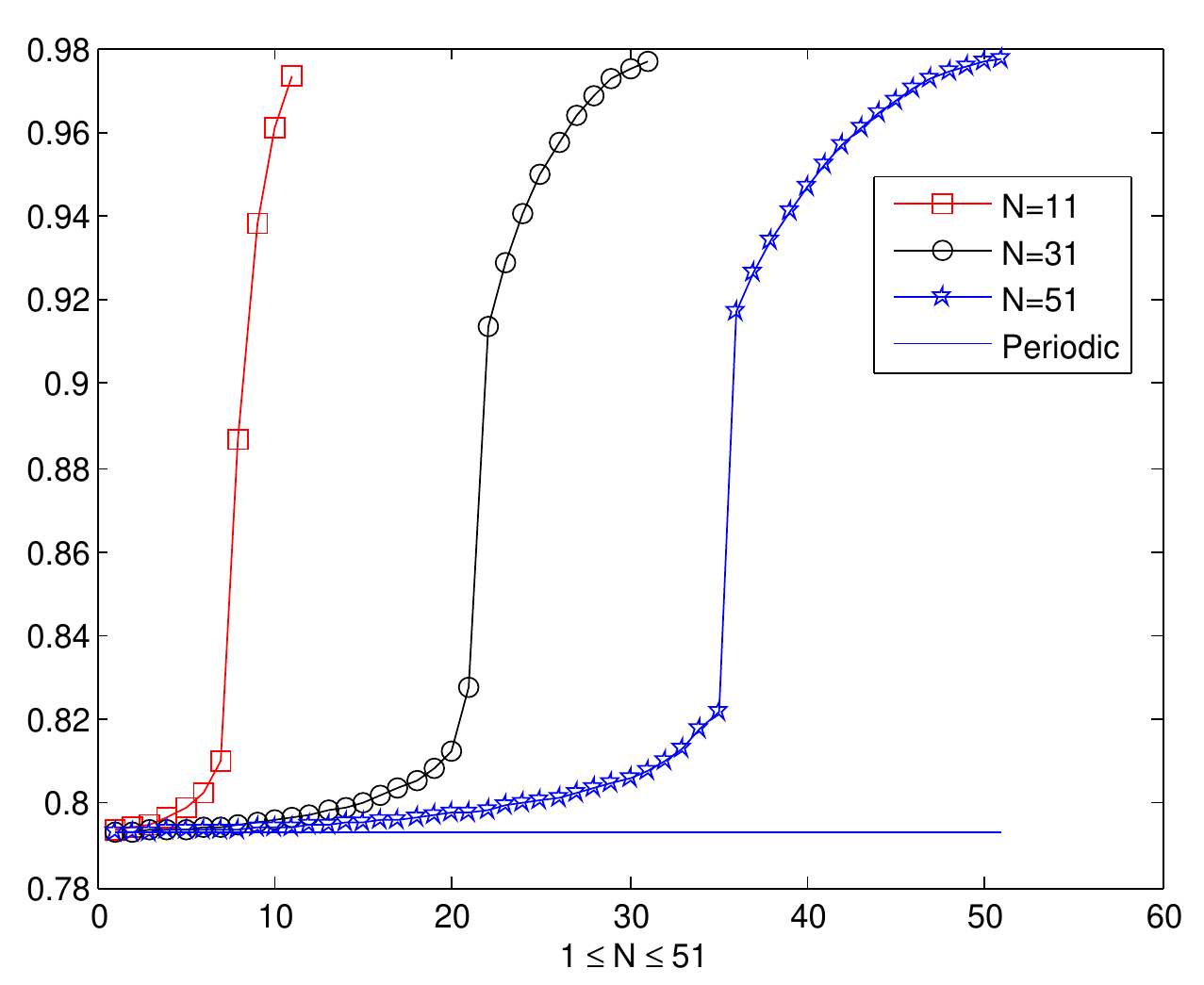} \\ %& \includegraphics[scale=.35]{all_eigenvalues_space3_zoom.pdf} \\
    %\multicolumn{2}{c} 
   \hskip .8in $space=3$  
  \end{tabular} }
  \caption{Vertical axis: resonant frequencies 
$\xi_j$ solving  $p(\xi^2)\tau_j (\xi^2) = q(\xi^2) $
	coupled to  the $N$ identical, equally spaced, building
problem (\ref{helmholtz_1}-\ref{eigenvector_equation}). Horizontal axis: value of integer $j$.
 Horizontal lines: value of $\xi_{per}$. Three cases corresponding to different spacings between buildings,
0.5, 1.7, and, 3, are shown. For each spacing we examined the case where the total number of buildings $N$ is
11, 31, or, 51.
%	
%	for the free-space problem (\ref{non_dim_helmholtz})-(\ref{non_dim_eigenvalue}) compared to the frequency of the periodic problem (\ref{helmholtz_per_1})-(\ref{eigenvalue_per}). In each case 
%	the building halfwidth is $l_b=1$; the distance between consecutive buildings ``$space$'' is first .5, then 1.7, and finally 3. The number of buildings $N$ is $11$, $31$, and $51$; $M=5$. 
 As previously, the physical parameters
%	are
 $l = l_j=1$, $\gamma_j = 1.5$, $f_j = 0.5$, $c_j = 1$, $r_j=0.1$, ${\cal{B}}_j=1.5$.}
  \label{tab:all_eigenvalues}
\end{table}

\section{The case of sets of different buildings}
\label{sec:nonhomogenuous_city}
%subsection{Small size city}
%\
In this section the physical parameters
 $m_{1,j}$, $m_{0,j}$, 
$l_j$, $h_j$, $\rho_j$, $\beta_j$ may depend on the building $j$, and the lengths
 $b_j - a_j$ and $a_{j+1} - b_j$ too.
Accordingly the two functions $p_j$ and $q_j$ defined in 
(\ref{p_and_q}) may also depend on $j$, and consequently the symmetric matrix method used in
section \ref{matrix_method} is no longer applicable.

\subsection{Case where the number of buildings $N$ is finite}

In this case we have to solve for $\xi$ the system of $N$ non linear equations

\begin{equation}
  \label{nonlin_system}
  \alpha_i q_i(\xi^2)  = p_i(\xi^2) \mbox{Re } \int_{\Gamma_i} {\frac{\partial \Psi} {\partial y} (s,0)} ds, \q i=1, ..., N
\end{equation}
where $\Psi$ and the real numbers $\alpha_i$
 are coupled through the PDE (\ref{helmholtz_1}-\ref{sommerfeld_cond}).
For each step in the search for $\xi$ satisfying
(\ref{nonlin_system}), that PDE is solved by the same integral equation method as previously.

%The previous results were obtained under stringent assumptions. It is more realistic to allow buildings to be of different size and rigidity. We now allow buildings to have different heights, foundation areas, distances between each other, and to be built from various materials. In terms of our physical model, it means that not only displacements, but all the other parameters defined in section \ref{sec:physical_model} depend now on the building number $j$. The same is true for non-dimensional parameters (\ref{non_dimensional_parameters}), and the polynomials $p$ and $q$ defined by (\ref{p_and_q}). It is reasonable to conjecture that eigenvalues and eigenvectors of (\ref{linear_operator}) do not solve problem (\ref{helmholtz_1}) - (\ref{sommerfeld_cond}). We will illustrate this for a 2-building city with building left endpoints $a=[-2.5; 1.5]$, building right endpoints $b=[-1.5; 3]$, and distance between the buildings $space=3$. \\
%\begin{table}[H]
 % \centering
 % \begin{tabular}{c c}
 %   \begin{cases} 
 %     q_1(\xi^2) = p_1(\xi^2) \tau_1(\xi^2) \\
 %     q_2(\xi^2) = p_2(\xi^2) \tau_1(\xi^2)
 %   \end{cases} &
 %   \begin{cases}
 %     q_1(\xi^2) = p_1(\xi^2) \tau_2(\xi^2) \\
 %     q_2(\xi^2) = p_2(\xi^2) \tau_2(\xi^2)
 %   \end{cases}
 % \end{tabular}
%\end{table}

\begin{figure}[H]
 \centering
 \includegraphics[scale=.5]{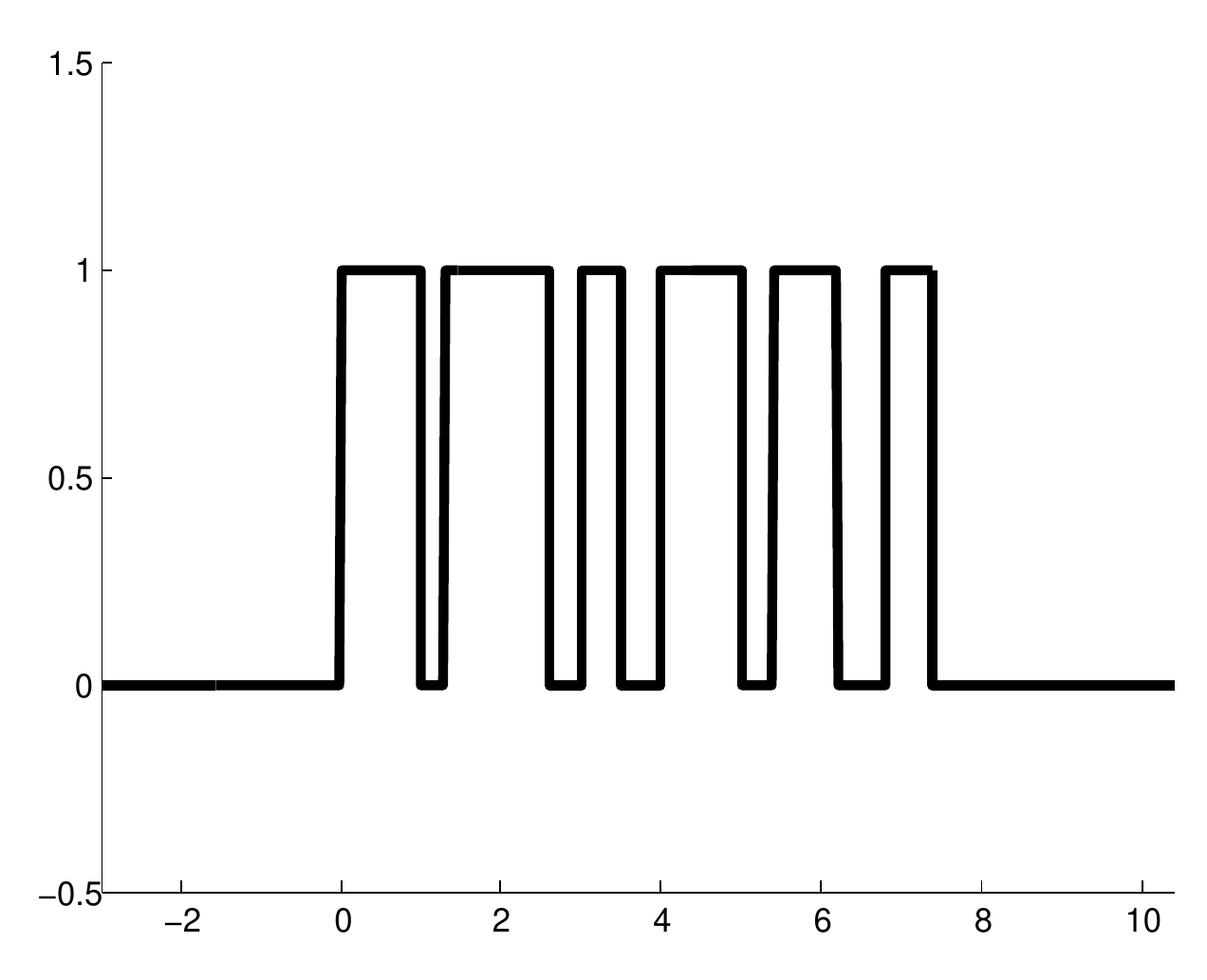}
 \caption{A sketch illustrating the relative size of the foundations and the spacing between the buildings for the 6-building city defined by 
$a =  ( 0, 1.3,  3,  4,  5.4, 6.8)$, $b =  (1,  2.6,  3.5,  5,  6.2, 7.4)$.}
 \label{fig:6build_city}
\end{figure}

Let us now present some numerical results for this new case.
We will vary building half-widths $l_{j}$ and distances from one building
to the next  $a_{j+1} - b_j $. 
%$space_j= a_{j+1} - b_j $. 
%This will lead us to different $p_j(\xi^2)$ and $q_j(\xi^2)$. 
%We recall that for the homogeneous case characteristic length $l$ was set up equal to $l_b$, such that $c_b=1$ (see (\ref{parameter_values})). We can not proceed in the same manner, because $l_{bj}$ are different, and $l$ should be unique for our problem. 
We choose the characteristic length $l=1$, so that  that $c_{j} = l_{j}$. 
%All the other parameters (\ref{parameter_values}) will be unchanged and equal for all the buildings, but we notice that altering them will not change anything in our methods and reasoning. 
To facilitate comparison to previous cases, we will let $l_j$ vary while the other physical parameters 
will remain constant from one  building to another and their values will be the same as previously. 
In other words, $\gamma_1 =\gamma_j= 1.5$, $f_1 = f_j= 0.5$,  $r_j=r_1=0.1$, ${\cal{B}}_j={\cal{B}}_1=1.5$.
Clearly, since  $l_j$ and $a_{j+1} - b_j$ are non constant,
this is will lead to different functions $p_j$ and $q_j$ as $j$ varies. \\
Suppose that we can find a solution to (\ref{helmholtz_1}-\ref{sommerfeld_cond}, \ref{nonlin_system})
such that for a particular index $j$, $\alpha_j \neq 0$.
Then by linearity we can assume that $\alpha_j =1$.
Let us now impose $\alpha_j =1$ in the search for a solution to
(\ref{helmholtz_1}-\ref{sommerfeld_cond}, \ref{nonlin_system}).
Our numerical simulations clearly indicate that solutions depend on the choice of such a $j$. 
In 
Table \ref{tab:six_building} we show solutions $(\xi, \alpha)$ for a 6-building geometry sketched in 
Figure \ref{fig:6build_city}, where 
$a = ( 0, 1.3,  3,  4,  5.4, 6.8)$, $b = (1,  2.6,  3.5,  5,  6.2, 7.4)$ and  $M=10$ ($2M$ was set 
to be the number of grid points per building in our numerical calculations). The values  $\alpha_j$ of $\Psi(x,0)$ on $\Gamma_j$
are shown as bar graphs, and coupling wavenumbers $\xi$ are given below each graph. 
We observe that the runs for $\alpha_2=1$, $\alpha_3=1$, and $\alpha_6=1$
all lead to the same eigenvalue $\xi= 1.3660$. A closer look at Table \ref{tab:six_building}
reveals that these three runs also lead, after rescaling,  to the same \textbf{eigenvector}. 
In other words, the solutions in these three cases are clearly multiple of one another.
The same remark can be made for the runs for $\alpha_1=1$ and $\alpha_5=1$.

%Remark. \\
%We use $Matlab$ to solve the nonlinear systems of equations. The termination conditions are set up in such a way that the final answer accuracy is $10^{-6}$. We notice that this is not the accuracy of our method. \\
%\end{remark_ch3}

\begin{table}[H]
{\scriptsize
  \begin{tabular}{c c c}
    \includegraphics[width=0.3\textwidth]{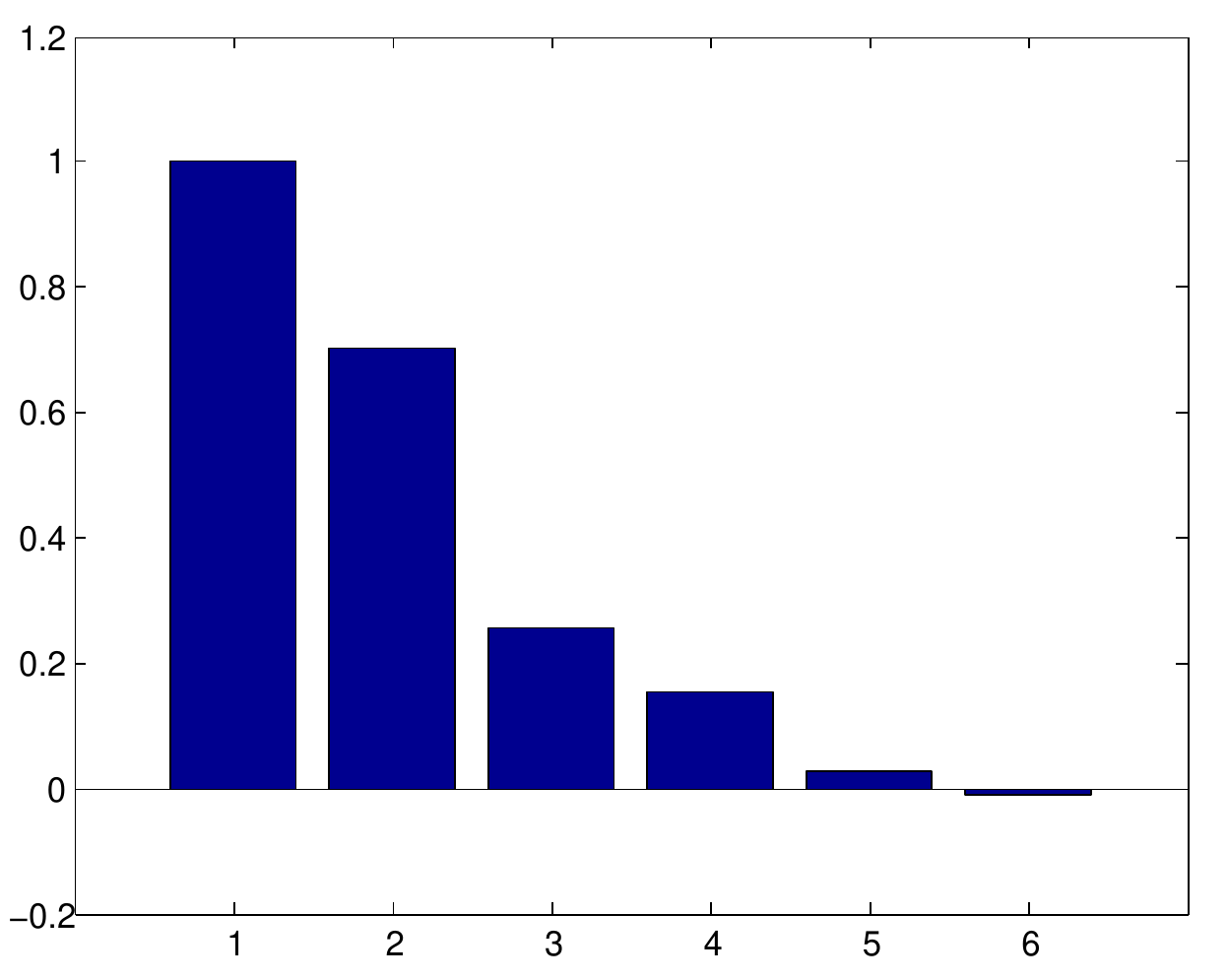} & \includegraphics[width=0.3\textwidth]{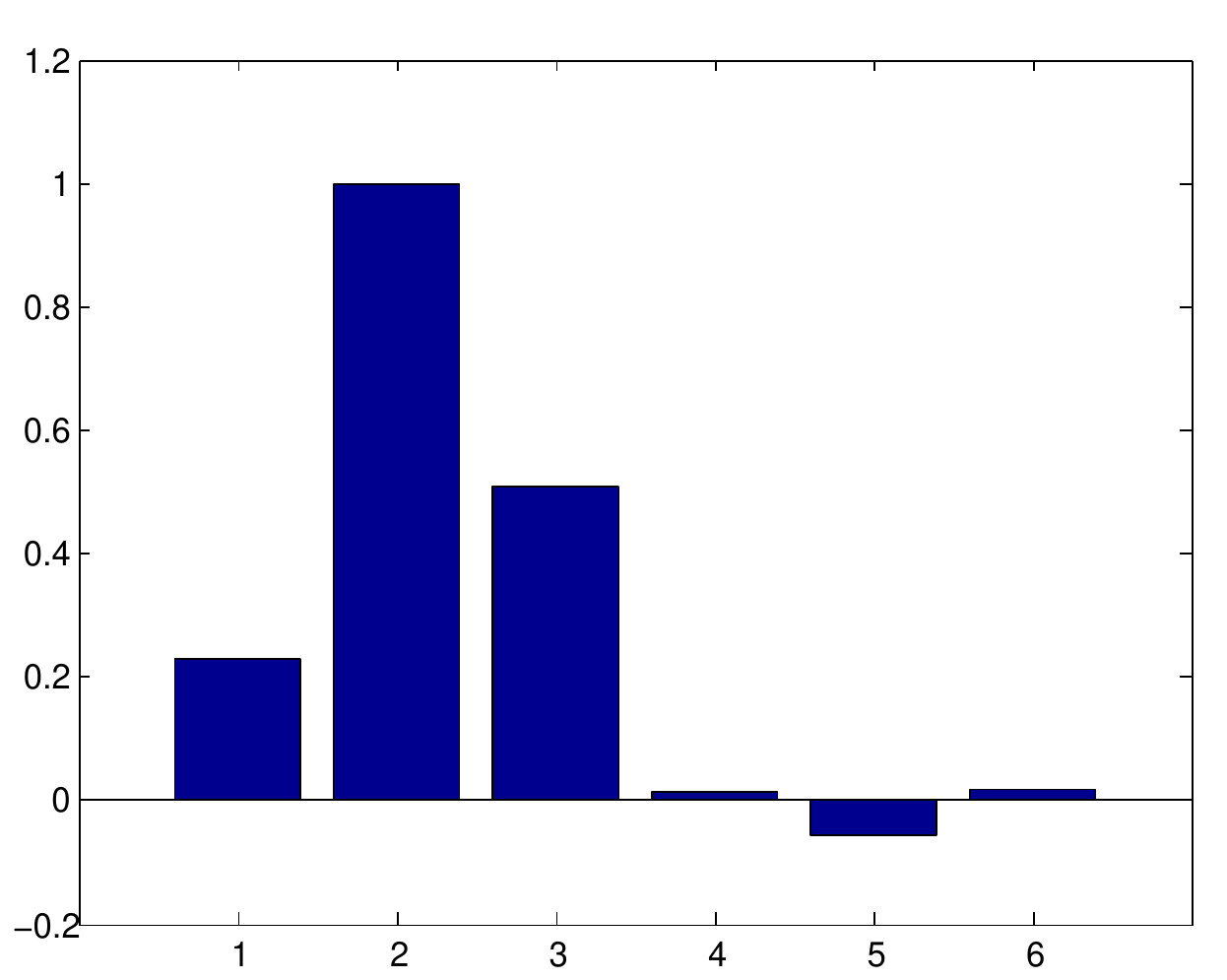} & \includegraphics[width=0.3\textwidth]{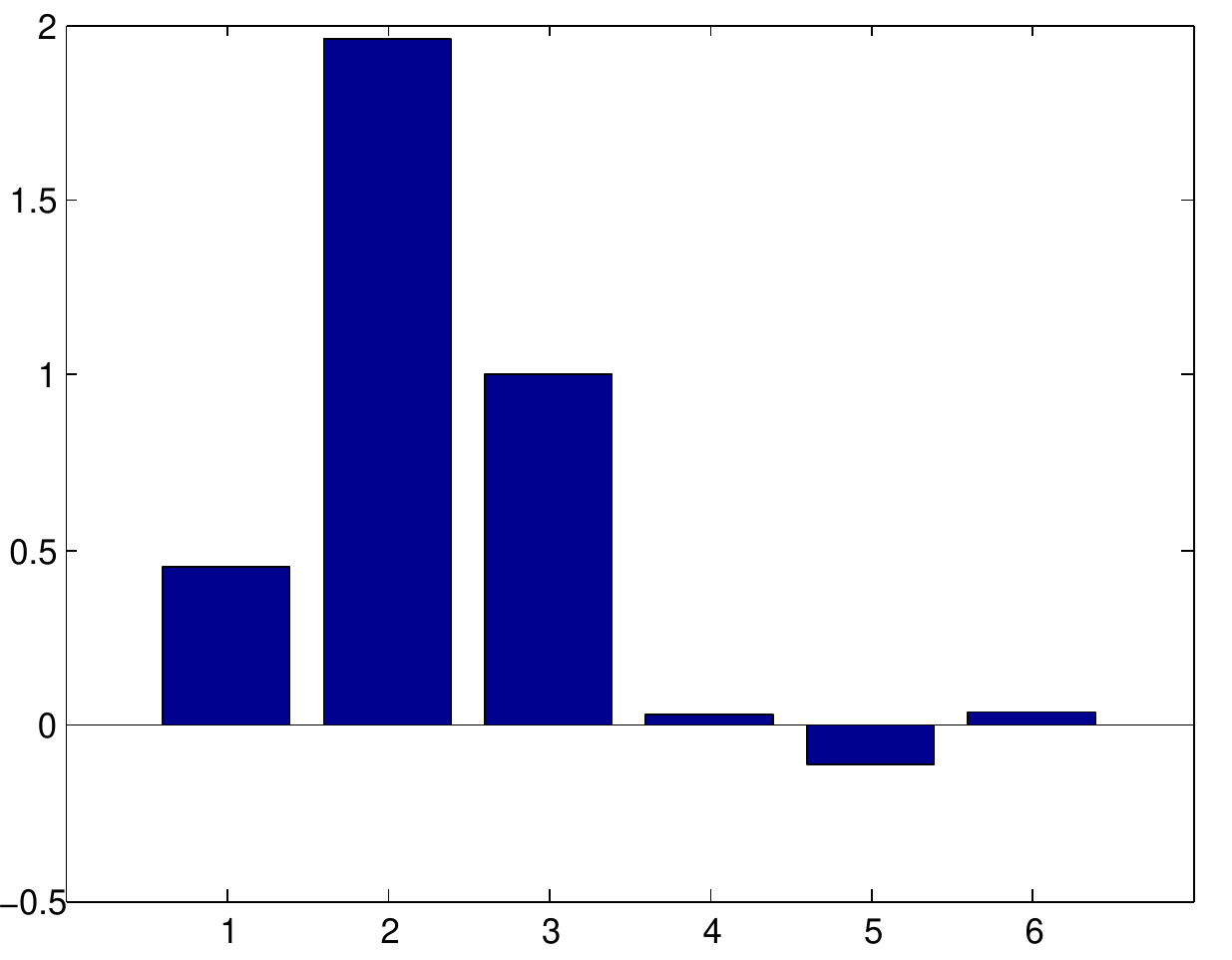} \\
    $\alpha_1 = 1, \; \xi = 1.7301$ & $\alpha_2 = 1, \; \xi =1.3660$ & $\alpha_3 = 1, \; \xi = 1.3660$ \\
    \includegraphics[width=0.3\textwidth]{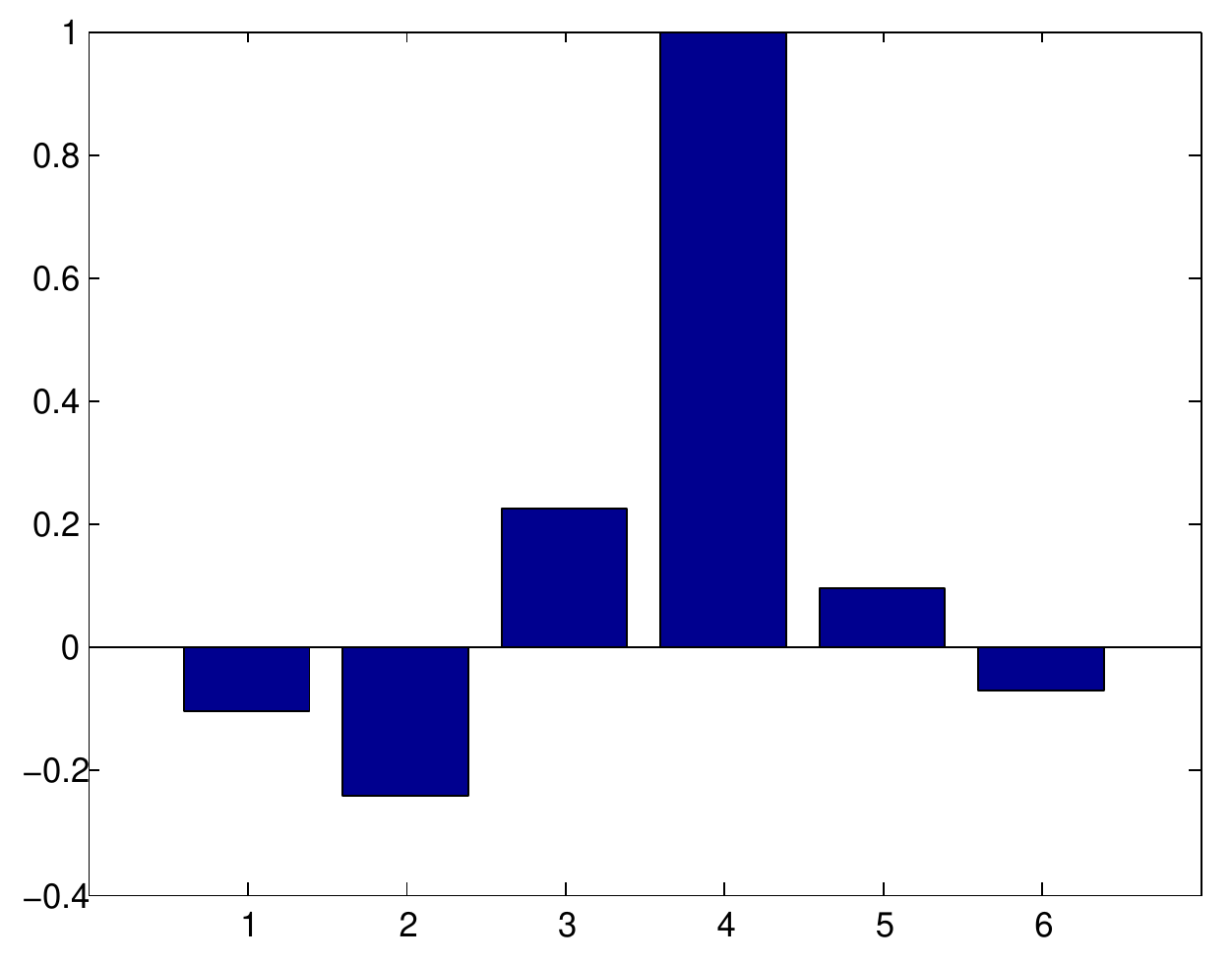} & \includegraphics[width=0.3\textwidth]{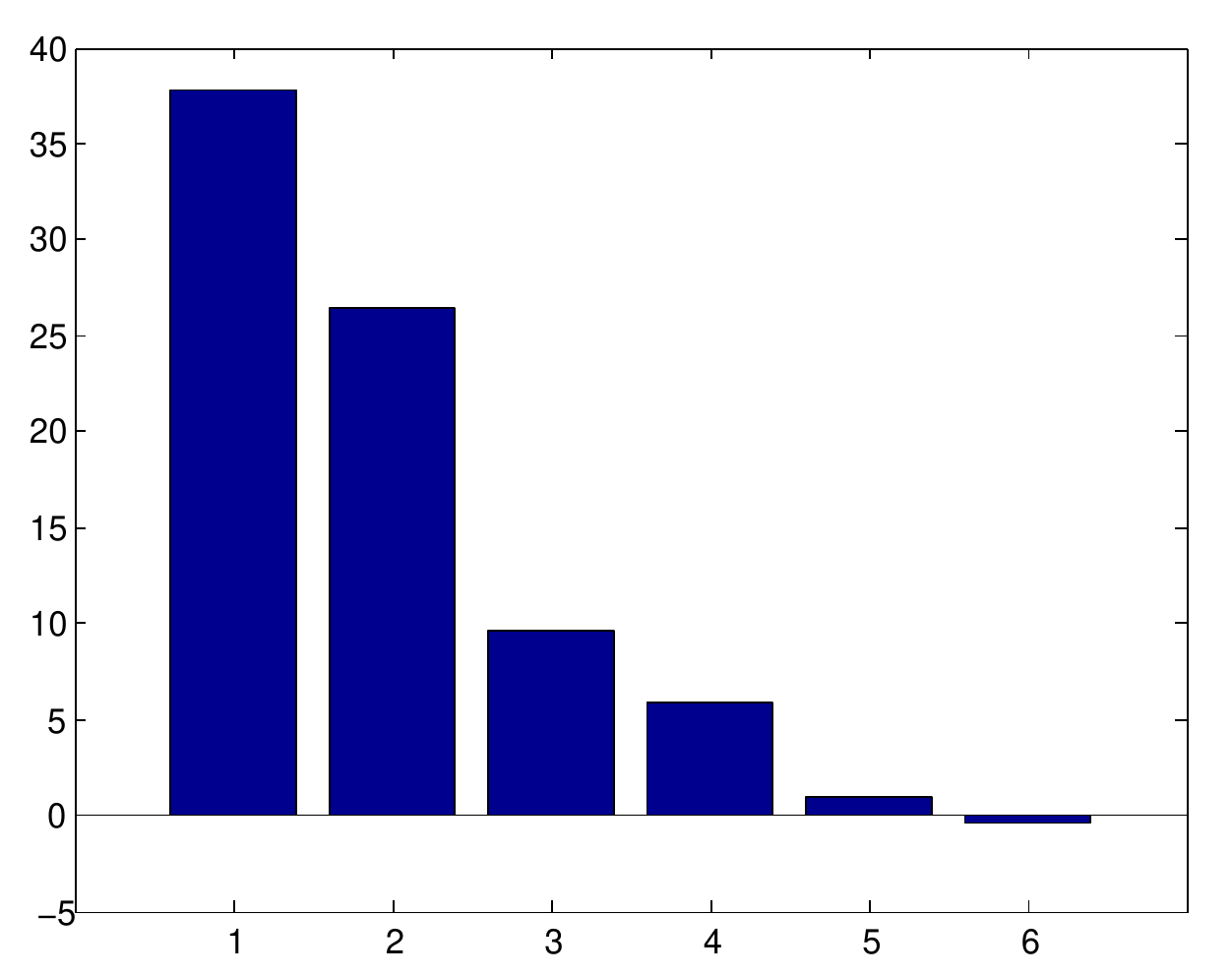} & \includegraphics[width=0.3\textwidth]{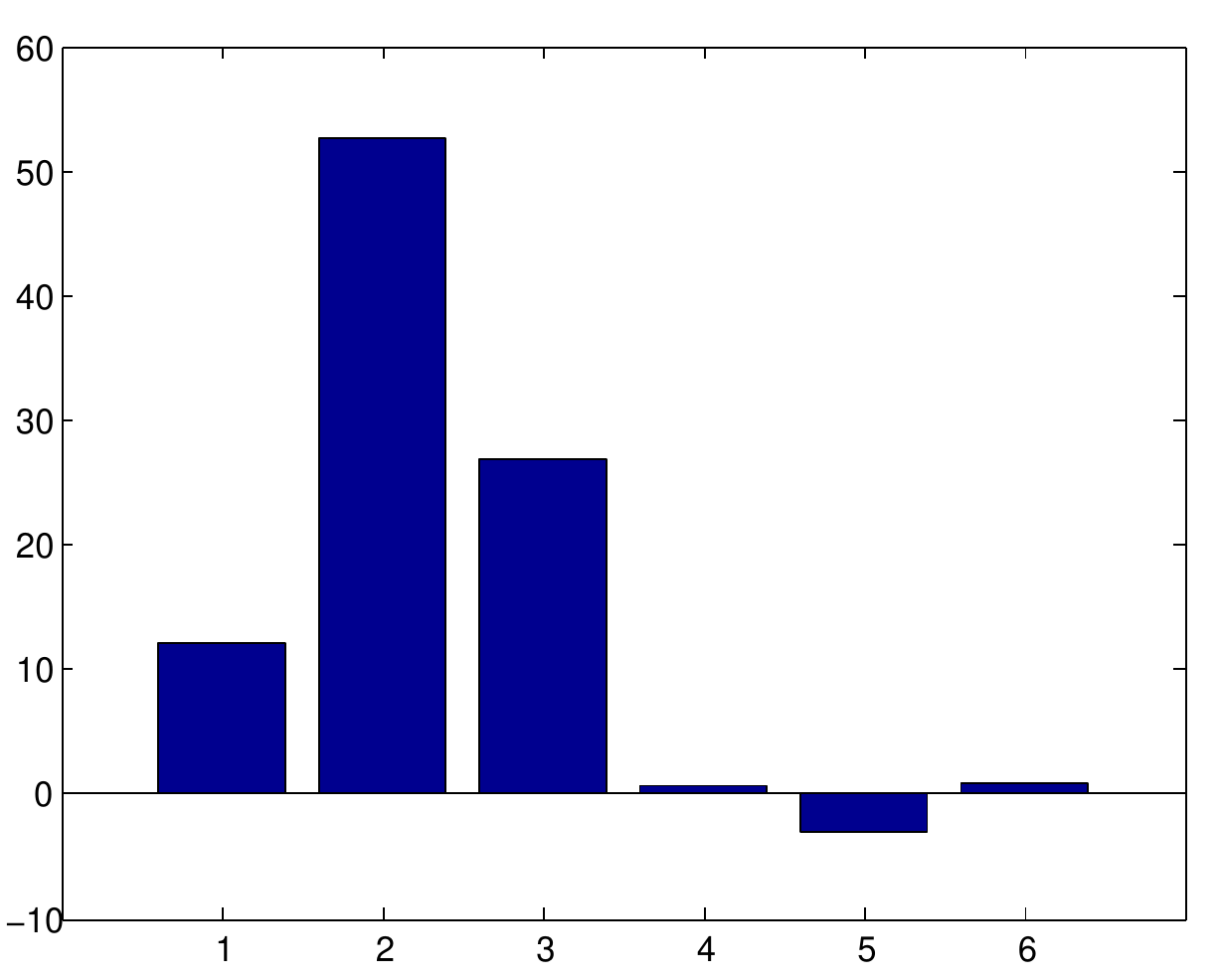} \\
    $\alpha_4 = 1, \; \xi = 1.7784$ & $\alpha_5 = 1, \; \xi =1.7301$ & $\alpha_6 = 1, \; \xi =1.3660$
  \end{tabular} }
\caption{Solutions for the 6-building city: $a =  ( 0, 1.3,  3,  4,  5.4; 6.8)$, 
$b =  (1,  2.6,  3.5,  5,  6.2, 7.4)$, $M=10$. The runs for $\alpha_2=1$, $\alpha_3=1$, and $\alpha_6=1$
all lead to the same eigenvalue $\xi= 1.3660$. A closer look 
reveals that these three runs also lead, after rescaling,  to the same eigenvector. }
\label{tab:six_building}
\end{table}

As expected, the solution for each case $\alpha_j=1$ is not  unique. The results in Table \ref{tab:six_building} 
were obtained for an initial guess for the wavenumber $\xi_0 = 1$. Choosing instead $\xi_0=2.5$, we obtain the following: \\
if we impose $\alpha_1 = 1$, or $\alpha_2 = 1$, or $\alpha_3 = 1$, then our computation results in $\xi = 1.7301$; \\ 
if we impose $\alpha_4 = 1$, or $\alpha_5 = 1$, then our computation results in $\xi=2.1861$; \\
if we impose $\alpha_6=1$, then our computation results in $\xi =2.8057$. \\

%Ultimately, we want to assume that some building cluster is repeated periodically in a city, and apply periodic Green's function to find city frequencies. That is why in the free-space case we will consider cities with several equal clusters of a few buildings. Then number of buildings in the city $N$ equals $N_b N_c$, where $N_b$ is the number of buildings in one cluster, and $N_c$ is the number of clusters. Below we present the results for three different geometries. First, it will be the 2-building city that we have already mentioned. Then we consider two 3-building cities. We will refer to them as city7.5, city7, because their lengths of the periodic clusters are 7.5, 7 . 
Let us now consider consider sets of buildings, referred to as cities, such that a pattern of buildings is repeated finitely many times.
As previously $N$ denotes the total number of buildings and we now set $N = N_c B$ where $N_c$ is the number of repeated patterns or cells,
and $B$ is the number of buildings in each cell. We now show numerical results for two examples of city
geometry.
They will be referred to as City 1 and City 2. \\
\begin{enumerate}
\item City 1 : B=2, $(a_1,a_2)=(-2.5, 1.5), \; (b_1,b_2)=(-1.5,3), \; a_{j+2} - a_j =b_{j+2} - b_j = 7.5 $. 
  
\item City 2: B=3, $(a_1,a_2,a_3)=(0, 2, 5), \;  (b_1,b_2,b_3)=(1.2,3, 6.7), \; a_{j+3} - a_j =b_{j+3} - b_j = 7$.
  
  \end{enumerate} 
\begin{figure}[H]
    \centering
    \includegraphics[scale=.5]{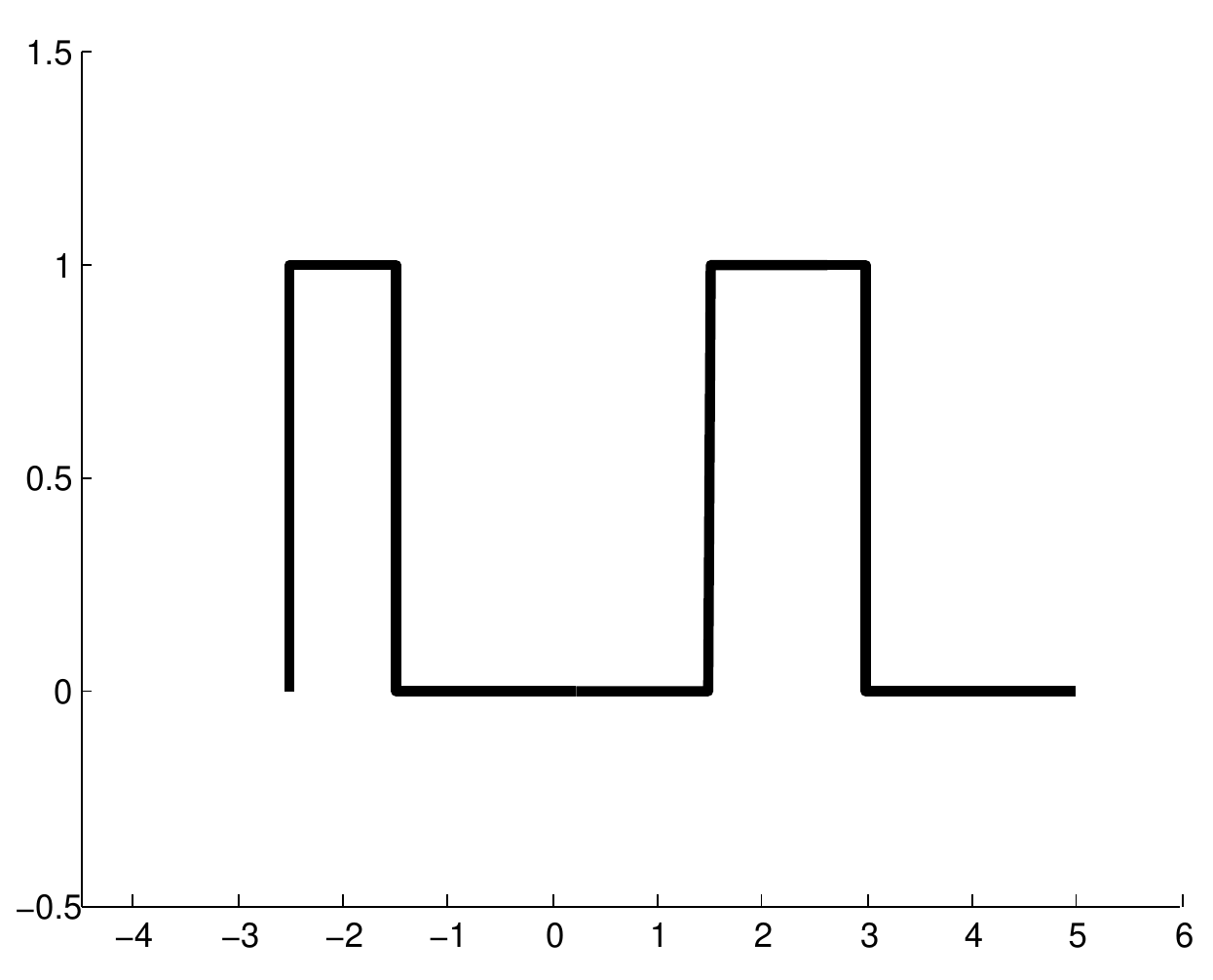}
    \caption{A sketch illustrating the relative size of the foundations and the spacing between the buildings for City 1 defined above.}
    \label{fig:city75}
  \end{figure}
	\begin{figure}[H]
    \centering
    \includegraphics[scale=.5]{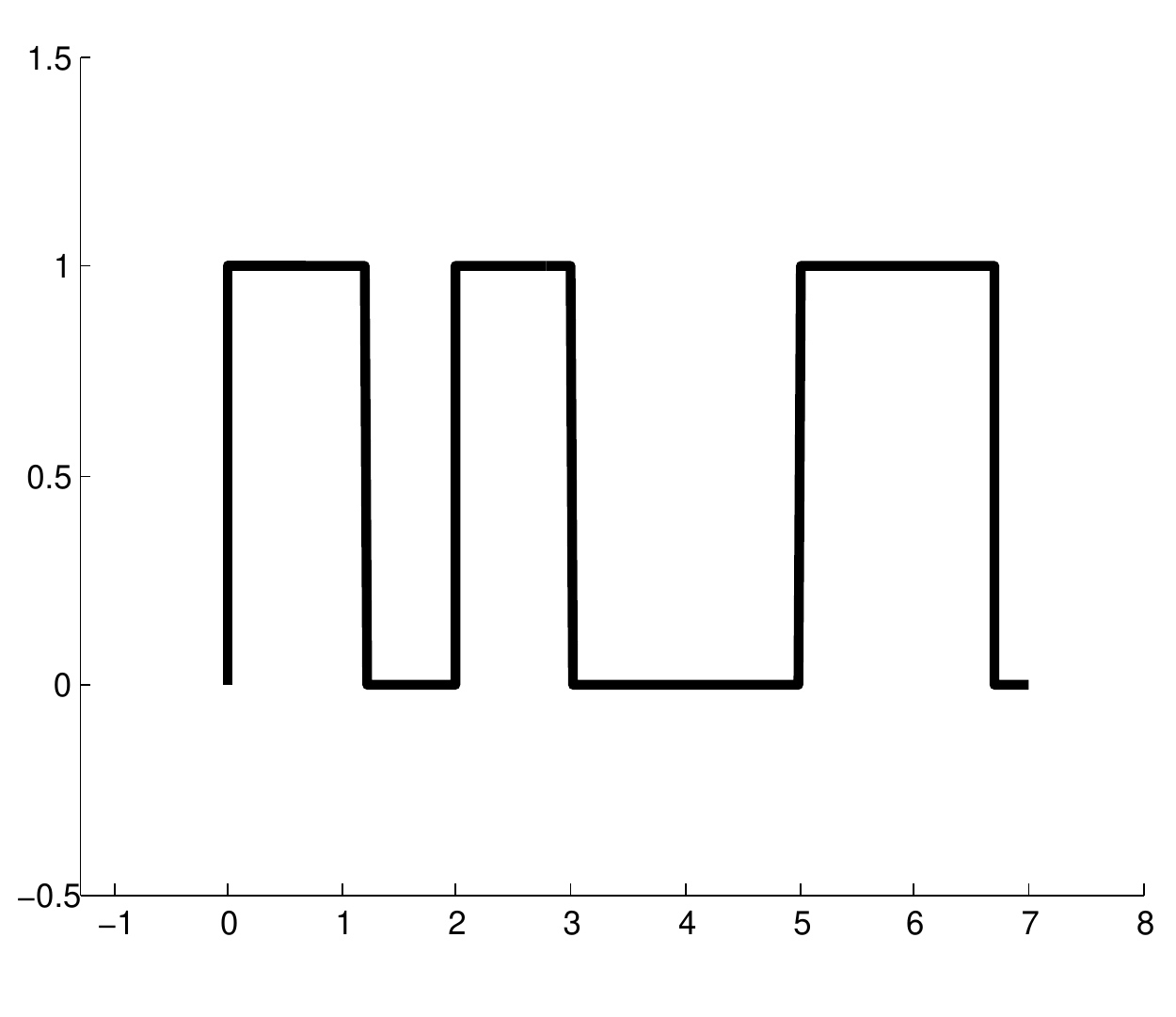}
    \caption{A sketch illustrating the relative size of the foundations and the spacing between the buildings for City 2 defined above.}
    \label{fig:city7}
  \end{figure}
Tables \ref{tab:city75_M5_N12} and \ref{tab:city7_M5_N15}  show some computed solutions.
They depend on the geometry of each city, including number of clusters $N_c$, on the the choice of building $j$
where the condition $\alpha_j =1$ is imposed and  on the initial value $\xi_0$ for $\xi$ for the search
of a solution to non linear equation (\ref{nonlin_system}).
Here too, we observe that runs that differ in the choice of the index $j$ for which 
we impose $\alpha_j=1$ may eventually lead to the same eigenvector (after rescaling).
In Table \ref{tab:city75_M5_N12} this occurs for example, 
$\alpha_6=1$, $\alpha_{10}=1$, $\alpha_{12}=1$.

%In the next subsection we will present the data that demonstrates that they converge as the number of clusters increases (see Tables \ref{tab:asymptot_for_city7.5} and \ref{tab:asymptot_for_city7}). Also, at this point we would like to emphasize $\xi=1.1584$ for city7.5 and $\xi=1.0349$ for city7. We will meet these again when we solve the periodic problem. 

\begin{table}[H]
{\scriptsize
  \begin{tabular}{c c c}
    \includegraphics[width=0.3\textwidth]{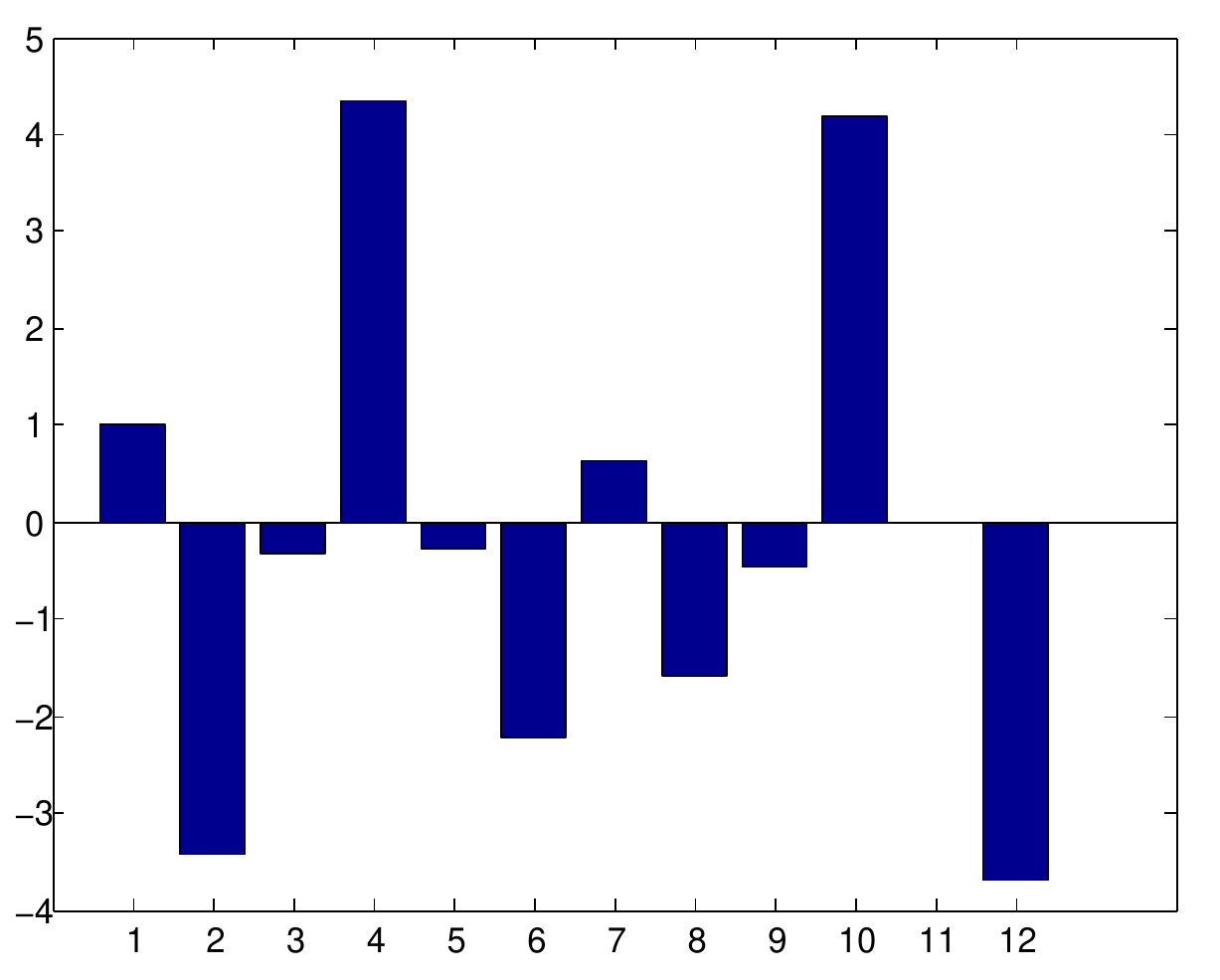} & \includegraphics[width=0.3\textwidth]{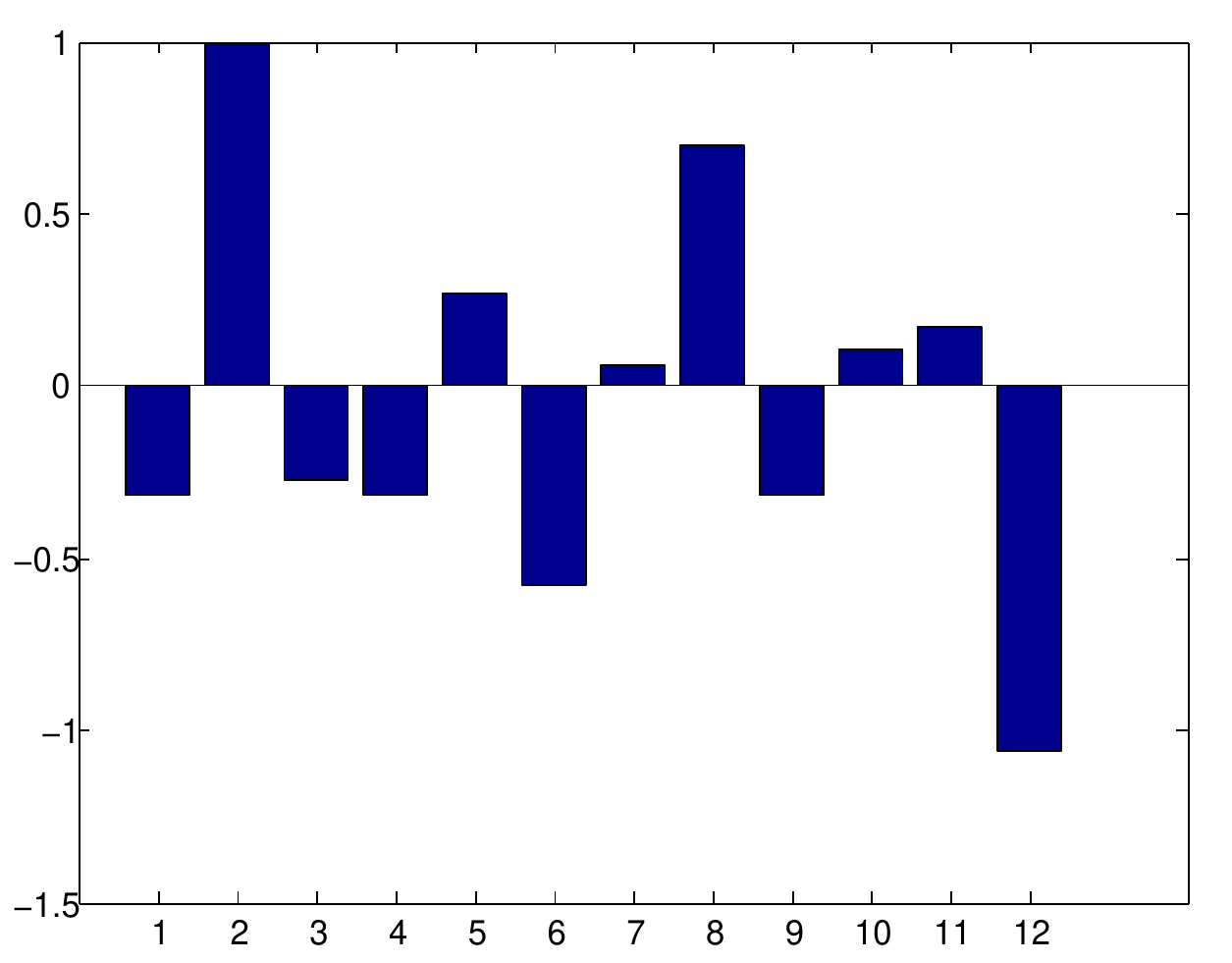} & \includegraphics[width=0.3\textwidth]{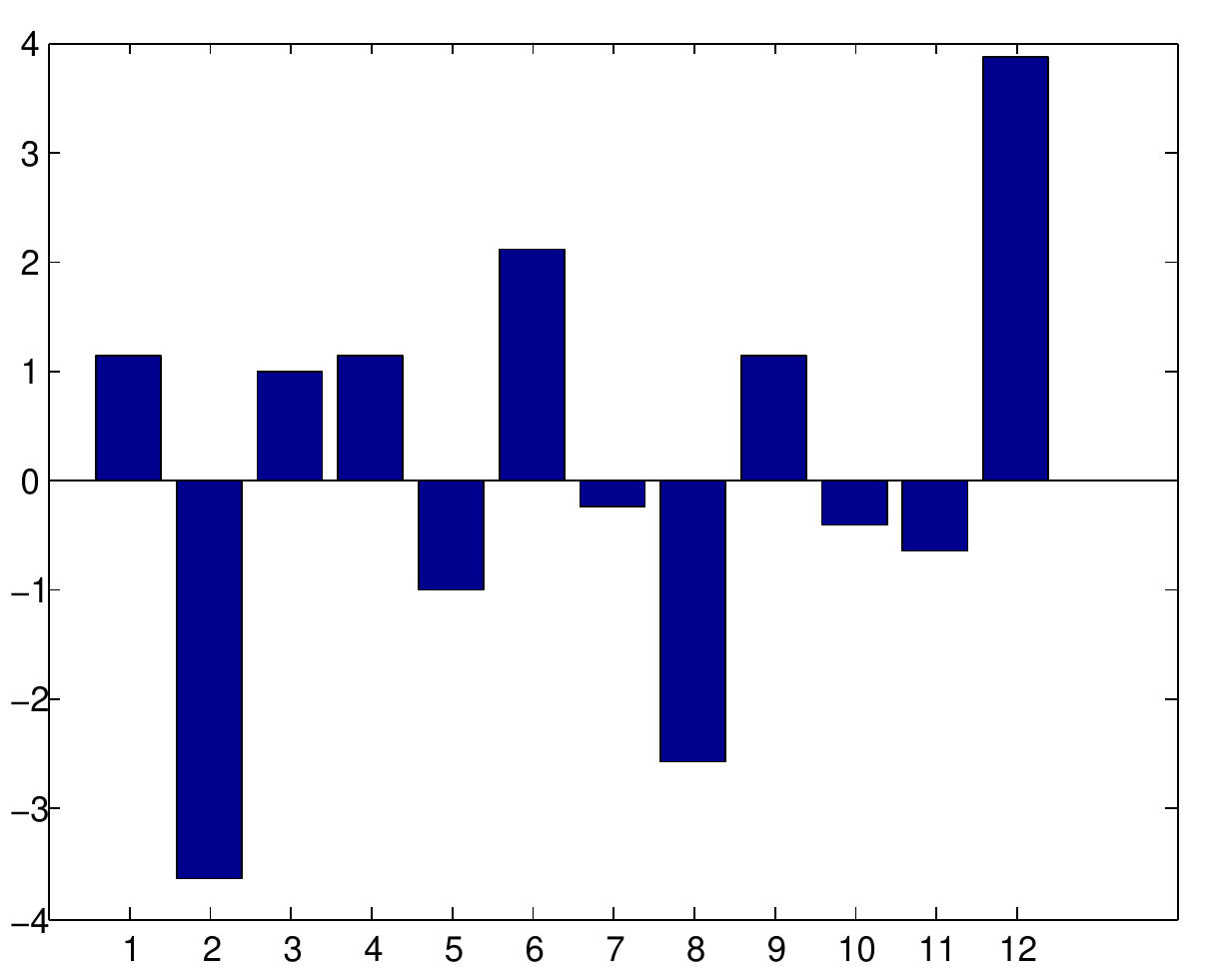} \\
    $\alpha_1 = 1, \; \xi = 1.0778$ & $\alpha_2 = 1, \; \xi =1.1212$ & $\alpha_3 = 1, \; \xi = 1.1212$ \\
    \includegraphics[width=0.3\textwidth]{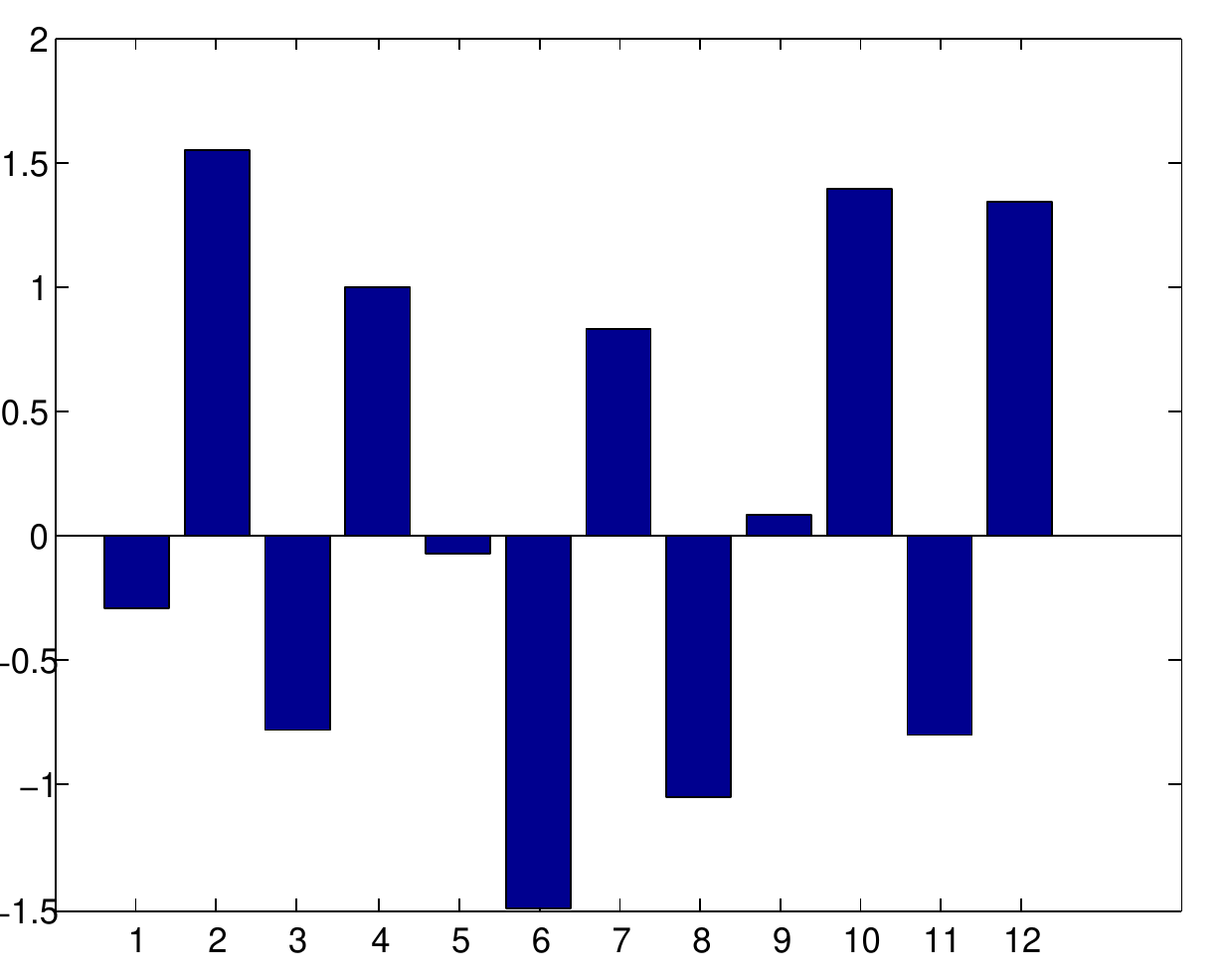} & \includegraphics[width=0.3\textwidth]{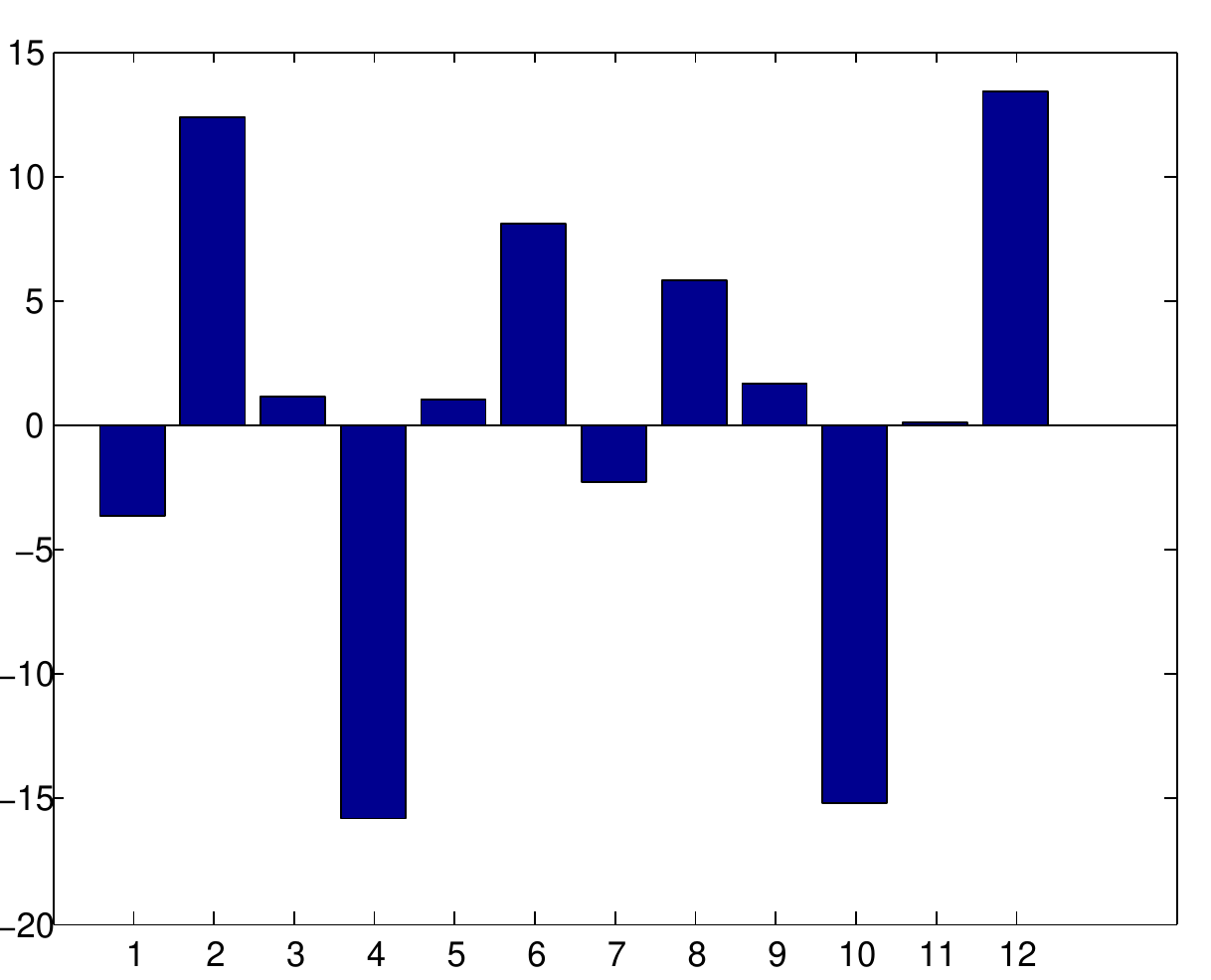} & \includegraphics[width=0.3\textwidth]{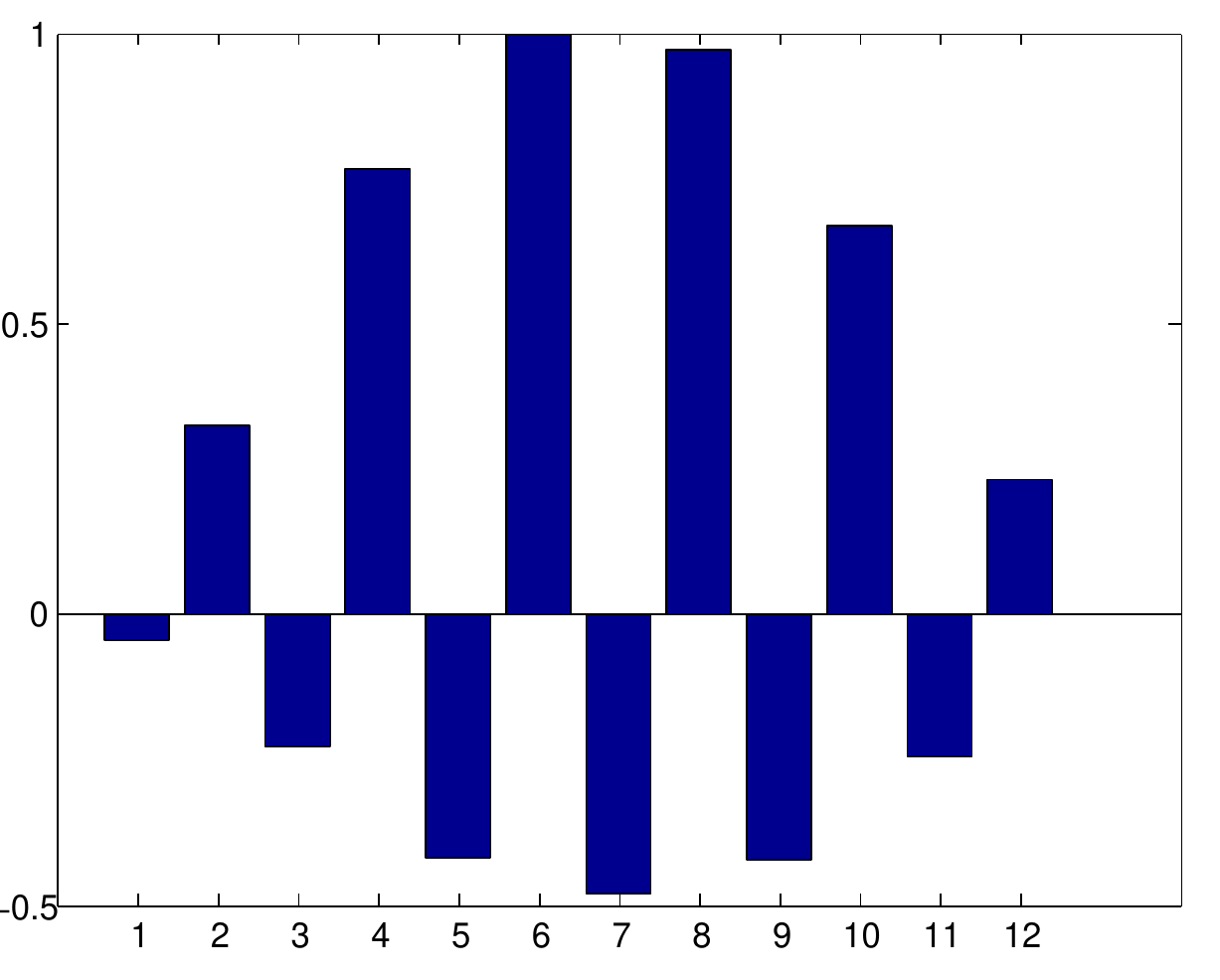} \\
    $\alpha_4 = 1, \; \xi = 1.1489$ & $\alpha_5 = 1, \; \xi =1.0778$ & $\alpha_6 = 1, \; \xi = 1.1584$ \\
    \includegraphics[width=0.3\textwidth]{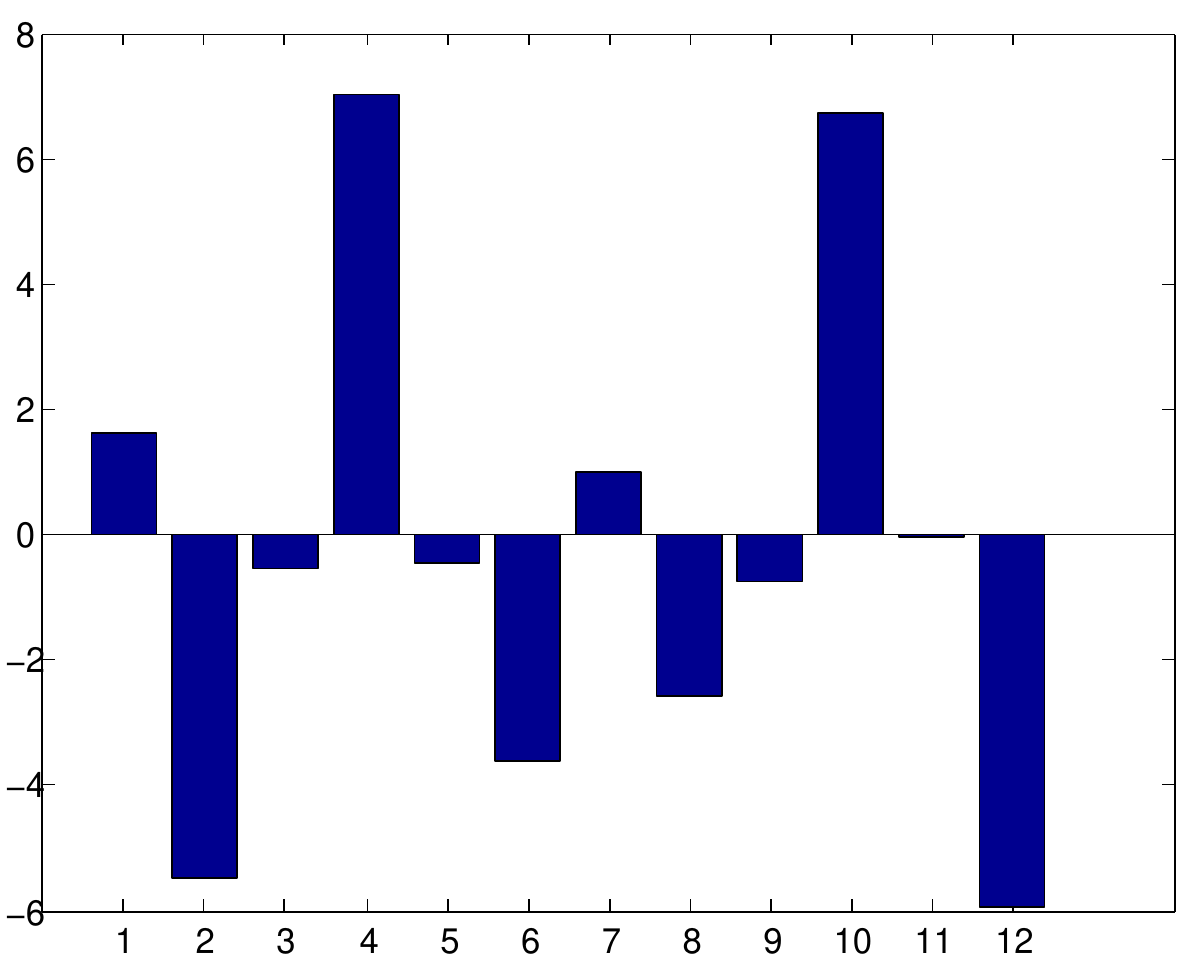} & \includegraphics[width=0.3\textwidth]{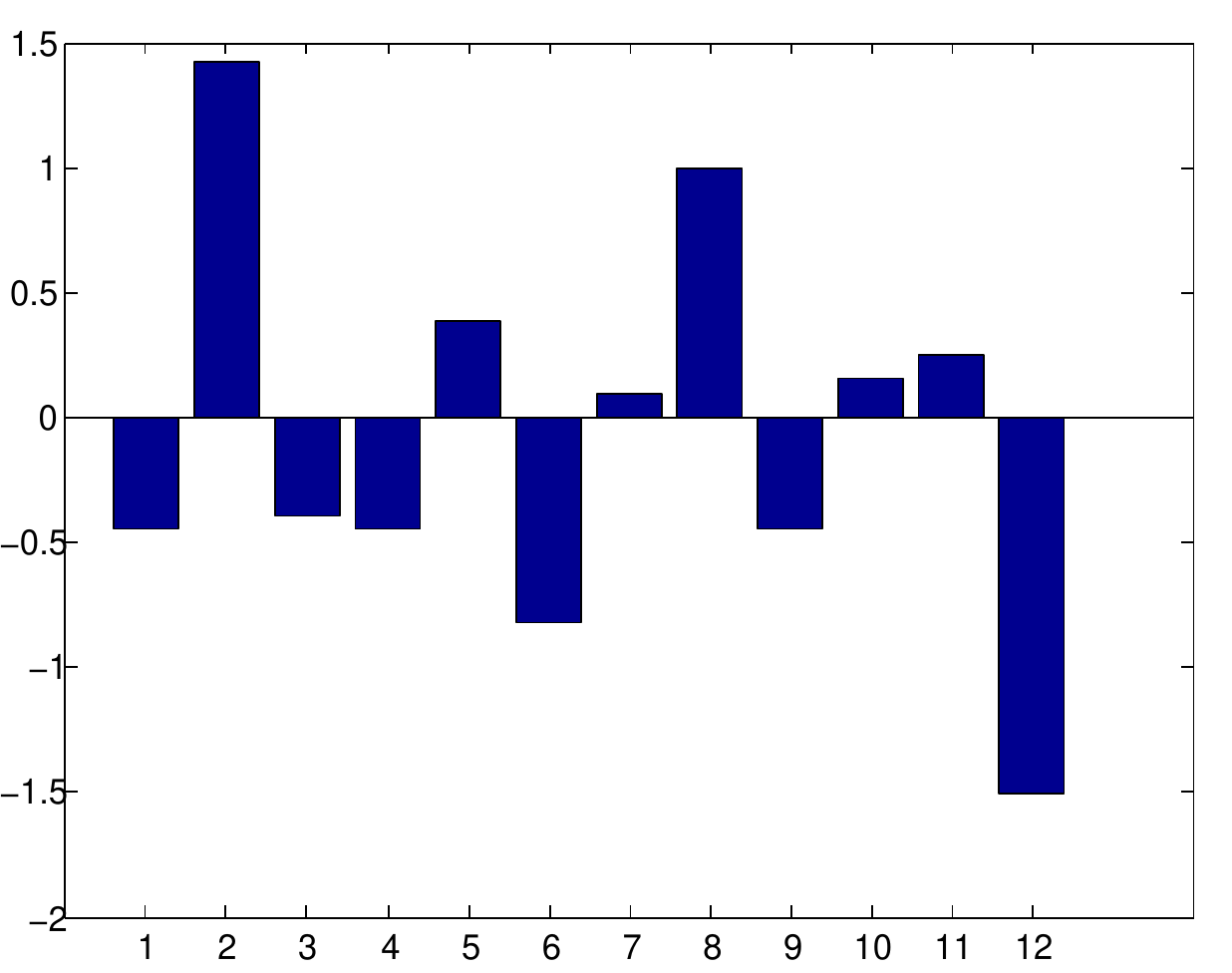} & \includegraphics[width=0.3\textwidth]{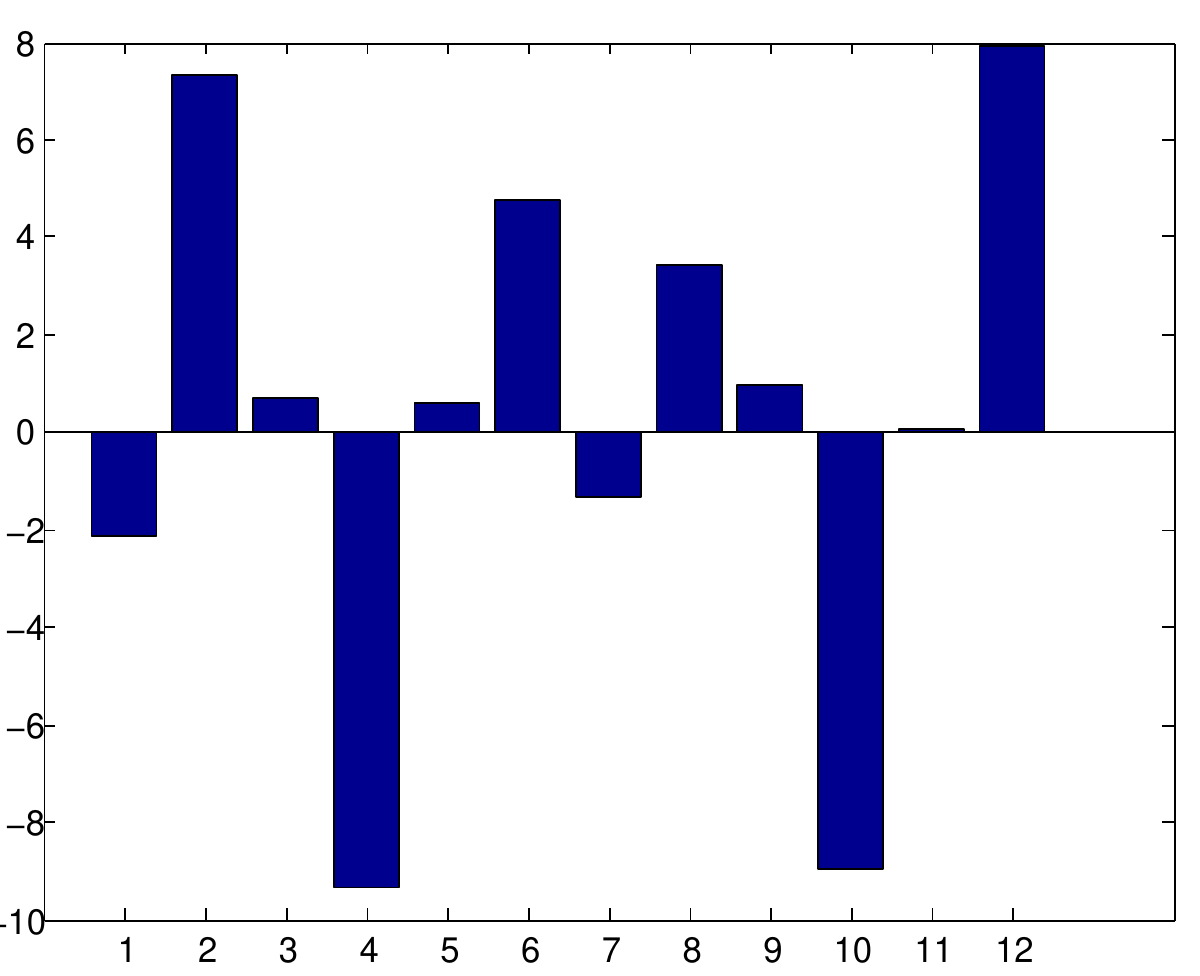} \\
    $\alpha_7 = 1, \; \xi =1.0778$  & $\alpha_8 = 1, \; \xi =1.1212$ & $\alpha_9 = 1, \; \xi =1.0778$ \\
    \includegraphics[width=0.3\textwidth]{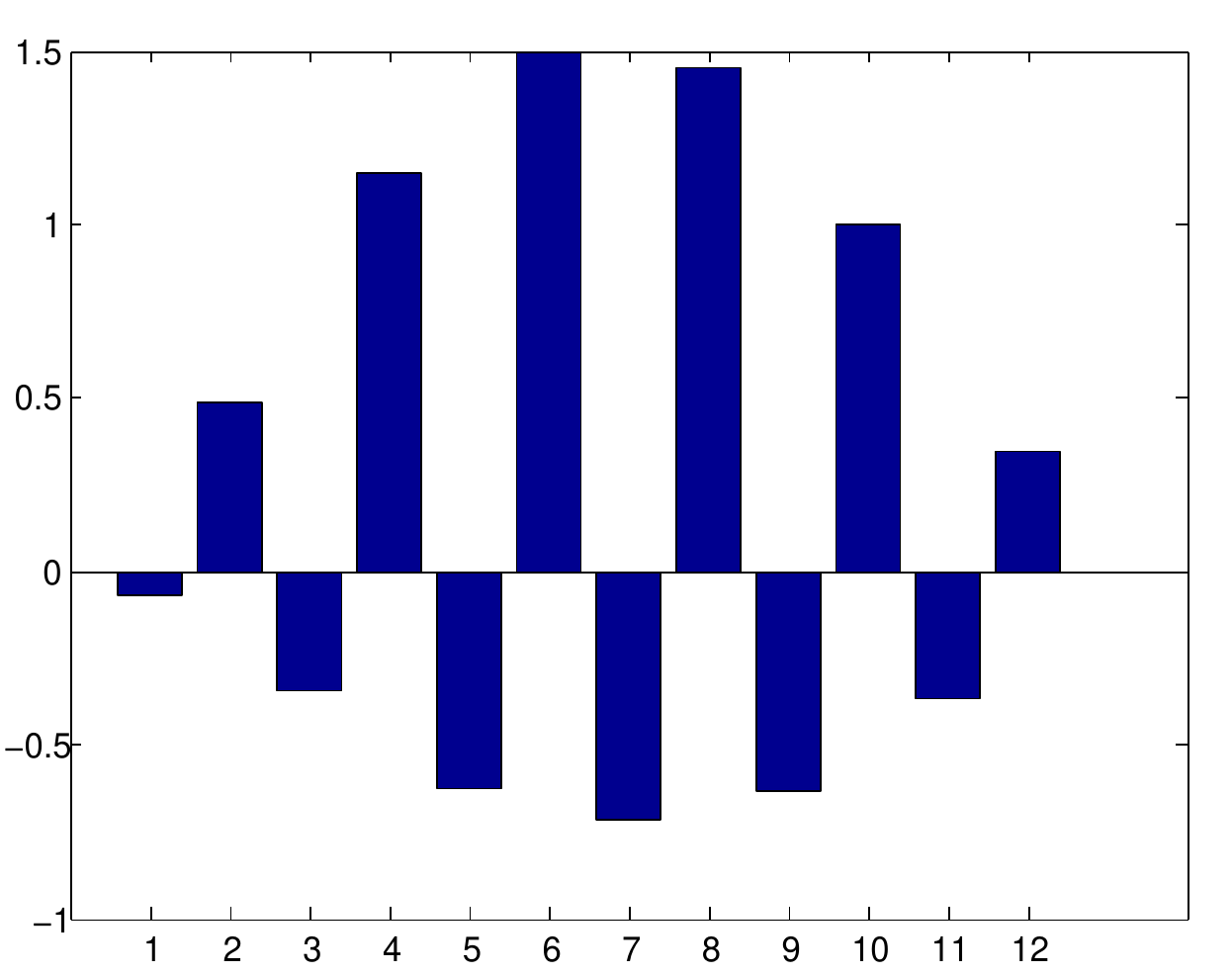} & \includegraphics[width=0.3\textwidth]{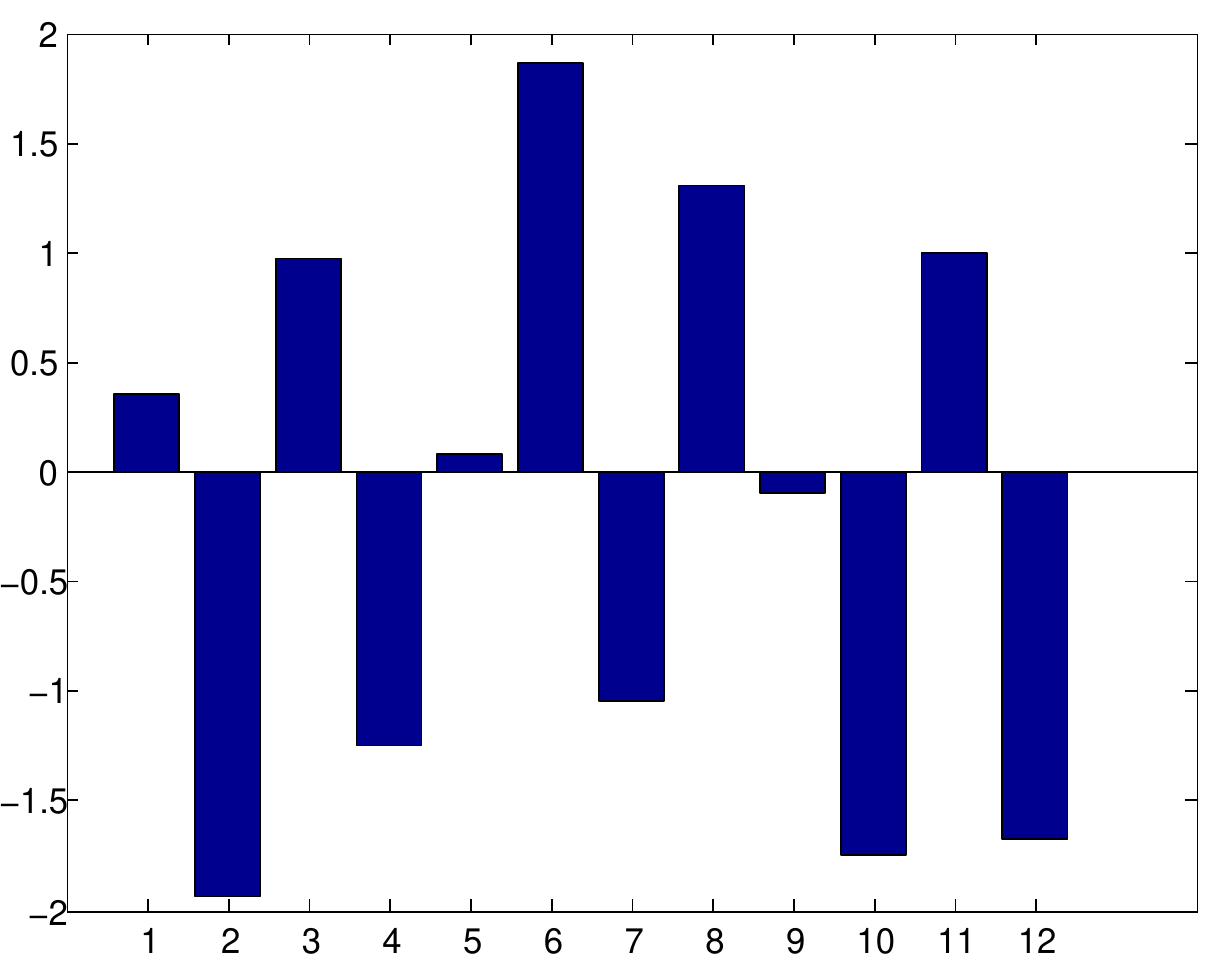} & \includegraphics[width=0.3\textwidth]{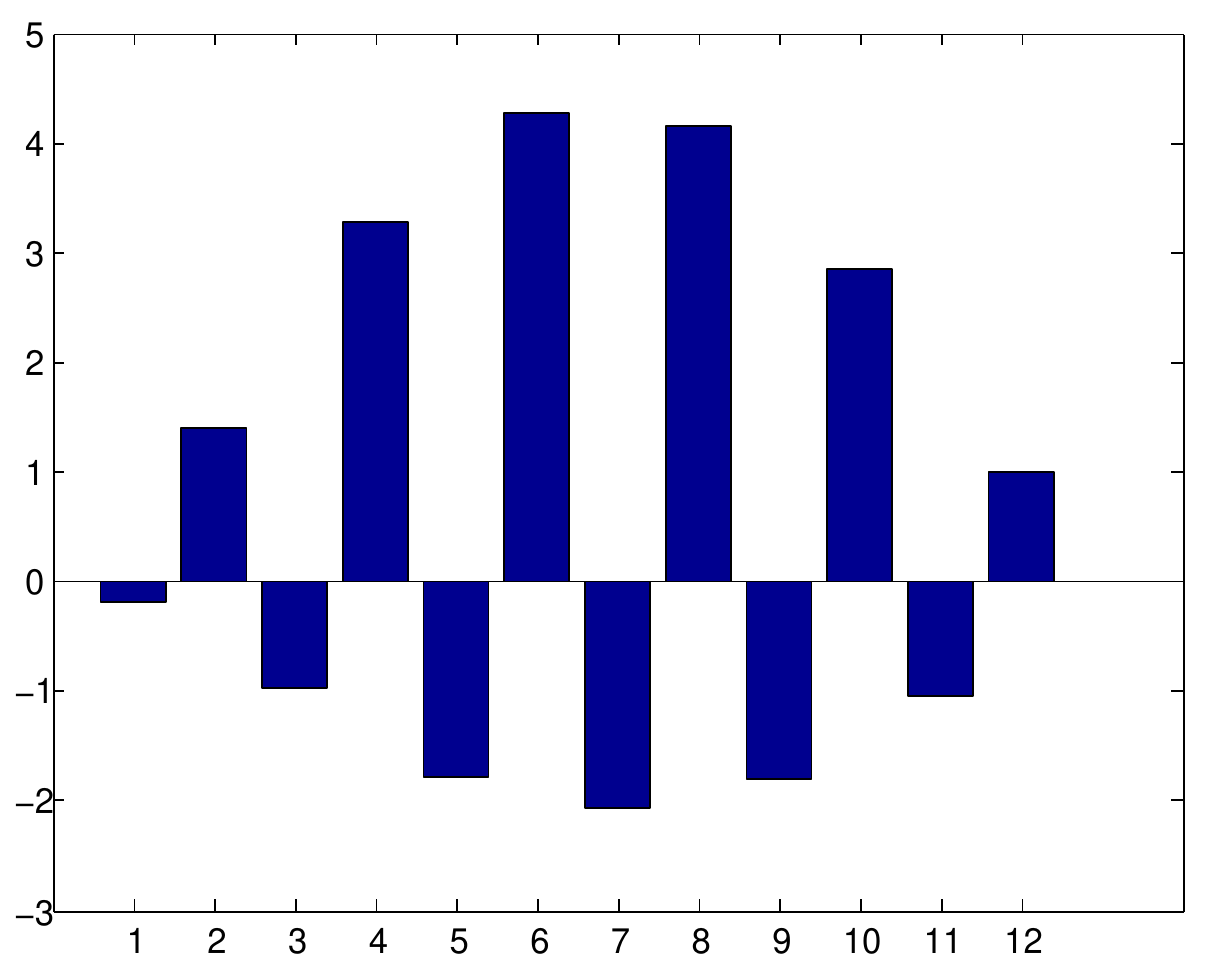} \\
    $\alpha_{10} = 1, \; \xi =1.1584$  & $\alpha_{11} = 1, \; \xi =1.1489$ & $\alpha_{12} = 1, \; \xi =1.1584$
  \end{tabular} }
\caption{City 1: number of cells is $N_b=6$. Number
of buildings per cell is $B=2$. $M=5$ ($2M$ is the number of computational points
per building). The foundation displacements $\alpha$ are depicted as bar graphs, $\xi$ are the coupling frequencies solving (\ref{nonlin_system}).
In each case  $\alpha_j =1$ is imposed, for the indicated $j$ -th building.
Note that runs that differ in the choice of the index $j$ for which 
we impose $\alpha_j=1$ may eventually lead to the same eigenvector (after rescaling).
For example, 
$\alpha_6=1$, $\alpha_{10}=1$, $\alpha_{12}=1$.}
\label{tab:city75_M5_N12}
\end{table}
   
\begin{table}[H]%\label{city2}
{\scriptsize
  \begin{tabular}{c c c}
    \includegraphics[scale=.3]{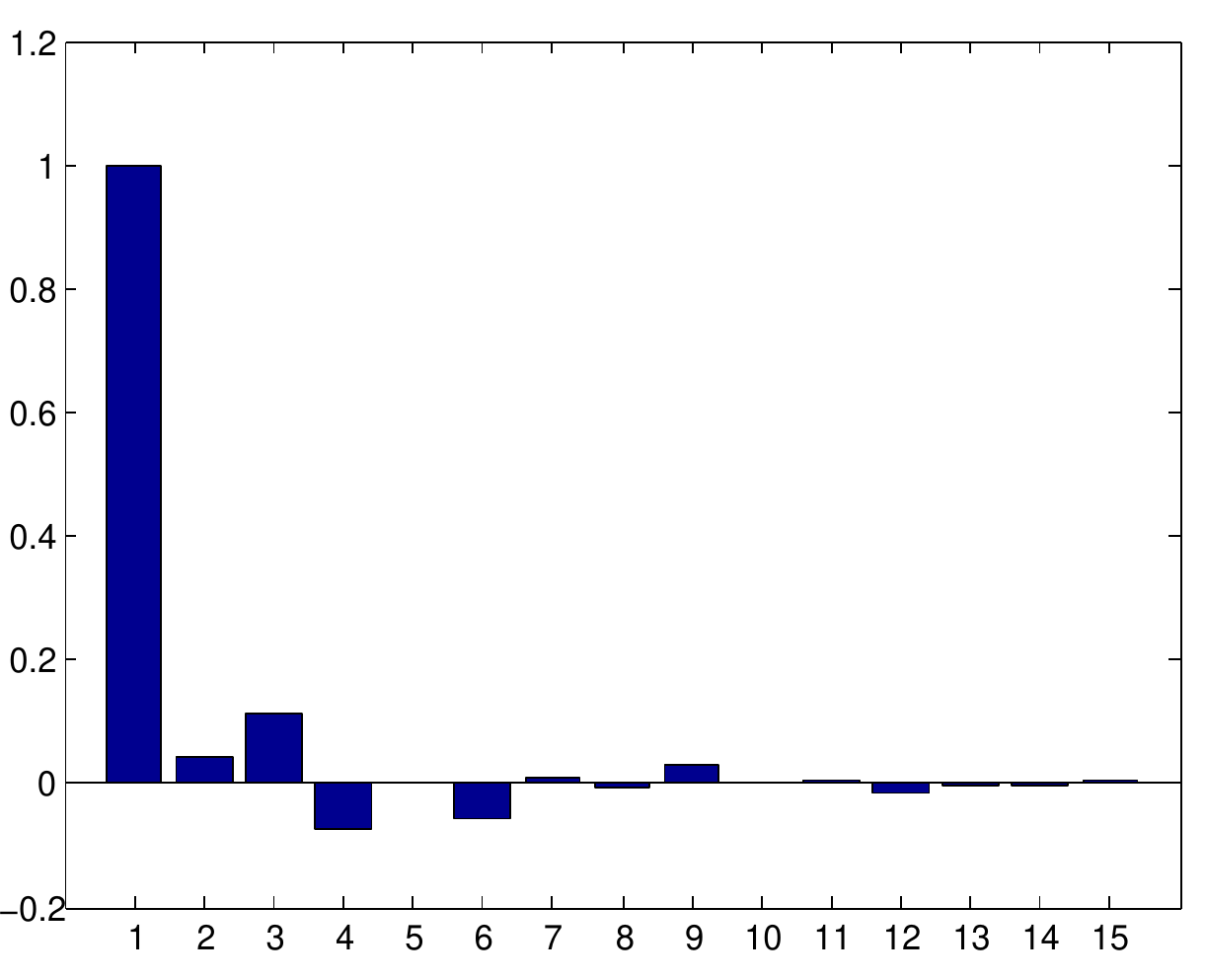} & \includegraphics[scale=.3]{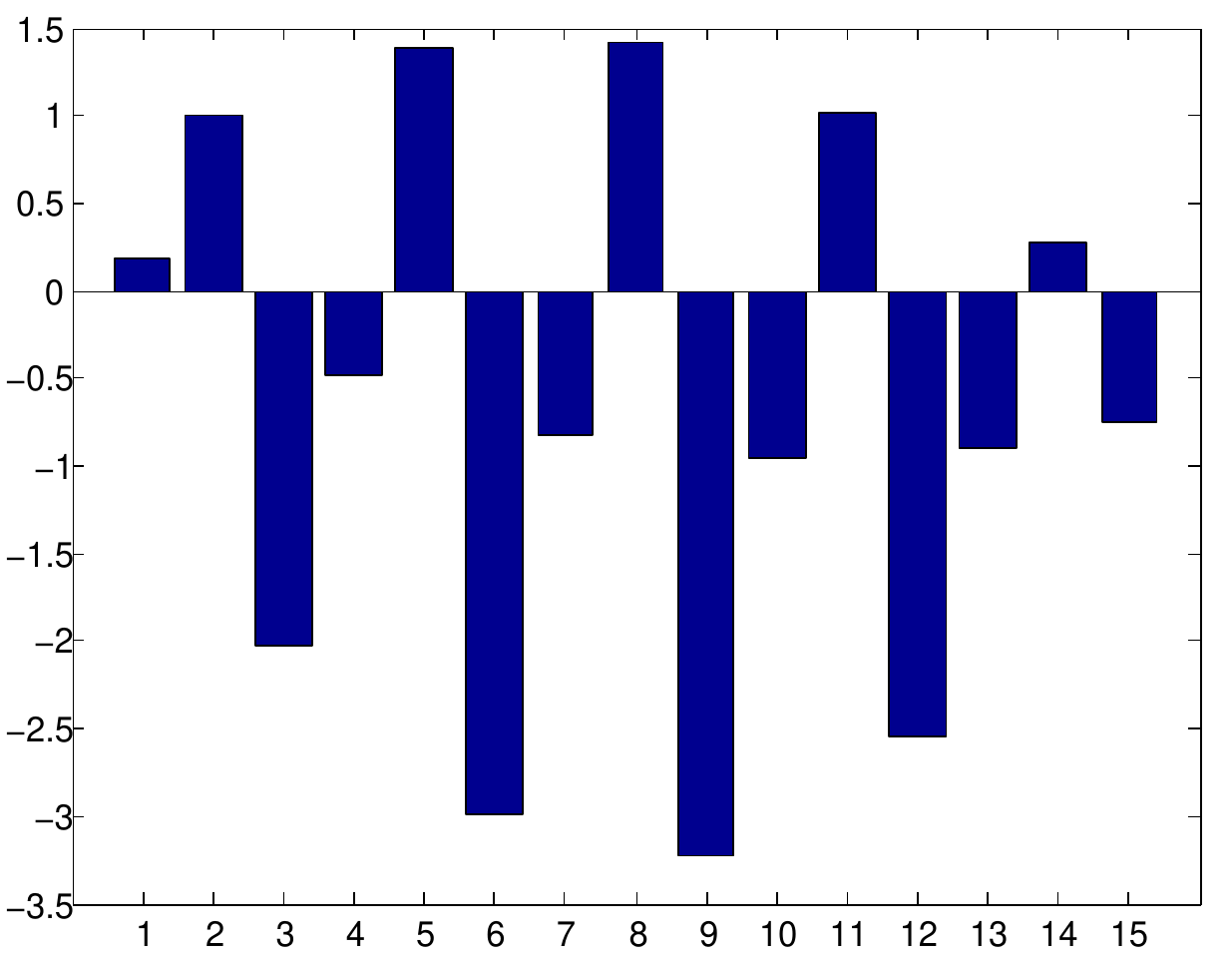} & 
		\includegraphics[scale=.3]{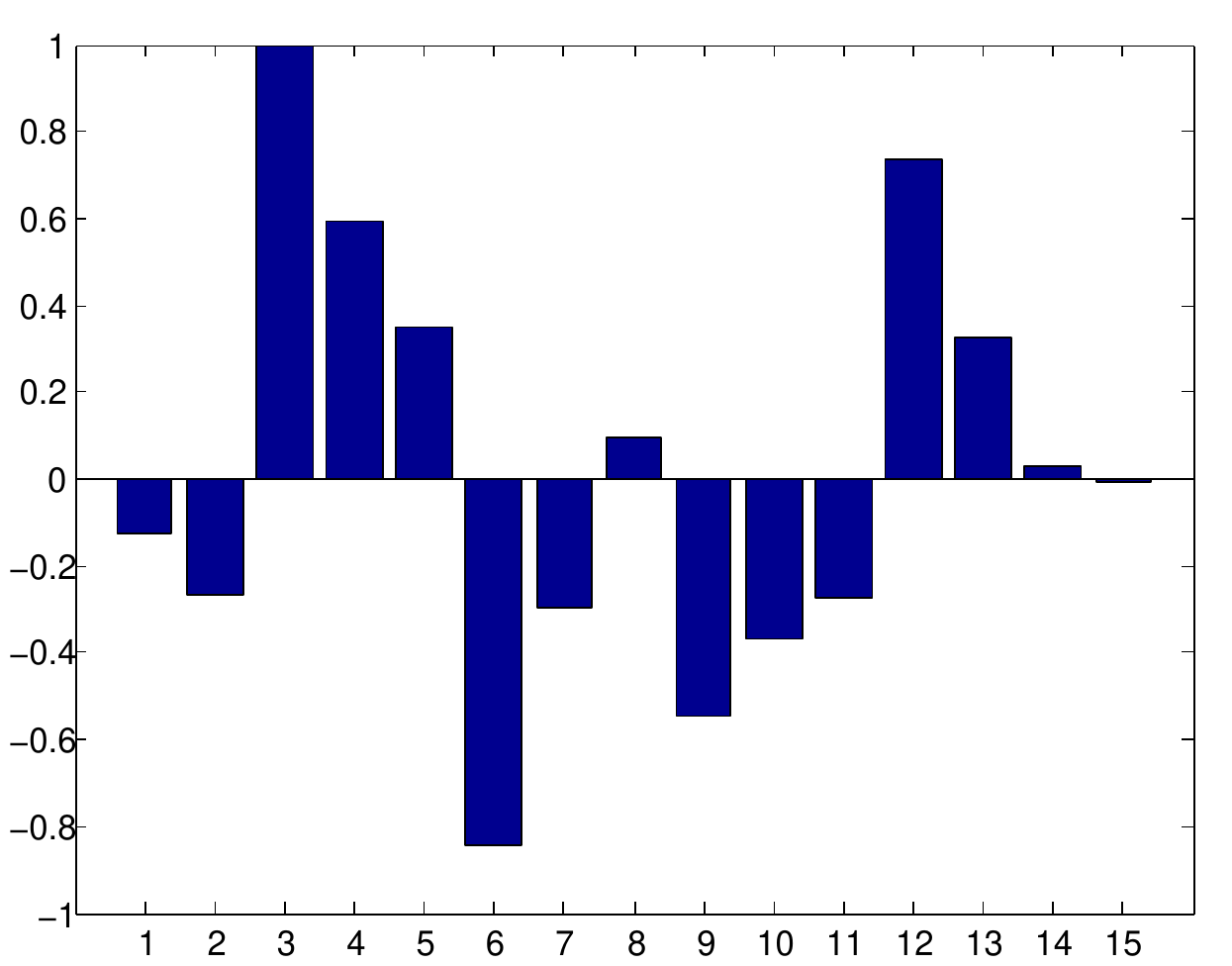} \\
    $\alpha_1 = 1, \; \xi = 1.4119$ & $\alpha_2 = 1, \; \xi =1.0349$ & $\alpha_3 = 1, \; \xi = 1.0043$ \\
    \includegraphics[scale=.3]{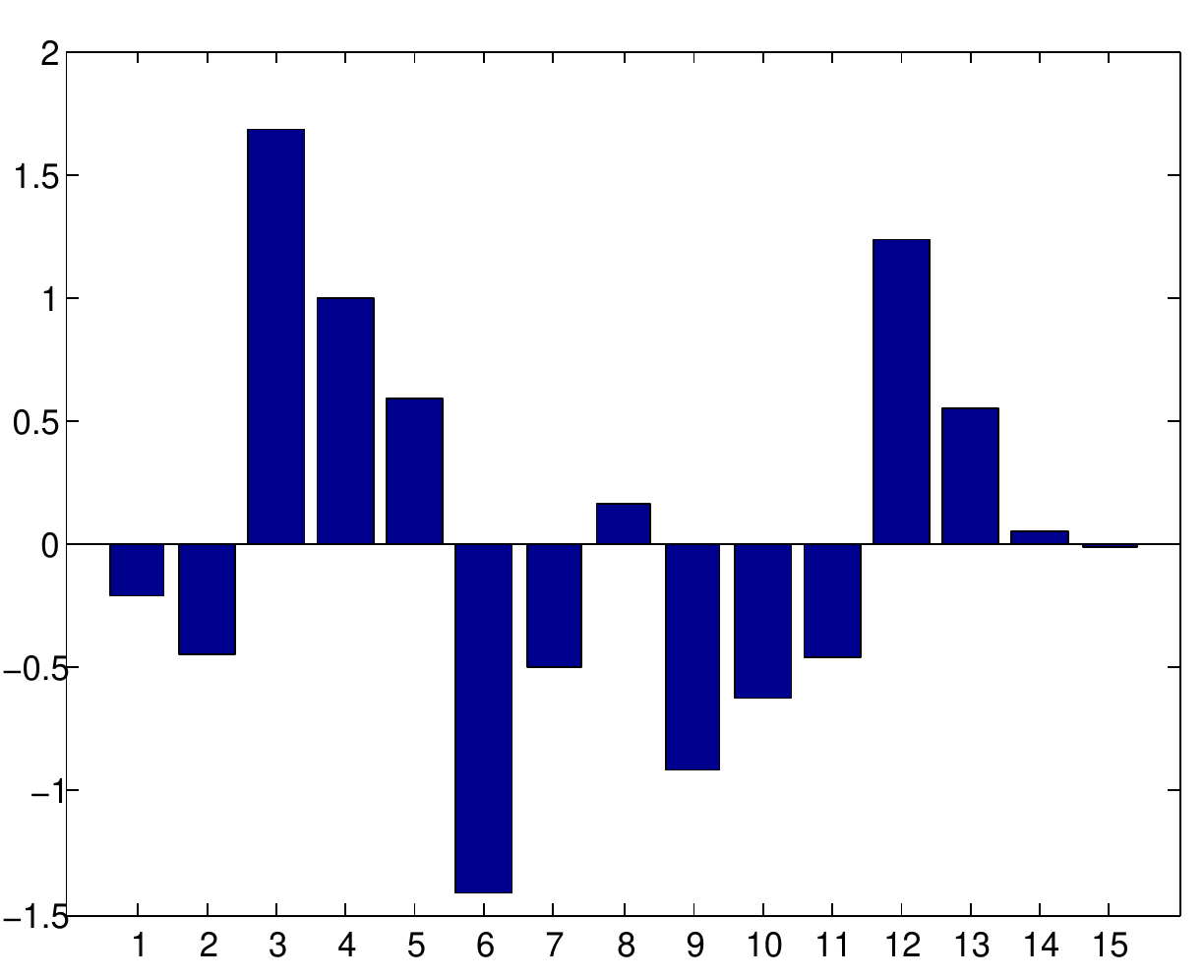} & \includegraphics[scale=.3]{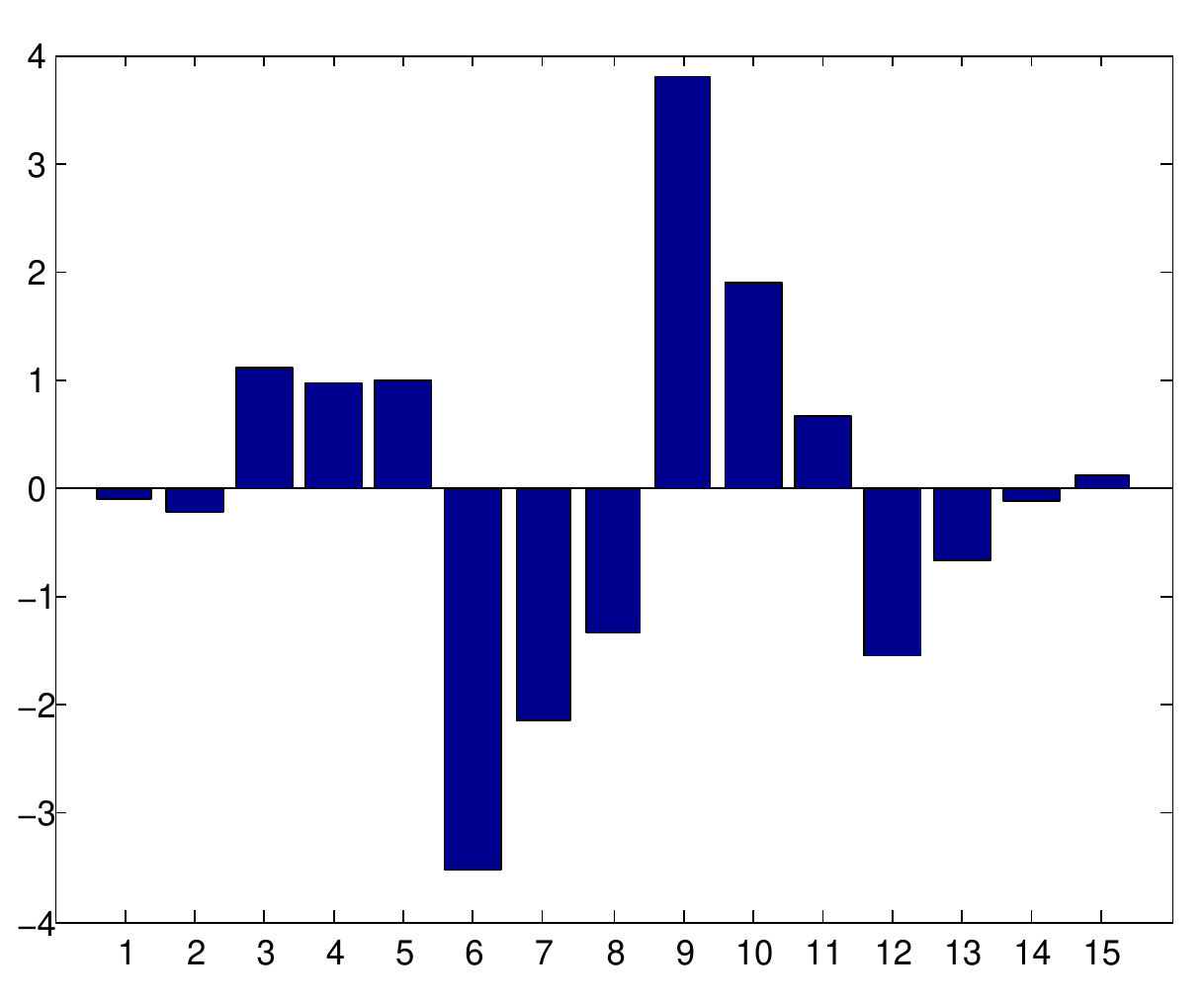} & 
		\includegraphics[scale=.3]{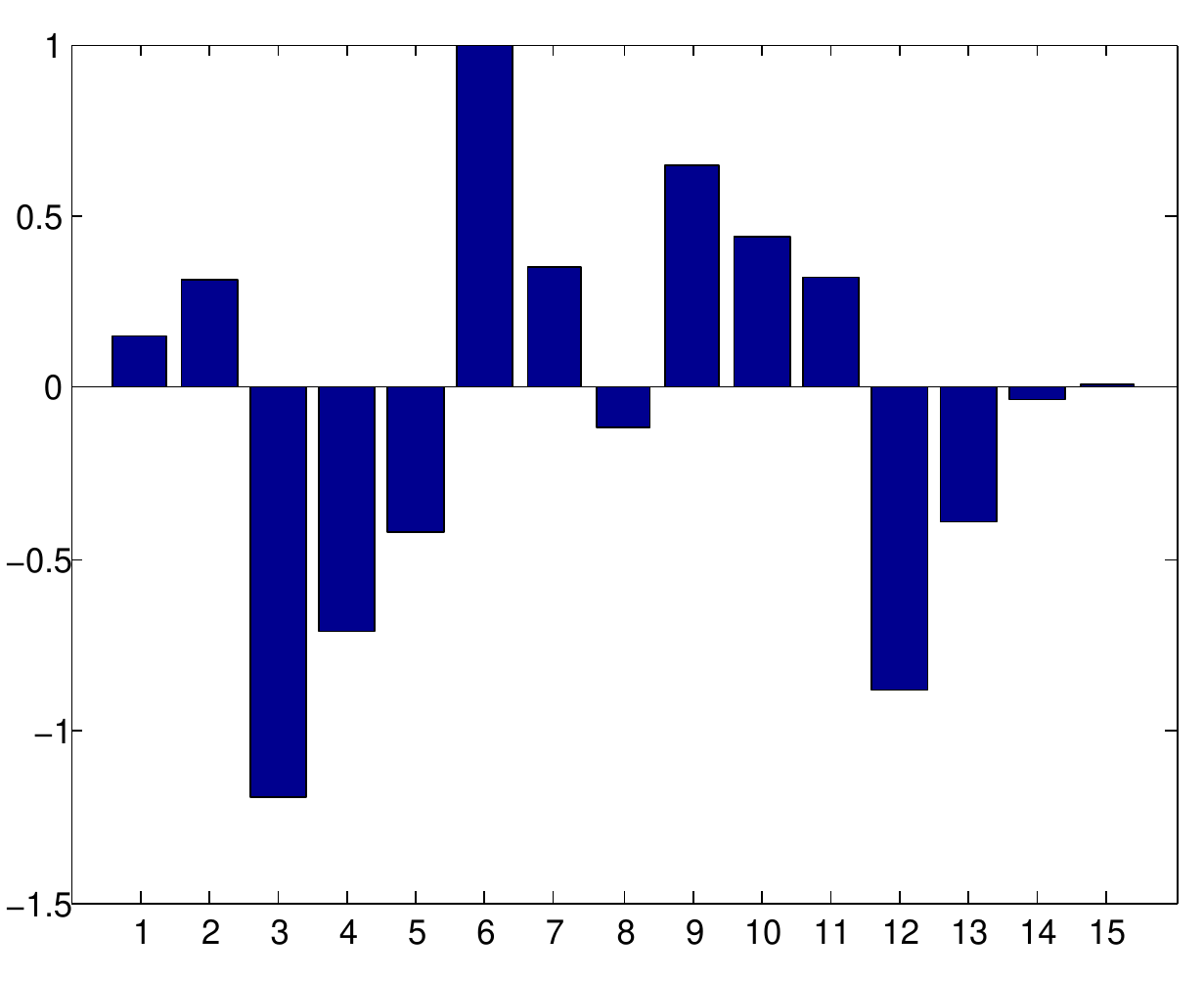} \\
    $\alpha_4 = 1, \; \xi = 1.0043$ & $\alpha_5 = 1, \; \xi =1.0048$ & $\alpha_6 = 1, \; \xi = 1.0043$ \\
    \includegraphics[scale=.3]{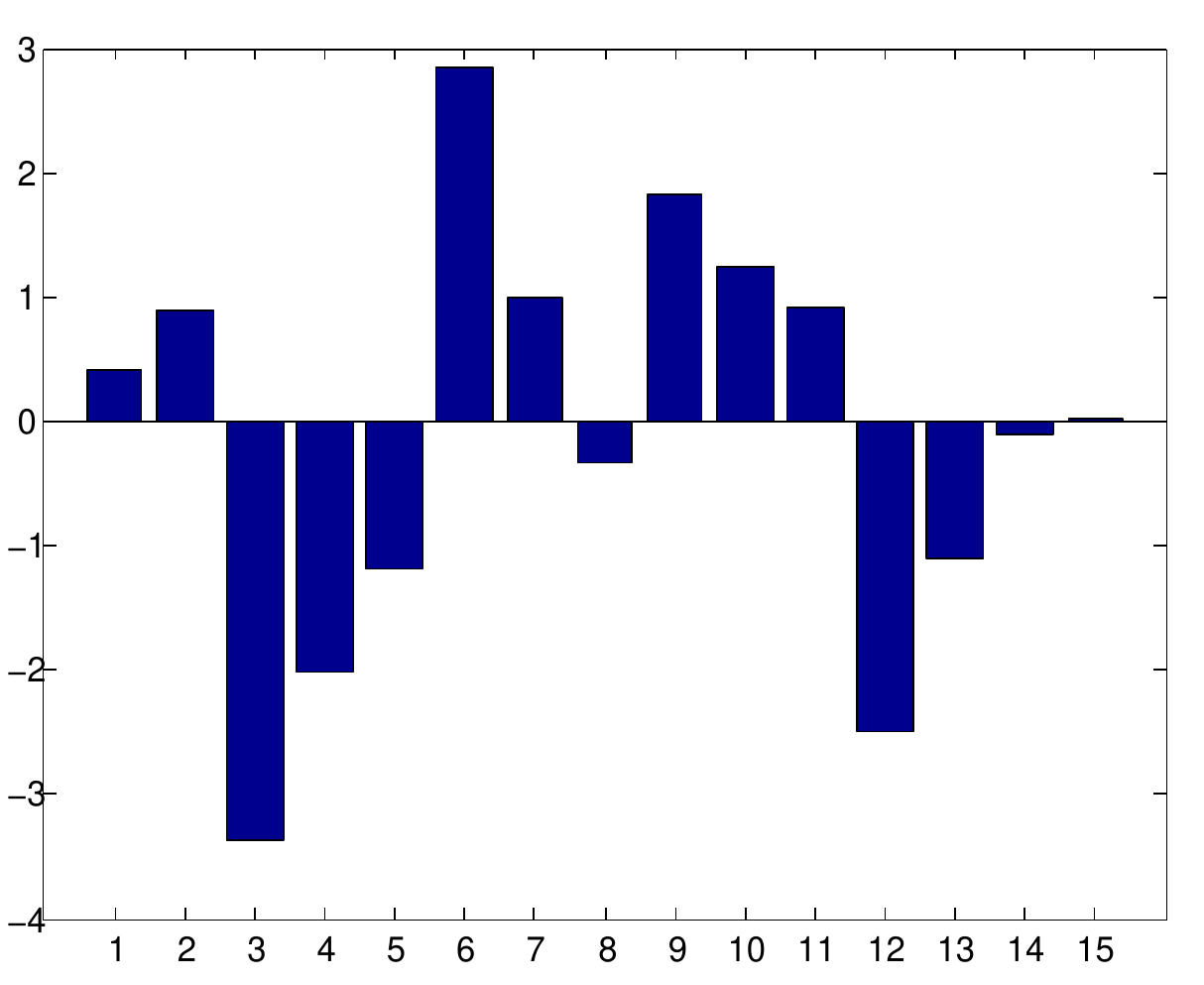} & \includegraphics[scale=.3]{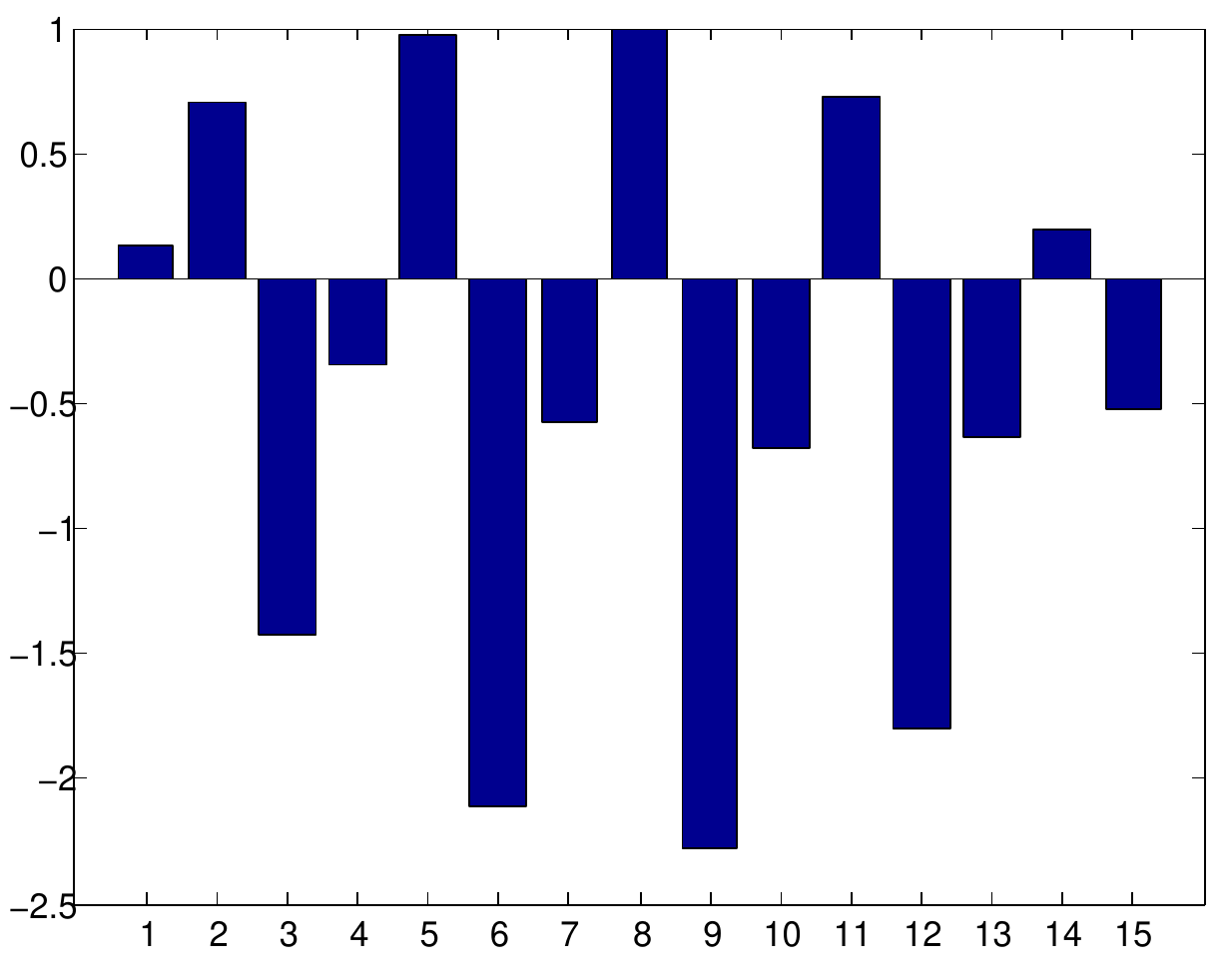} &
		\includegraphics[scale=.3]{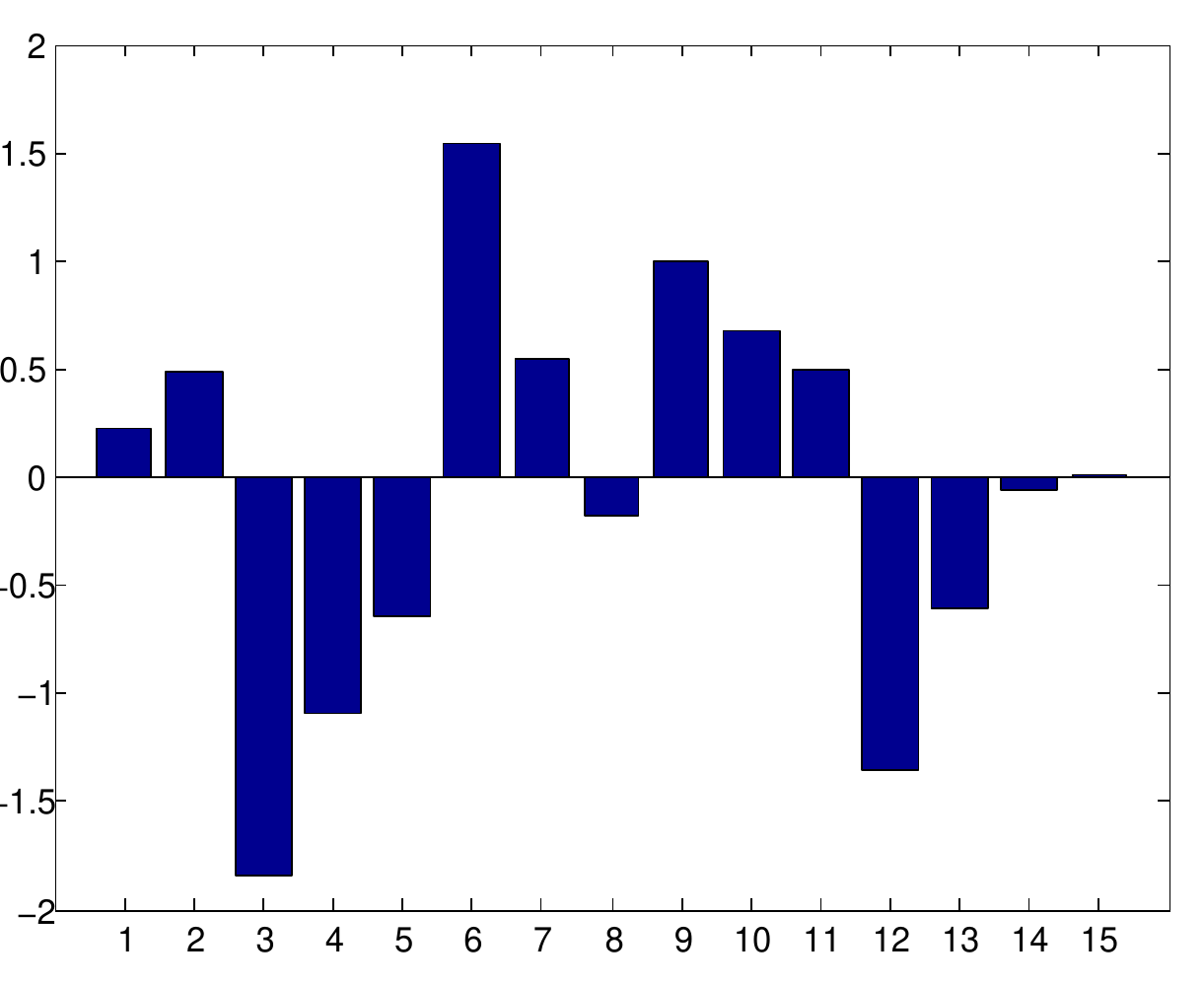} \\
    $\alpha_7 = 1, \; \xi =1.0043$  & $\alpha_8 = 1, \; \xi =1.0349$ & $\alpha_9 = 1, \; \xi =1.0043$ \\
    \includegraphics[scale=.3]{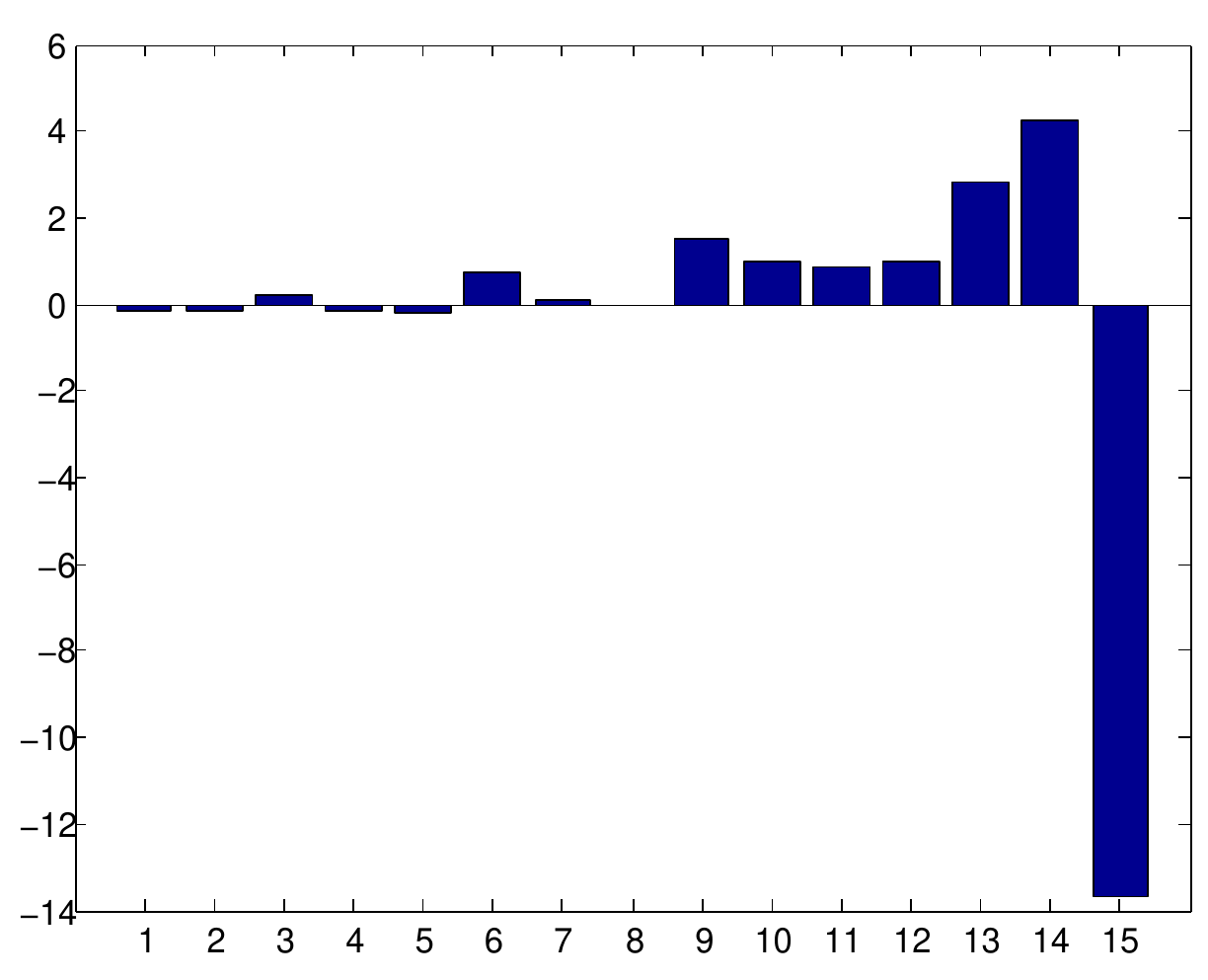} & \includegraphics[scale=.3]{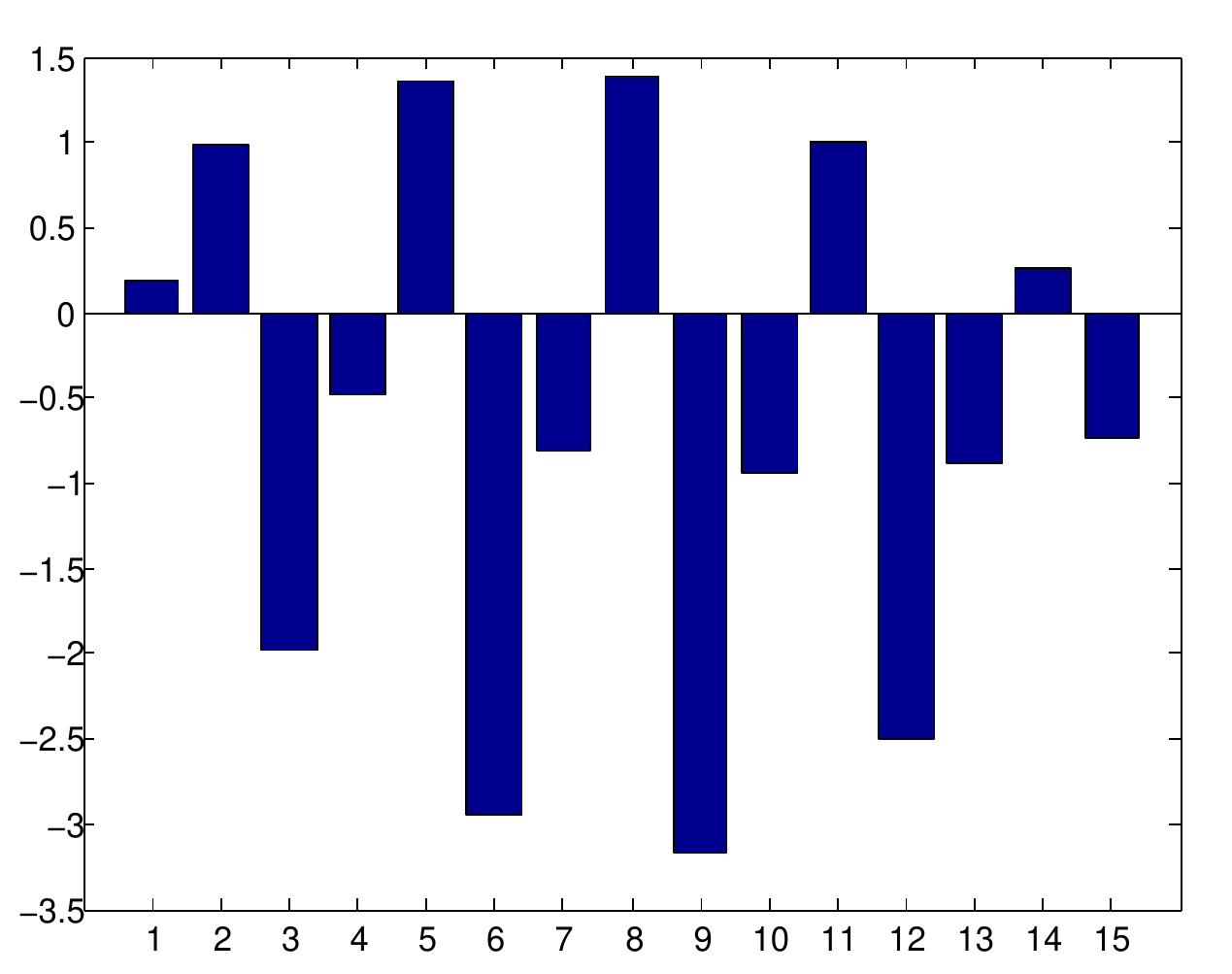} & \includegraphics[scale=.3]{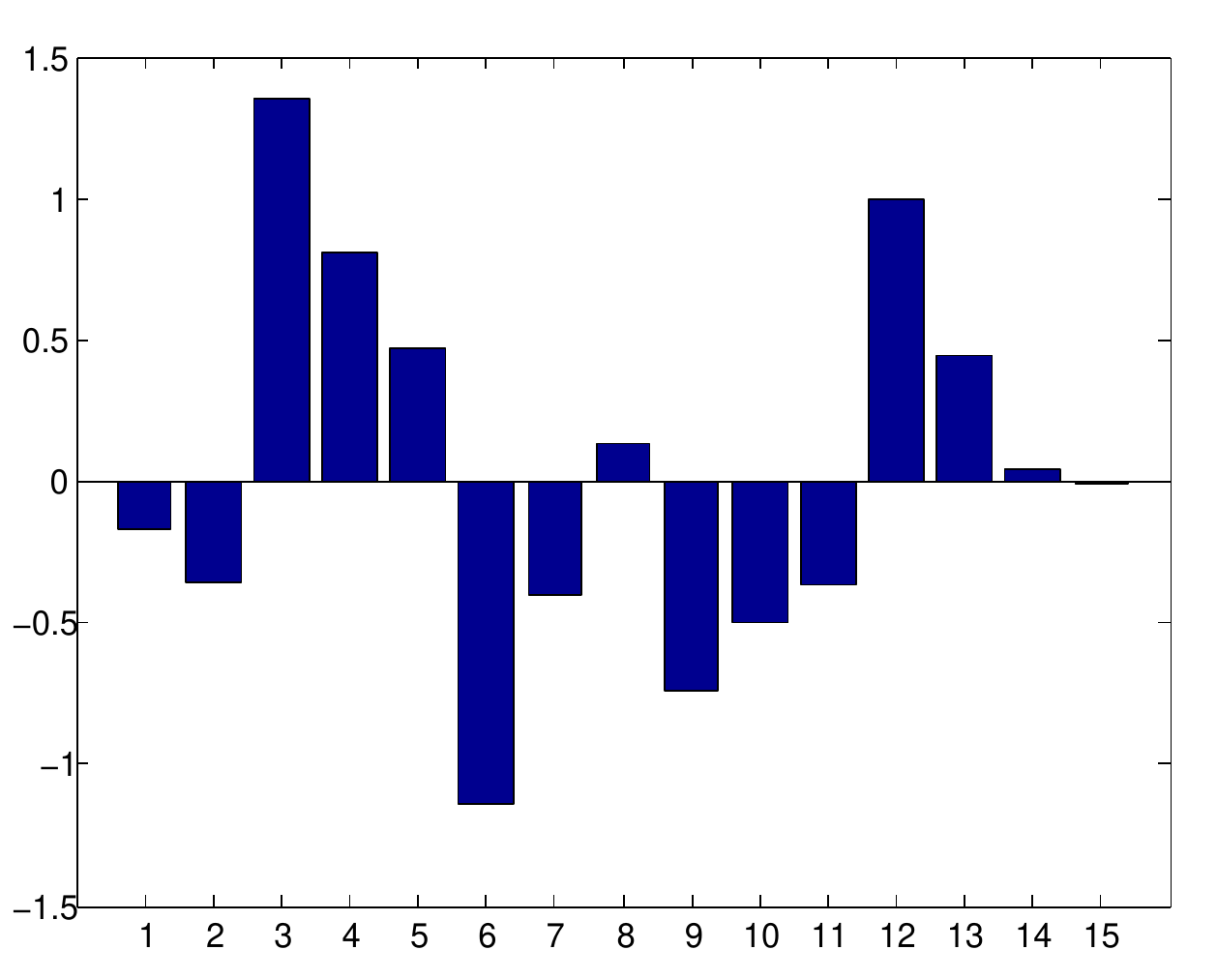} \\
    $\alpha_{10} = 1, \; \xi =0.9837$  & $\alpha_{11} = 1, \; \xi =1.0349$ & $\alpha_{12} = 1, \; \xi =1.0043$ \\
    \includegraphics[scale=.3]{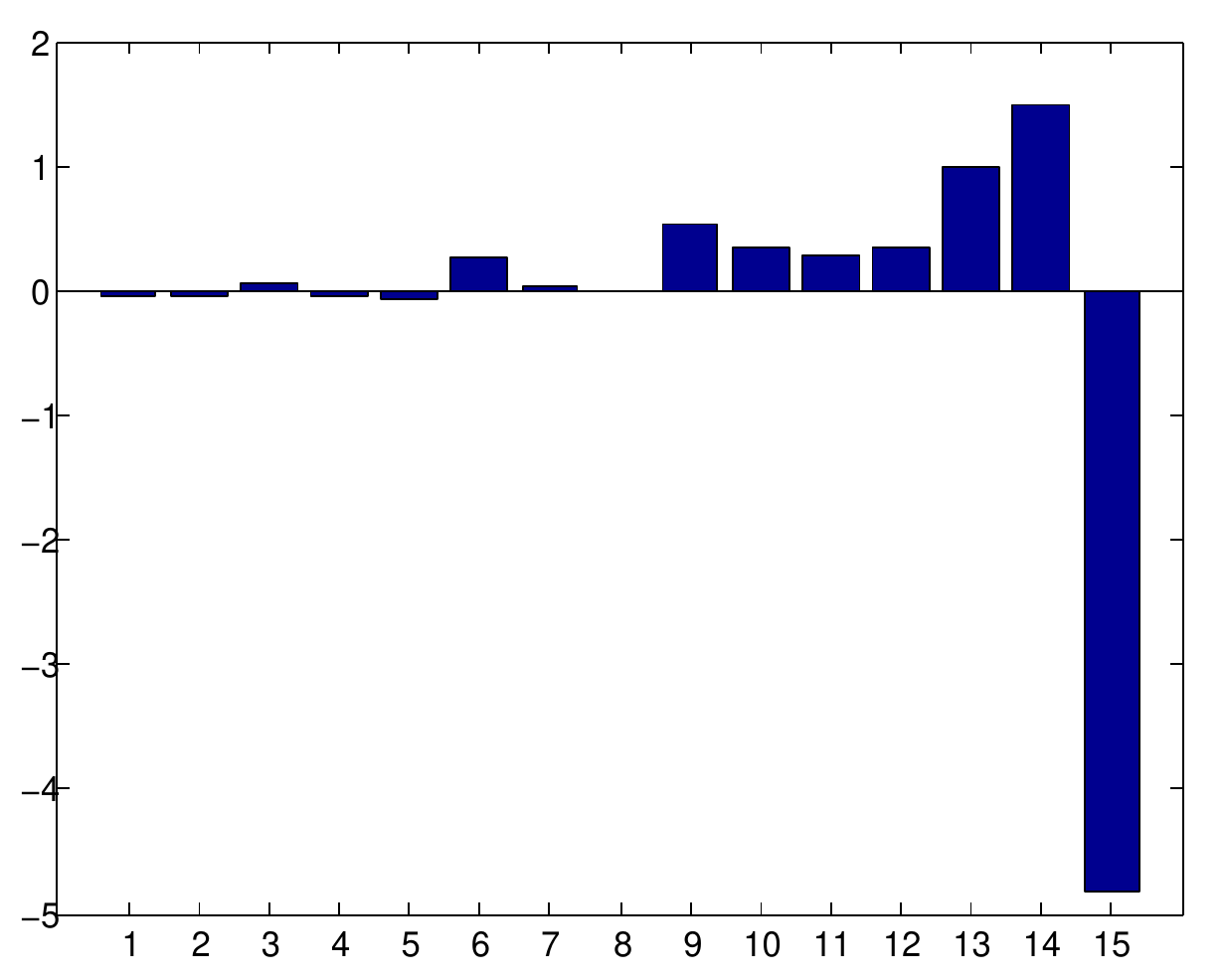} & \includegraphics[scale=.3]{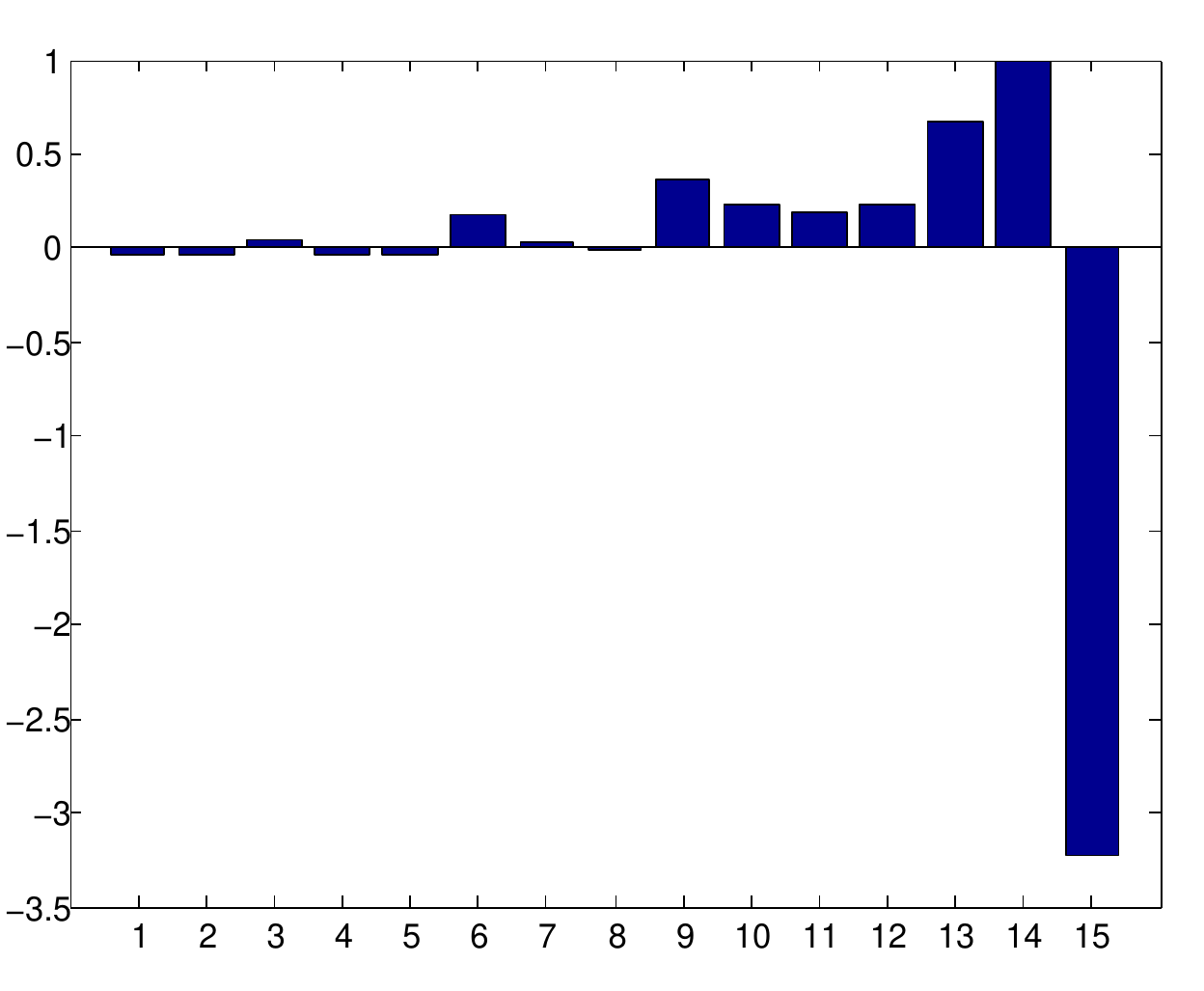} & \includegraphics[scale=.3]{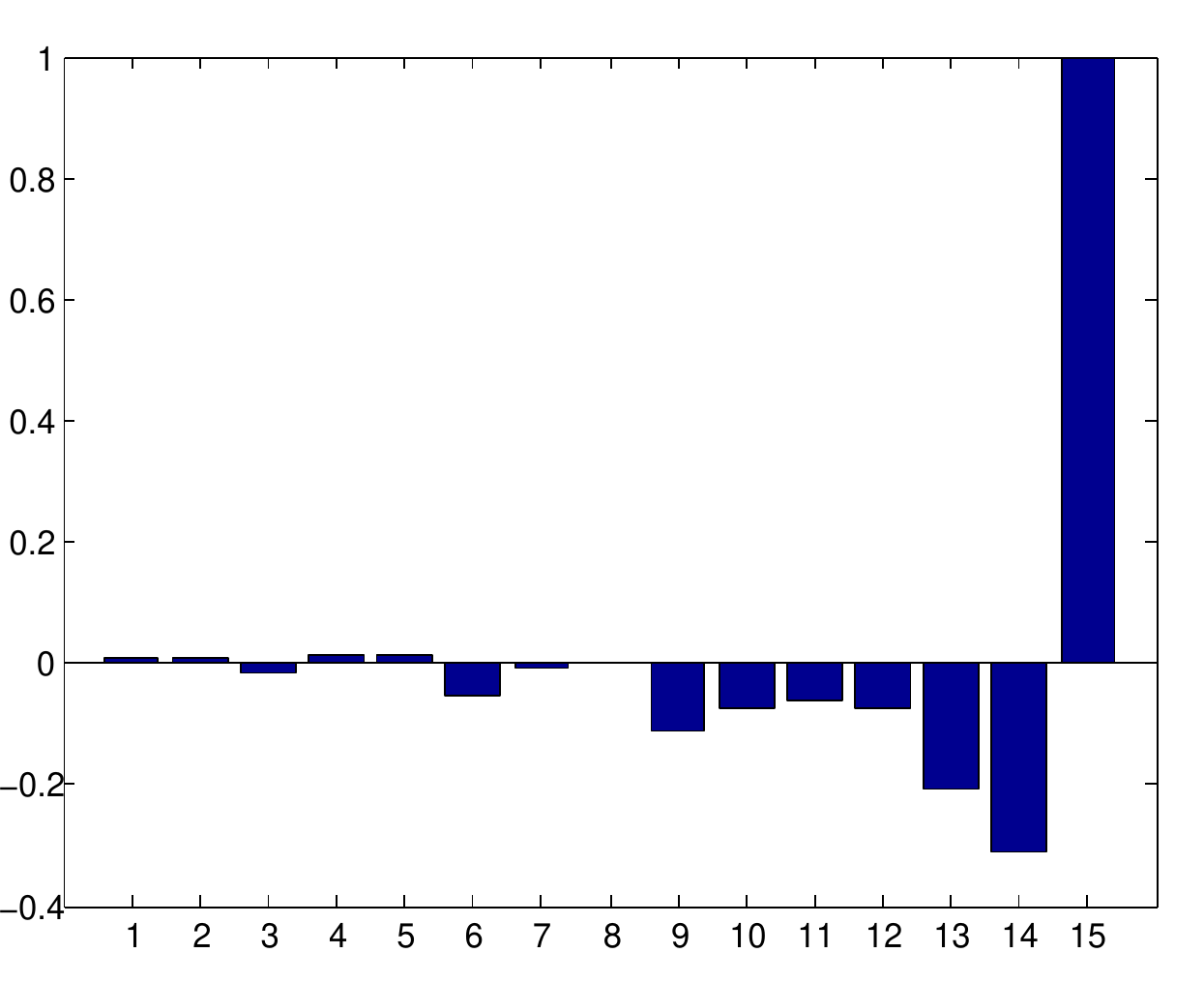} \\
    $\alpha_{13} = 1, \; \xi =0.9837$  & $\alpha_{14} = 1, \; \xi =0.9837$ & $\alpha_{15} = 1, \; \xi = 0.9837$
  \end{tabular} }
\caption{City 2: number of cells is $N_b=5$.  Number
of buildings per cell is $B=3$. $M=5$ ($2M$ is the number of computational points
per building). The foundation displacements $\alpha$ are depicted as bar graphs, $\xi$ are the coupling frequencies solving (\ref{nonlin_system}).
In each case  $\alpha_j =1$ is imposed, for the indicated $j$ -th building.}
\label{tab:city7_M5_N15}
\end{table}

%--------------------------------------------------------------------------------------------------- 
\subsection{Case of repeated patterns of buildings}
Assume that the number of buildings $N$ is now infinite and that there is a periodic pattern of buildings.
We will denote by $B$  the number of buildings in the pattern (as previously), and by $2P$ the length of each periodic cell.
Accordingly,
\begin{itemize}
\item all the physical parameters of the buildings $m_{1,j}$, $m_{0,j}$, 
$l_j$, $h_j$, $\rho_j$, $\beta_j$ are periodic in $j$, with period $B$
\item $a_{j+B} = a_j + 2P $, $b_{j+B} = b_j + 2P $
\end{itemize}
Assume that the foundations of the buildings $1, ..., B $ are included in some interval
$[I_1, I_2]$, where $I_2 - I_1 =2 P$. 
Introduce the following notations:
$\Omega_{per} = (I_1, I_2) \times (0, \infty)
$,  $\Gamma_{per,j} =  [a_j, b_j]$,
$\Gamma_{per}^{free} = (I_1, I_2) \setminus \ds \cup_{j=1}^B \Gamma_{per,j}$.
In this new case %finding coupling frequencies for the underground and buildings system
%involves solving the periodic PDE 
the analog of (\ref{helmholtz_per_1}-\ref{per_decay}) 
is
\begin{align}
     \label{helmholtz_per_1_mult} \Delta \Psi + \xi^2 \Psi &= 0 \mbox{ in } \Omega_{per},
  \\ \label{helmholtz_per_2_mult} \displaystyle \Psi = \alpha_j \mbox{ on } \Gamma_{per,j}, 
  j=1, ..., N
  & \; \; \frac{\partial \Psi} {\partial y} = 0 \mbox{ on } \Gamma_{per}^{free}.
\end{align}
augmented by the decay condition (\ref{per_decay}).
The corresponding periodic boundary condition is
\be \label{per_cond2}
\Psi(I_1, y) = \Psi(I_2, y)  \mbox{ for all } y \geq 0
\ee
To
find coupling frequencies for the underground and buildings system
%involves solving the periodic PDE
we have to solve the  system of $B$ non linear equations
\be
  \label{eigenvalue_per_sys}
 \alpha_j  \ds q_j(\xi^2) = p_j(\xi^2) \mbox{Re }\int_{\Gamma_{per,j}} {\frac{\partial \Psi} {\partial y} (s,0)} ds,
  \q j=1, ..., B
\ee
where these equations are coupled through the PDE
(\ref{per_decay}, \ref{helmholtz_per_1_mult}-\ref{per_cond2}).
As previously, one of the $\alpha_j$'s in (\ref{helmholtz_per_2_mult}) may be set to 1, while the
others will have to be determined.\\
Let us now examine a numerical example. The patterns of buildings are similar to those
from the previous section, to facilitate comparison.
We discuss two cases, City 1-per and City 2-per (where "per" is for periodic):
\begin{itemize}
\item 
City 1-per: $B=2, P=7.5$,
$(a_1, a_2)=(-2.5, 1.5), \q(b_1, b_2)=(-1.5, 3)$,
\item  City 2-per: $B=3, P=7$,
$(a_1, a_2,a_3)=(0,2,5), \q(b_1, b_2)=(1.2,3,6.7)$. 
\end{itemize}
Accordingly the geometry of each periodic cell of these two cities is the same 
as the ones sketched in Figures 
 \ref{fig:city75} and \ref{fig:city7}.
 As previously, we pick $\gamma_1 =\gamma_j= 1.5$, $f_1 = f_j= 0.5$,  $r_j=r_1=0.1$, ${\cal{B}}_j={\cal{B}}_1=1.5$,
 but the lengths of the building foundations, $b_j -a_j$ are variable.
The table below 
gives computed values of coupling frequencies
$\xi$, given an initial search value $\xi_0$. We varied $2M$, the number
of grid points on each building to illustrate numerical convergence.

\begin{enumerate}
  \item City 1-per \\
  $\xi_0 = 1, \; M=5: \; \xi = 1.1594, \; \alpha = (1, -2.1171)$, \\
  $\xi_0 = 1, \; M=10: \; \xi = 1.1583, \; \alpha = (1, -2.1222)$, \\
  $\xi_0 = 1, \; M=20: \; \xi = 1.1580, \; \alpha = (1, -2.1241)$.
  \item City 2-per \\
  $\xi_0 = 1, \; M=5: \; \xi = 1.0420, \; \alpha = (1, -1.5703, 3.5458)$, \\
  $\xi_0 = 1, \; M=10: \; \xi = 1.0382, \; \alpha = (1, -1.5103, 3.4610)$, \\
  $\xi_0 = 1, \; M=20: \; \xi = 1.0368, \; \alpha = (1, -1.4874, 3.4288)$.
\end{enumerate}

% We notice that some of the free-space solutions go to the periodic solutions as the number of buildings grows, just as in the case of a homogeneous city, wavenumbers $\xi$ almost coincide. Tables \ref{tab:asymptot_for_city7.5} and \ref{tab:asymptot_for_city7} illustrate wavenumber convergence, showing one of the free-space solutions for different number of clusters in the city7.5 and city7. 
Next we report that given a pattern of $B$ buildings, repeated $N_c$ times,  
some computed coupling frequencies exhibit a convergence trend as $N_c $ grows large, 
and the limit value equals a  
coupling frequency for the periodic problem. 
Let us examine the case of the geometry 
given by City 1.
In the periodic case City 1-per we found for $M=5$ and the initial values
$\xi_0 = 1, \, \alpha_1 = 1$, the final values
$ \xi_{per} = 1.1594, \, \alpha_2 =-2.1222 $. 
This is clearly close to the two cases shown in Table \ref{tab:city75_M5_N12}
for, say,  $\alpha_6 =1$. 
 If we keep increasing the number of clusters $N_c$ we find at least one solution $\xi$ which approaches 
$\xi_{per} $, as shown in Table \ref{tab:asymptot_for_city7.5}. 
We also notice convergence, in some sense, of the parameter $\alpha$.
In Table \ref{periodic_outcome}, we sketched the computed value  $(\alpha_1, \alpha_2)=(1, -2.1171)$
for City 1-per (when the initial guess for $\xi$ is $\xi_0=1$),
which is clearly close to a multiple of $(\alpha_7, \alpha_8)$ 
obtained in  the case  $\alpha_6 =1$ in Table \ref{tab:city75_M5_N12}. \\
 Similar observations can be made in the case of the geometry of City 2: see 
Table \ref{tab:asymptot_for_city7}.
This time $(\alpha_1, \alpha_2, \alpha_3) $ computed in the periodic case, given in Table
\ref{periodic_outcome} compares to $(\alpha_7, \alpha_8, \alpha_9) $ 
in Table \ref{tab:city7_M5_N15}, case $\alpha_8=1$.

\begin{table}[H]
\centering
  \begin{tabular}{|c|c|c|c|c|}
    \hline
    & $N_c=4$ & $N_c=5$ & $N_c=6$ & Periodic \\
    \hline
    $\xi$ & 1.1572 & 1.1579 & 1.1584 & 1.1594 \\
    \hline
  \end{tabular}
  \caption{Convergence to $\xi_{per}$, for the geometry given by City 1.
 $M=5$.}
  \label{tab:asymptot_for_city7.5}
\end{table}

\begin{table}[H]
\centering
  \begin{tabular}{|c|c|c|c|c|c|}
    \hline
    & $N_c=2$ & $N_c=3$ & $N_c=4$ & $N_c=5$ & Periodic \\
    \hline
    $\xi$ & 1.0116 & 1.02160 & 1.0299 & 1.0349 & 1.0420 \\
    \hline
  \end{tabular}
  \caption{Convergence to $\xi_{per}$, for the geometry given by City 2.
 $M=5$.}
  \label{tab:asymptot_for_city7}
\end{table}

\begin{table}[H] 
\centering
\includegraphics[scale=.5]{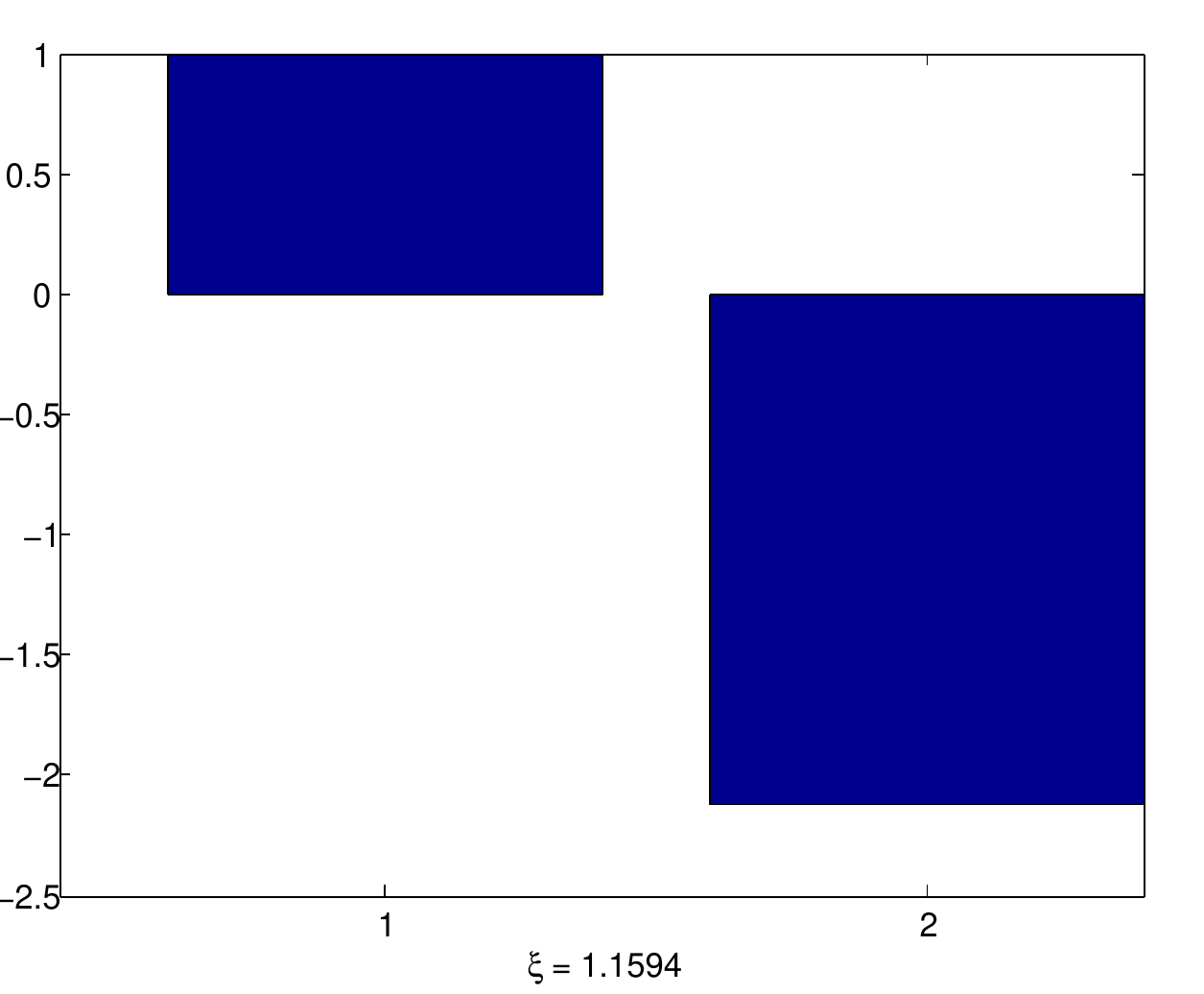}
\includegraphics[scale=.5]{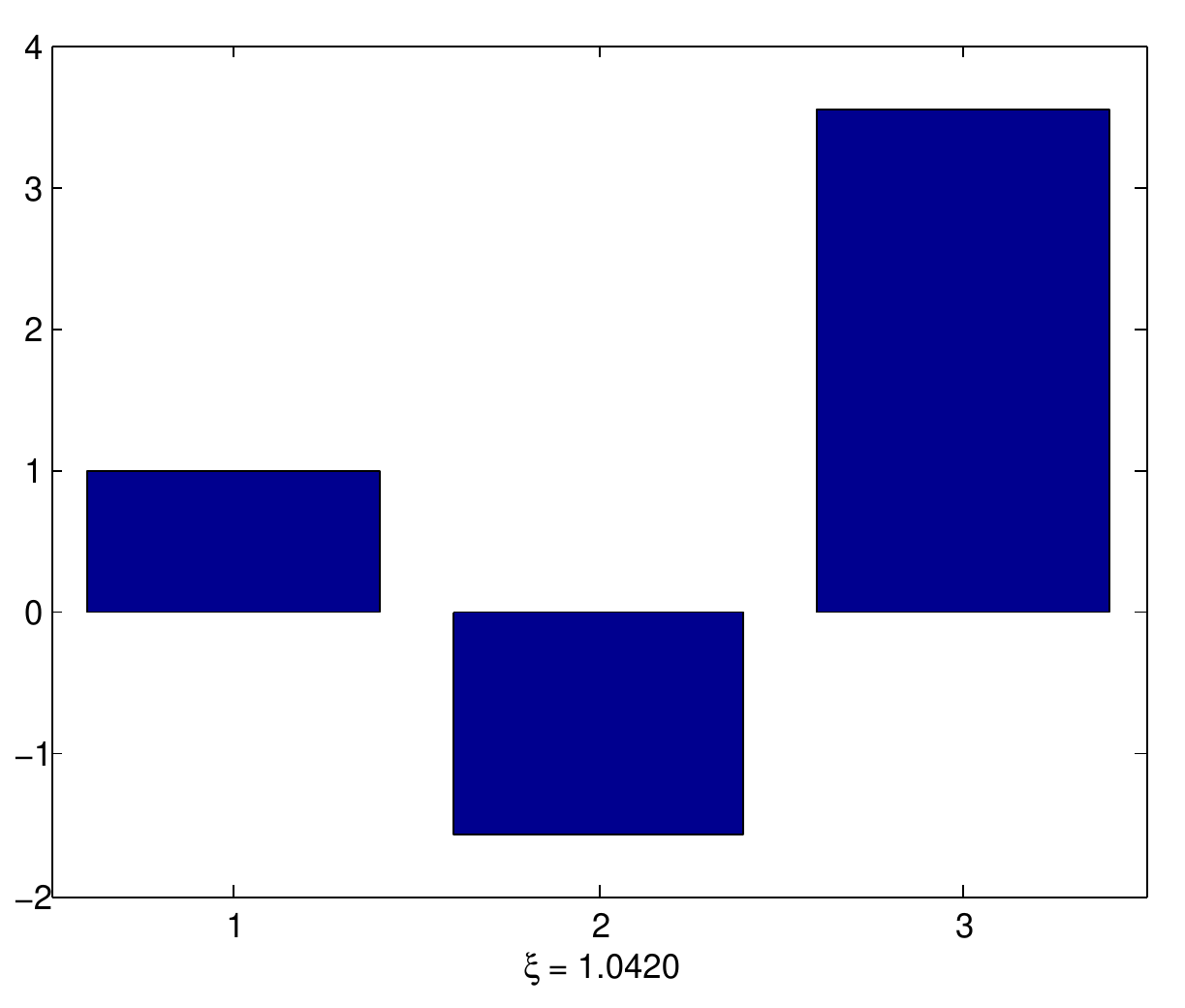} \\
\caption{Left: computed values of $(\alpha_1,\alpha_2)$ for the periodic case City 1-per.
Right:  computed values of $(\alpha_1,\alpha_2, \alpha_3)$ for the periodic case 
City 2-per.}
\label{periodic_outcome}
\end{table}
%----------------------------------------------------------------------------------------------------
\subsection{Conclusion and perspectives}
Using a model involving  vibrating tall buildings and harmonic elastic (anti plane) 
displacements of the underground, we have introduced in this paper  methods 
for computing   frequencies that achieve the coupling of these two vibrations.
Earlier computations by Ghergu  and Ionescu, \cite{ghergu}, were limited to the case
of a finite set of identical, equally spaced buildings. 
Their results are quite instructive since they point to a collective response of buildings  
to seismic waves, a phenomenon that authors have named "city-effect", \cite{urbanization_effect}.
We noted, however, that the  method for finding coupling frequencies introduced by Ghergu  and Ionescu
becomes quickly computationally expensive as the number of buildings grows large,
and may be hard to extend to a fully three dimensional setting.
Given that weakness, we have therefore resorted to the use of periodic 
formulations. We have   first obtained results that coincide with those of Ghergu et al.'s,
but that can be obtained at a much lower computational cost.
From there we have shown how our method can be extended to repeated sets of
non identical buildings: more realistic city lay outs can thus be modeled. 
Interestingly, in the case of non identical buildings, our simulations indicate that the response 
to this coupling phenomenon
may differ drastically from one building to another.\\
Our next endeavor, which will be the subject of a forthcoming publication, will be 
to study fully three dimensional models. 
Generalizing our work to that case will certainly prove to be quite challenging since
the 3D elastic half space Green's tensor is already very involved, and we would need
to compute its periodic analog.

\section{Appendix A: proof of symmetry for matrix $T$ defined by formula (\ref{operator_matrix})
}

Let $\psi_k$ solve (\ref{helmholtz_1}-\ref{sommerfeld_cond})
with $\alpha_j=1$, if $j=k$, and $\alpha_j=0$ otherwise.
It suffices to show that  
\bea \label{toshow}
\nonumber \int_{\Gamma_i} {\Psi_i \frac{\partial \Psi_j}{\partial y}}ds = \int_{\Gamma_j} { \Psi_j \frac{\partial \Psi_i}{\partial y}}ds.
\eea
Denote by  $D$ the set  $ \{(x,y): x^2+y^2 \leq r^2, y \geq 0 \}$. Due to Green's theorem, 
\begin{equation}
\nonumber
  \int_{D} \left( \Psi_i \Delta \Psi_j - \Psi_j \Delta \Psi_i \right) dx dy = \int_{\partial D} \left(\Psi_i \frac{\partial \Psi_j}{\partial \nu} - \Psi_j \frac{\partial \Psi_i}{\partial \nu} \right) ds,
\end{equation}
where $\nu$ is the unit normal vector to $\p D$ pointing outward.
Using (\ref{helmholtz_1}-\ref{sommerfeld_cond})
it follows that
\be
\nonumber
    \int_{\Gamma_i} {\Psi_i \frac{\partial \Psi_j}{\partial y}}ds  - \int_{\Gamma_j} {\Psi_j \frac{\partial \Psi_i} {\partial y}} ds
		= \int_{\{|x|=r, y>0\}} {\Psi_i \frac{\partial \Psi_j}{\partial r}}
		-  {\Psi_j \frac{\partial \Psi_i} {\partial r}} ds \\
 =   \int_{\{|x|=r, y>0\}} (\Psi_j - \Psi_i ) o(r^{-1/2})  ds
\ee
Due the decay of the fundamental solution of the Helmholtz equation
in a half plane, we have that $ \Psi_j = O(r^{-1/2})$ and 
$ \Psi_i = O(r^{-1/2})$, so if we let $r \ri \infty$ we arrive at identity (\ref{toshow}).

\section{Appendix B: the numerical solutions to equations (\ref{single_layer_potential2})
and (\ref{single_layer_potential2_per})}
\subsection{The free space case: equation (\ref{single_layer_potential2})}
It is well known in the literature that $\psi$ solution to (\ref{single_layer_potential2})
must have square root singularities at the edges of
$\Gamma_j$, therefore we set
\begin{equation} 
  \label{psi}
  \displaystyle \psi(s) = \frac{\phi(s)}{\sqrt{(s-a_j)(b_j-s)}}, \mbox{ where } s\in(a_j,b_j), \; 1\leq j \leq N, 
\end{equation}
and $\phi$ is a smooth function in $[a_j,b_j]$. 
Recalling (\ref{single_layer_potential}) we have
\begin{equation}
  \label{helmholtz_solution}
  \Psi(x,y) = \frac{i}{4} \sum_{j=1}^{N} \int_{a_j}^{b_j} {H_0^{(1)} (\xi \sqrt{(x-s)^2 + y^2}) \psi(s)} ds,
\end{equation}
We use potential theory to assert that 
$
  \lim_{t \rightarrow 0^{+}} \frac{\partial \Psi}{\partial y} (s,t) = - \frac{1}{2} \psi(s)$ for all $a_j < s < b_j$.
To employ the same numerical mesh for each building foundation $\Gamma_j$, we 
set for $t$ in $[-1,1]$
\begin{equation}
  \label{variable_t}
  \displaystyle s = g_j(t) = \frac{b_j-a_j}{2} t + \frac{b_j+a_j}{2}
\end{equation}
%In terms of $\phi$ and $t$ we rewrite operator $T$:
%\begin{equation}
%  \label{operator_T} \displaystyle \big( T(\xi^2) \alpha \big)_j = \frac{1}{2} \int_{-1}^{1} {\frac{\varphi \big( g_j(t)\big) } {\sqrt{1-t^2}} } dt.
%\end{equation}
%Using the first part of boundary conditions (\ref{helmholtz_2}), we get $N$ equations
Substituting in (\ref{helmholtz_solution}) we obtain
\begin{equation}
  \label{egn24}
  \displaystyle \frac{i}{4} \sum_{j=1}^{N} \int_{-1}^{1} {H_0^{(1)}(\xi |x-g_j (t)|) \frac{\varphi \big( g_j(t)\big) } {\sqrt{1-t^2}} } dt = \alpha_k, \; 1 \leq k \leq N, \; x \in \Gamma_k.
\end{equation}
We then solve for $\varphi \big( g_j(t)\big) $ following the numerical method
introduced in the appendix of \cite{ghergu}. Here we just recall that this numerical method relies on the
fact that (see \cite{abramovitz})
\be \label{splitup}
\f{i}{4} H_0^{(1)} (z) = A(z) \ln \f{z}{2} + B(z),
\ee
for any non zero complex number $z$, where $A$ and $B$ are two entire functions.

\subsection{The periodic case: equation (\ref{single_layer_potential2_per})}
%he analog of (\ref{bessel_decomp}) is
The periodic Green's function relative to problem (\ref{helmholtz_per_1}-\ref{per_cond})
can be written out as
\be
G_{per} (x,y) = \sum_{n=-\infty}^\infty  G(x- 2nP ,  y )
\ee
%% notation problem!!!!!!!!!!!!!!!!!!!
The analog of decomposition (\ref{splitup}) is now
\be
  \label{periodic_decomp}
 G_{per} (x,y) = A(\xi \sq{x^2 + y^2}) \ln \f{\xi \sq{x^2 + y^2}}{2} + \tilde{B}(x,y),
%  \ds \frac{i}{4} \sum_{n=-\infty}^{\infty}{H_0^{(1)}(\xi r_n)} = A_0(r_0) \ln{\frac{\xi r_0}{2}} + B_0(z_1, z_2),
\ee
%where $z_1=y$, $z_2 = x - s$, $\ds r_n = \sqrt{y^2 + (x-s-nd)^2}$, $\ds A_0(r_0) = -\frac{1}{2\pi} J_0(\xi r_0)$, and
%It is useful to factor out singularities as in (\ref{splitup}).
%Doing so we obtain, setting $z_n=\xi \sq( (x- 2nP)^2 + y^2)$
%\be \label{splitup}
%G_{per} (x,y) = A(z_0) \ln \f{z_0}{2} + B(z),
%\ee
where 
\be
  \label{periodic_func_B}
  \tilde{B}(x,y) = \begin{cases}
           \ds \sum_{n=-\infty}^\infty  G(x- 2nP ,  y )
           -  A(\xi \sq{x^2 + y^2}) \ln \f{\xi \sq{x^2 + y^2}}{2} , & \mbox{ if } x^2 + y^2 \neq 0, \\ \\
           \ds \frac{i\pi - 2\cal{C}}{4\pi} + \sum_{n \neq 0} G(- 2nP ,  0 ), & \mbox{ otherwise }
         \end{cases} 
\ee
and $\cal{C}$ is the Euler constant. 
Note that $\tilde{B}$ is real analytic in $(x,y)$: this is due to the asymptotics of the Hankel function
$H_0^{(1)}$, see \cite{abramovitz}.
At this stage we see that the numerical method for the integral equation for the free space case
(\ref{single_layer_potential2}) can be extended to  the integral equation (\ref{single_layer_potential2_per})
 for the periodic case: all we need to do is to replace $B(\xi |x-g_j(t)|)$
by $\tilde{B}(\xi (x-g_j(t)),0)$.
It is crucial to be able to compute $\tilde{B}(0,0)$ for this numerical method to be applicable:
see \cite{darko} for a more detailed account of a comparable calculation and computational method
for a related Green's function.
\\
As explained above, we need to be able to  efficiently compute the two slowly convergent sums
$\ds \sum_{n=-\infty}^\infty  G(x- 2nP ,  y )$  and 
$\ds \sum_{n \neq 0} G(- 2nP ,  0 )$. The first sum can be efficiently computed
by Ewald's method. For that subject, we refer the reader to 
\cite{linton}.
We did not find in the literature any results on the computation of the second sum, however, it can be inferred
from the first sum. Here is how:
setting in this appendix only, $p=\pi/P$, $r_m=\sq{(x- 2mP)^2 + y^2}$ and
\bea
\gamma_m=\sq{m^2 p^2 - \xi^2}, \mbox{ if } m^2 p^2 - \xi^2>0, \\
\gamma_m= i \sq{-m^2 p^2 + \xi^2}, \mbox{ if } m^2 p^2 - \xi^2<0 .\\
\eea
Then applying Ewald's formula (see \cite{linton}),
\begin{equation}
  \label{periodic_Green_func_Ewald}
  \begin{array}{c}
    \ds G_{per}(x,y) = \frac{1}{8P} \sum_{m=-\infty}^{\infty} {\frac{e^{ipmx}} {\gamma_m}} \left[ e^{\gamma_m y}
		\erfc \left( \frac{\gamma_m P}{a} + \frac{ay}{2P}\right) + e^{-\gamma_m y}
		\erfc \left(\frac{\gamma_m P}{a} - \frac{ay}{2P}\right) \right] \\
  \ds + \frac{1}{4\pi} \sum_{m=-\infty}^{\infty} \sum_{n=0}^{\infty} {\frac{1}{n!} 
	\left( \frac{\xi P }{a} \right)^{2n} E_{n+1} \left( \frac{a^{2}r_{m}^{2}} {4 P^2} \right)},
  \end{array}
\end{equation} 
where $a >0$ is called ``splitting parameter'',
\begin{equation}
  \label{error_func}
  \erfc(z) = \ds \frac{2}{\sqrt{\pi}} \int_{z}^{\infty}e^{-t^2} dt
\end{equation}
is the complementary error function, and 
\begin{equation}
  \label{E_integral}
  E_{n}(z) = \ds \int_{1}^{\infty}t^{-n}e^{-zt}dt
\end{equation}
is the exponential integral. Note that the present method is valid only if $\gamma_m \neq 0$ for all
integers  $m$. (Note that if $\gamma_m = 0$ for some $m$ then the system of equations
(\ref{helmholtz_per_1}-\ref{per_cond}) is non uniquely solvable). \\
We now proceed  to find a new series formula for the
expression $\ds \sum_{m \neq 0}{\frac{i}{4} H_{0}^{(1)}(\xi r_{m})}$ for 
 $(x,y) =(0,0)$. %(\ref{periodic_Green_func_Ewald}).
 The only singular term in (\ref{periodic_Green_func_Ewald})
as $x=y=0$  is $E_{1}\left(\frac{a^{2}r_{0}^{2}}{4 P^2}\right)$:
it appears in the second sum for $m=n=0$.  Note that  
\begin{equation}
  \label{singular_part}
  \ds E_{1}\left(\frac{a^{2}r_{0}^{2}}{4P^2}\right) = - \left\{ {\cal{C}}
	+ \ln \frac{a^{2}r_{0}^{2}}{4P^2}  + \sum_{k=1}^{\infty}{\frac{(-1)^k}{k \cdot k!} \left(\frac{a^{2}r_{0}^{2}}{4P^2} \right)^{k}} \right\},
\end{equation}
according to  formula 5.1.11 in \cite{abramovitz}. 
We may also expand $\ds \frac{i}{4} H_{0}^{(1)}(\xi r_0)$
as follows
\begin{equation}
  \label{H_0_series}
  \begin{split}
    \ds \frac{i}{4} H_{0}^{(1)}(\xi r_0) = & \frac{i}{4} \big( J_0(\xi r_0) + iY(\xi r_0) \big) \\
    = & \frac{i}{4} J_{0}(\xi r_{0}) -   \frac{1}{2\pi} (\ln\frac{\xi r_0}{2} + 
		{\cal{C}})J_{0}(\xi r_{0}) + \frac{1}{2\pi} \sum_{m=1}^{\infty}{a_{m} \frac{(-1)^m}{(m!)^2} \left( \frac{\xi r_0}{2} \right)^{2m} } , 
  \end{split}
\end{equation}
where $ \ds a_{m} = \sum_{j=1}^{m}{\frac{1}{j}}$ (see formulas  9.1.13 in \cite{abramovitz}). \\
Finally, we substitute $E_1$ in (\ref{periodic_Green_func_Ewald}) by (\ref{singular_part}), and we subtract (\ref{H_0_series}). 
We note that the logarithmic singularities in $\ln r_0$ cancel out.
As $x=y=0$, $r_m =2|m|P$ and we obtain, 
\be
\ds \sum_{n \neq 0} G(- 2nP ,  0 ) = 
\f{1}{4P} \sum_{m=-\infty}^{\infty} \f{1}{\gamma_m} \erfc (\f{\gamma_m P}{a})
+ \f{1}{4 \pi} \sum_{m \neq 0}\sum_{n=0}^\infty \f{1}{n!} (\f{\xi P}{a})^{2n} E_{n+1} (a^2 m^2) \nonumber\\
+ \f{1}{4 \pi} \sum_{n=1}^\infty \f{1}{ n \,n!} (\f{\xi P}{a})^{2n} 
+ \f{1}{2 \pi} \ln \f{\xi P}{a} \nonumber
\ee

\end{document}